\documentclass[12pt]{amsart}
\usepackage{amsthm,amsmath,amssymb,amscd}
\usepackage{amssymb}
\usepackage{graphics}
\usepackage{epsfig}

{
   \newtheorem{theorem}{Theorem}[subsection]                     
   \newtheorem{proposition}[theorem]{Proposition}     
   \newtheorem{lemma}[theorem]{Lemma}
   \newtheorem*{claim}{Claim}
   \newtheorem{corollary}[theorem]{Corollary}

}
{\theoremstyle{definition}

   \newtheorem{definition}[theorem]{Definition}
   \newtheorem{remark}[theorem]{Remark}
}
\newcommand{\RR}{{\mathbb{R}}}

\newcommand{\CC}{{\mathbb{C}}}
\newcommand{\QQ}{{\mathbb{Q}}}

\newcommand{\PP}{{\mathbb{P}}}
\newcommand{\ZZ}{{\mathbb{Z}}}

\newcommand{\bbA}{{\mathbb{A}}}
\newcommand{\OO}{{\mathbb{O}}}

\newcommand{\bM}{{\overline {\mathbf{M}}}}
\newcommand{\bfM}{{{\mathbf{M}}}}
\newcommand{\bmu}{{\boldsymbol{\mu}}}
\newcommand{\cA}{{\mathcal A}}

\newcommand{\cC}{{\mathcal C}}
\newcommand{\cD}{{\mathcal D}}
\newcommand{\cE}{{\mathcal E}}

\newcommand{\cI}{{\mathcal I}}

\newcommand{\cL}{{\mathcal L}}
\newcommand{\cM}{{\mathcal M}}
\newcommand{\cN}{{\mathcal N}}
\newcommand{\cO}{{\mathcal O}}

\newcommand{\cQ}{{\mathcal Q}}

\newcommand{\cS}{{\mathcal S}}
\newcommand{\cU}{{\mathcal U}}
\newcommand{\cX}{{\mathcal X}}
\newcommand{\cY}{{\mathcal Y}}
\newcommand{\cZ}{{\mathcal Z}}
\newcommand{\cV}{{\mathcal V}}

\newcommand{\down}{\downarrow}
\newcommand{\Spec}{\operatorname{Spec}}

\newcommand{\dar}{\downarrow}

\newcommand{\isom}{\operatorname{Isom }}

\newcommand{\double}{\genfrac..{0pt}1
{\raise -1pt\hbox{$\scriptstyle\to$}}{\raise 3pt\hbox {$\scriptstyle\to$}}}

\begin{document}

\title[ Explicit stable models of elliptic surfaces with sections] 
{Explicit stable models of elliptic surfaces with sections}
\author[Gabriele La Nave] {Gabriele La Nave} 
\address{Courant Institute of Mathematical Sciences\\ New York University\\ 
251 Mercer Street\\ New York, NY 10012\\ USA} 
\email{lanave@cims.nyu.edu}

\maketitle

\vspace{-,15in}

\begin{abstract}

In this note we propose to show how to find stable models of
a one-parameter family of elliptic surfaces.  
The strategy we use was initiated by
 Abramovich-Vistoli in \cite{A-V}: that is to say to consider a fibred 
surface as a map from the base curve to the moduli 
of stable n-pointed curves, and to consider then the Kontsevitch space of such 
maps that are stable. In our case we will then be describing the moduli 
space of stable surfaces via the moduli space of Kontsevitch stable maps 
to $ \bM _{1,1},$ the moduli space  of stable 
one-pointed elliptic curves. Then, following Koll\'ar--Shepherd-Barron
(\cite{K-S}) and V. Alexeev (\cite{A1}) we apply Mori's Minimal Model
Program in an explicit manner by means of toric geometry.
 
\end{abstract}

\section{Introduction}
\subsection{Description of the problem}

When speaking of a {\bf surface} we will always mean a reduced,
integral, normal, projective noetherian scheme of dimension $2.$ 

For the purpose of studying the boundary of the moduli space, we need
to introduce the following generalization of the well-known concept of
{\bf elliptic surface} (for convinience, we will mantain the name):  

\begin{definition}
 An {\bf elliptic surface with zero section} is the datum of
 a surface $X$ 
together with a proper map, 
$\pi : X \to C$ to a proper curve $C$ and a section
 $\sigma : C \to X,$  
called the {\bf zero section}, such that:
\begin{enumerate}

\item  the 
generic 
fibre of $\pi$ is a stable complete curve of arithmetic genus $1;$

\item the zero section is not contained in the singular locus of $\pi.$ 
\end{enumerate}

Such an object is called {\bf relatively minimal} if it is smooth and 
there is no 
$(-1)\text{-curve}$ 
in any fibre. Furthermore, we say it is {\bf minimal} if it is smooth 
and contains no 
$(-1)\text{-curve}$ at all. Notice that if the base curve is not rational, 
the last 
two notions coincide.
 
\end{definition} 

\begin{remark}
Note that our notion of elliptic surface differs from the usual one in that 
in point $1$ we ask the generic fibre to be stable, as opposed to the usual 
notion where the generic fibre is required to be smooth, and in that we do
not ask the map $\pi$ to be flat.
\end{remark}

Whenever we are given such an object, we automatically get a rational map 
$C \dashrightarrow \bM _{1,1}$ to 
the 
moduli space of $1$-pointed elliptic curves: $  {\bM} _{1,1}$. If the base 
curve were smooth, we would then extend the map to a regular map on the whole
 $C\to \bM _{1,1}$.   
 Composing this map with $\pi$ provides us with a regular map 
$X\to C \to {\bM}_{1,1} .$

Let $S$ be a fixed base scheme.
Let $\pi: X\to C \to S$ be a proper $S$-map of relative dimension $2$
from a proper scheme $X$ to a proper scheme $C,$ with $X \to C$
generically fibred in stable curves of genus $1.$

 Following Alexeev and Koll\'ar--Shepherd-Barron (cf. \cite {A1}, \cite {A2}
and \cite {K-S}), 
we are naturally led to the following: 
 
\begin{definition}\label{logstabsurf}

A {\it pair } $(\pi :\cX \to \cC\to S,\cQ)$ consisting of an $S$-morphism  
$X \to C$ fibred in generically
stable curves of genus $1$ 
and a {\em section } $\cQ$ is called {\bf stable } if, for each 
geometric point $s$ in $S$:

\begin{enumerate}
\item $\cX\to S$ and $\cC \to S$ are flat morphisms of relative dimension 
$2$ and $1$ respectively;
\item the pair $(\cX _s , \cQ _s )$ has semi-log-canonical 
singularities (see section \ref{prelim} 
for a definition);
\item the relative log-canonical sheaf $ \omega _{X/S} (Q)$ is $\QQ$-Cartier;
\item the relative log-canonical sheaf $ \omega _{X/S} (Q)$ is $S$-ample.

\end{enumerate}

\noindent
If $S=Spec(k)$ for any field $k, $ we simply say that the pair is {\bf stable}
  
\noindent
Similarly, a {\it triple } $(\cX \to \cC \to S,\cQ, 
f: \cX \stackrel{ \pi}{\to}\cC \stackrel {j}{\to} {\bM}_{1,1})$ consisting of 
an $S$ morphism $\cX \to \cC$ fibred in generically smooth 
elliptic curves, 
a {\em section } $\cQ$ and a map $f$ as above
will be called {\bf stable} if conditions ~1, ~2 and ~3 above hold and if 
  $ \omega _{X/S} (Q)
\otimes f^* \cO (3)$ is {\it ample}.

\noindent
If $S=Spec(k)$ for any field $k, $ we simply say that the triple is {\bf stable}
\end{definition}

Given a relatively minimal elliptic surface $X$ with section $Q,$ it 
is natural to consider its associated Weierstrass model $Y$: this is roughly 
obtained by contracting all the components of the fibres that do not meet 
$Q$ (see section \ref{wform} for a more precise definition). This surface 
$Y$ has possibly a finite number of cuspidal fibres, with the property 
that near each such fibre $Y_0,$ the surface $Y$ has a local equation of the 
form $y^2 = x^3 + a x +b$ with $a,b \in \cO _{C,0}$ 
with $min (\nu _0(a^3),\nu _0 (b^2)) <12 .$ Here $\nu _ 0 (g )$ 
denotes the order of vanishing of $g \in \cO _{C,0}$ at $0.$ 
 Such equations are called {\bf minimal Weierstrass equations} and 
$Y$ is said to be in {\bf minimal Weiertrass form}.

 Loosely speaking Weierstrass models, within the theory of elliptic surfaces,
 play the role that canonical models do for surfaces of general type. Indeed 
if the base curve $C$ is smooth and  non-rational, they 
constitute the log-canonical models for elliptic surfaces with sections having 
$C$ as base curve (see Corollary \ref{semiampl}).
  
For definition \ref {logstabsurf} to be of any use, we want the moduli 
functor it defines 
to include at least a large part of the locus of minimal 
Weierstrass equations.
In other words, we want most pairs
 $(\pi:X \to C,Q)$ consisting of an ellitpic surface in {\it minimal 
Weierstrass equation} $X\to C$ with section $Q,$ to be stable 
according definition \ref{logstabsurf}.
 Similarly, we want most triples  $(\pi:S \to C,Q, X\stackrel{\pi}{\to} C 
\stackrel{j}{\to} \bM _{1,1})$  to be stable according to definition 
\ref{logstabsurf}.
 In section \ref{wform} we will show that 
a pair $(X\to C, Q)$ with $X\to C$ in minimal Weierstrass form is {\it
  stable } if $C$ is not rational, and 
that a triple $(\pi:S \to C,Q, X\stackrel{\pi}{\to} C 
\stackrel{j}{\to} \bM _{1,1})$ with $X\to C$ in minimal Weierstrass
form is {\it stable } either if $C$ is not rational, or, if 
it is rational, if 
$\pi: X \to \PP ^1$ is not an 
isotrivial family (i.e., mapped to a point in $\bM _{1,1}$ via $j$).

This sets the course for us: the moduli problems we will be looking at are  
the ones 
defined by the functors:

$$\begin{array}{cccc} \cM _{pairs}: &{\mathcal Sch} &\longrightarrow
           & {\mathcal Sets} \\
                
                & S & \to & \cM_{pairs}(S)
\end{array}$$
where:
$$
 \cM _{pairs}(S):= \left\{
\begin{array}{ccc}  
 &\text{ pairs over } S \ (\pi:{\mathcal X}\to {\mathcal C},
\mathcal Q,): \text{ which }\\ 
&\text{ are stable according to definition \ref{logstabsurf} }
  
 \end{array}  
 \right\}  /  {\text{isomorphisms} }
$$
and:

$$\begin{array}{cccc} \cM _{triples}: &{\mathcal Sch} &\longrightarrow
           & {\mathcal Sets} \\
                
                &S & \to & \cM_{triples}(S)
\end{array}$$
where:
$$
 \cM _{triples}(S):= \left\{
\begin{array}{ccc}  
 &\text{ triples over } S \ (\pi:{\mathcal X}\to {\mathcal C},
\mathcal Q,
 f: {\mathcal X}\stackrel {\pi}{\to} \cC \stackrel {j}{\to}
 { {\bM}}_{1,1}):\\ 
&\text{ which are stable according to definition \ref{logstabsurf} }
  
 \end{array}  
 \right\} /{\text{isomorphisms} }
$$
To render our functors of finite type, it is essential to fix 
some numerical invariants.
 The pair of rational numbers:
$$A =  c_1(\omega_X(Q))^2;\quad \chi = \chi (\cO _X)$$

for the moduli of pairs and the triple of rational numbers: 
$$A =  c_1(\omega_X(Q))^2;\quad B= c_1(\omega_X(D))\cdot 
c_1(f^*\cO _{\bM _{1,1}}(1));
\quad
C = c_1(f^*\cO _{\bM _{1,1}}(1)))^2,$$
for the moduli of triples.
 Once a moduli problem is defined, there are two crucial questions one
 must pose and answer (preferably positively): 1)Is the associated functor proper?
 Is it a Deligne-Mumford stack?

Here we will be addressing only the question of properness. Specifically 
we will prove that the functors $\cM _{pairs}$ and  $\cM _{triples}$
satisfy the valutative criterion of properness.
Such a criterion is usually called ``stable reduction theorem'' for 
moduli problems.  In our case one has to prove that in 
a family of stable pairs (resp. triples) over the punctured disk, one 
can replace the central fibre of any compactification over the whole 
disk by a stable one, possibly after a base change. 

\subsection{Previous work}
It must be said that the questions we are addressing here have been addressed 
and answered in much greater generality by J. Koll\'ar and N. Shephard-Barron 
(cf. 
\cite {K-S}) and by V. Alexeev  (cf. \cite {A1}, \cite {A2}).
 In \cite {K-S} J. Koll\'ar and N. Shephard-Barron prove a stable 
reduction theorem for the moduli of
pairs $(X,D)$ consisting of a {\it semi-log-canonical } surface $X$ and a 
$\QQ$-Cartier divisor $D.$ 
In \cite {A1} and \cite {A2}, V. Alexeev generalizes this to moduli 
of triples $(X,D, f:X \to M)$ where $X$ is a surface with at most 
semi-log-canonical 
singularities, 
$D \subset X$ is a $\QQ \text{-Cartier}$ divisor, $f:X \to M$ is a proper 
morphism to a projective scheme $M$ and 
$\omega _{X} (Q) \otimes f^* A$ is $\QQ$-Cartier and ample for every choice of 
a sufficiently ample divisor $A$ on $M.$ Moreover he proves that the corresponding 
moduli functor is bounded (in particular for $M=\Spec (k)$ a point, one gets 
that the moduli space of pairs is bounded). 
 Their proofs make use of the full strenght of Mori's MMP 
(Minimal Model Program).

Also, in the case of Weiestrass fibrations over $\PP ^1,$ in 
\cite{M}, R. Miranda constructs
 a proper moduli space by identifying the GIT semi-stable points of the 
action of $k ^* \times SL(V_1 )$ on a suitable subset $T_N \subset V_{4N} 
\oplus V_{6N}.$ Here $V_1 = H^0(\PP ^1, \cO _{\PP ^1} (1)),$ 
$V_k = H^0(\PP ^1, \cO _{\PP ^1} (k))=Sym ^k (V_1 ),$ and the component 
$k^*$ of $k ^* \times SL(V_1 )$ acts on $T_N$ by 
$\lambda (A,B) = (\lambda ^{4N} ,
\lambda ^{6N})$ and $SL(V_1 )$ acts on $T_N$ by 
the action induced on it by the 
natural action of $SL(V_1 )$ on $Sym ^k(V_1 ).$

What we set out to do here, though, is to give an
 explicit 
description of the stable reduction process 
(in the sense of Alexeev--Kollar--Shephard-Barron, ie., in the sense 
of definition \ref{logstabsurf}) in ~1-parameter families 
and of the possible surfaces we may
 get 
at the boundary of the moduli spaces of elliptic pairs and of elliptic triples
respectively. In fact, the work of Kollar--Shepherd-Barron and of Alexeev  proves 
the valuative criterion for properness abstractly by means of the MMP. 

In order to explain better the significance of such an explicit
description and how it compares with the construction of
Kollar--Shepherd-Barron and of Alexeev, we may
refer to the case of the space $\bM _{g,n}$ of {\bf Deligne- Mumford}
({\bf DM} for short)
{\bf stable curves} of a given genus g with sections $\sigma _i$ with
$i=1,...,n,$ which compactifies the moduli space of {\it smooth}
curves of genus $g$ with $n$ marked points $\bfM _{g,n}$ (cf. \cite
{D-M} or \cite {H-M}). A (geometrically connected and proper) curve $\cC \to S$ over a scheme
$S$(resp. with sections $\sigma _i : S \to \cC$) can be defined to be
a D-M
stable curve in two
different and equivalent ways: 

{\bf the abstract description:} \begin{enumerate}

\item  the singularities of $\cC _s$ are at worst nodal for every
  closed point $s\in S;$
\item the relative dualizing sheaf $\omega _{\cC/S}$ (resp. $\omega
  _{\cC/S}(\sum _{i=1} ^n \sigma _i)$ ) is ample.

\end{enumerate}

{\bf the combinatorial description:} \begin{enumerate}

\item as $1$ above, 

\item if there is a rational irreducible component $R$ of $\cC _s$
  (for some closed point $s\in S$) then $R$ must meet the rest of $\cC
  _s$ in at least $3$ points (resp. $R$ must contain at least three
  points that are either nodes of $\cC _s$ or {\it marked points}
  coming from the sections $\sigma _i$); if there is a component $E$ of
  arithmetic genus $1,$ then $E$ must meet the rest of $\cC _S$ in at
  least one point. 
\end{enumerate}

An analogous picture holds in the case of {\it Kontsevich stable maps}.

 In this work we try to transpose the work of J. Koll\'ar and
 N. Shepard-Barron and of V. Alexeev into an incarnation that would
 correspond to the combinatorial picture of the DM-stable curves given
 above.

Up to now there are very few cases in which the degenerations 
in a ~1-parameter 
family have been described explicitly. To name one of these,
recently B. Hasset \cite{Ha} 
gave a description for the boundary of the moduli of pairs 
$(\PP ^2, C)$ where $C$ is a plane quartic, and built an isomorphism between
this space and $\bM _3,$ the moduli space of 
Deligne-Mumford stable genus $3$ curves (this is of importance also
for the undersdanding of the locus
of {\it limiting plane curves} in $\bM _ g,$ given that this locus
coincides with the whole $\bM _3$ for degree $4$ curves).

The present work adds to the list of cases that have been worked out.

\subsection{The strategy}

The strategy for the moduli of triples is the one initiated by 
Abramovich 
and Vistoli in \cite{A-V}: we consider an elliptic surface as a map from 
the base curve 
to $\bM _{1,1}.$ By marking the points on the base curve corresponding to the 
cuspidal fibres, thanks to the {\bf Purity Lemma} of Abramovich and Vistoli 
(cf. \cite {A-V}), 
we can replace these fibres by finite cyclic quotients of 
stable curves (we will refer to such curves as {\bf twisted curves}).
The idea is to now make the map from the base curve to $\bM _{1,1}$ 
Kontsevich stable. 

The result is that at the limit there are surfaces 
$X$ that map to curves $C;$ these curves come endowed with a
Kontsevich stable maps $C \to \bM _{1,1},$ wich correspond to the map 
to moduli; the general fibre of $X \to C$ is a stable curve and 
the componts of $X$ meet along fibres that are either stable or twisted.
The problem is that one had to replace the cuspidal curves in 
our original elliptic surfaces: the new surfaces are not in Weierstrass form 
anymore.

 In order to deal with this problem, we prove an extension lemma 
(see lemma \ref{Extension Lemma}), 
that allowes us to place back the cuspidal 
fibres, and so we can remove the extra marked points of the base. 
 
The double curves of $X$ are either stable or the {\it twisted fibres} 
of Abramovich--Vistoli (cf \cite{A-V2}).
 Looking at the irreducible components, we are therefore led to 
enlarging the class of {\it minimal Weierstrass surfaces} to what we call 
{\bf quasiminimal elliptic surfaces}. Roughly speaking, these are elliptic 
surfaces $X\to C$ over smooth curves $C,$ which are in Weierstrass form 
away from a finite number of {\it twisted fibres} and such that the local 
Weierstrass equation away from the twisted fibres is minimal.

 The price we 
have to pay in removing the extra marked points, is that the map from the 
base curve to moduli might not be Kontsevich stable anymore: there may
be isotrivial components mapping to a rational component of 
the base curve $C$, 
that meet the rest of $C$  
in one or two points. It turns out (proposition \ref{extr ray})
that the componets of $X$ dominating these isotrivial rational components 
are unstable. To be precise, the zero section of these components is 
an {\it extremal ray } if the base curve meets the rest of $C$ in only 
one point, 
and it is contracted by the log-canonical map if it meets the rest of $C$ in 
two points. 

To deal with these components, we start by performing a few explicit 
birational transformations by means of local Weierstrass 
equations; in doing so we are led to enlarging the class of {\it quasiminimal 
elliptic surfaces} to what we call {\bf standard elliptic surfaces}, 
which are, roughly speaking, elliptic surfaces with a finite number 
of twisted fibres,  
for which the local Weierstrass equations away from the twisted fibres 
are of the form $y^2 = x^3 + a x +b,$ for $a,b \in \cO _{C,p},$ 
with the minimum of the order of vanishing of $a^3$ and $b^2$ not greater than 
$12$ at each point $p \in C.$ Note that if the above mensioned minimum does 
achieve $12$ at some point $p \in C,$ the surface has an 
{\it elliptic singularity} at the point $x=y=0$ over $p$ (see section \ref{LCWE}).

We continue by following the steps of the Minimal Model Program 
(MMP). By means of toric geometry we are able to explicitly perform the necessary 
log-flips and small log-contractions. 

 The forgetfull functor sending, for each 
scheme $S,$ an $S$-triple $(\cX, \cQ,f:\cX \to \cC \to \bM _{1,1})$ to the 
$S$ pair $(\cX, \cQ),$ is not well-defined at the level of moduli: 
it does not always produce a stable pair out of a stable triple. 
 In fact, every component $X\to \PP^1$ of a geometric fibre 
$\cX _s \to \cC _s ,$ for some geometric point $s$ in $S,$  that meets 
$\cX _s$ in less than three fibres is {\it unstable}. Therefore, we need to 
perform more steps of the MMP. It will turn out that the same explicit 
operations described in the case of triples, go through for these more general 
settings in which the $j$-invariant associated to such components 
$X \to \PP ^1$ is not constant. The price we have to pay is that 
the map $\cX \to \bM _{1,1} $ is no longer regular.

When the general base curve $\cC _{\eta}$ is isomorphic to $\PP ^1 ,$ 
the log-canonical bundle for the pair $(\cX, \cQ)$ is not nef: it is 
negative on $\cQ _{\eta}.$ We therefore need to contract the zero-section 
in the central fibre, to obtain a {\it stable pair}. In order to do so we 
need to perform some extra birational operations on the total space of 
the family $\cX .$
   
\subsection{The result}

We need some definitions first.

Loosely speaking a {\bf log-standard} elliptic surface $X\to C$ (see 
definition \ref{logstand}) is a certain explicit locally toric blow-up 
(the unique such {\it semi-logcanonical} blow-up) of a 
pair $(Y\to C, G+F+Q)$ consisting of a standard elliptic surface 
$Y\to C$ together with the marking of the zero section $Q$ and of some 
(possibly 
twisted) fibres $G= \bigcup _i G_i$ and $F= \bigcup _i F_i;$ 
the centers of the blow-ups are supported on points in $G_i \cap Q . $ 
 We call {\bf splice} the proper transform of a fibre $G_i$ of $Y \to C$ 
on which the point we blow-up is supported, and {\bf scion } the proper
transform of each of the $F_i $'s. Such a {\it log-standard elliptic
  surface} is called {\bf strictly stable}, if $Q ^2 < -1 .$

Our results are stated in Theorems \ref{stabredpairs}, \ref{stabred}, 
\ref{rationalbase} and \ref{ellipticbase}.

Theorems \ref{stabredpairs}  and \ref{stabred} roughly speaking say 
that, in the case in which the base curve has genus $g\geq 2,$ at the
boundary we get a union of surfaces 
$X =\cup X_i$
 attached one to the other along curves which are either stable or
 twisted fibers or along splices, mapping to a curve 
$\pi: X \to C ,$ with a section $Q .$ 
For the {\bf moduli of pairs}, $C$ is a Deligne-Mumford stable curve and
 the components of $X$ are (possibly) of two 
kinds:
\begin{enumerate}
\item  strictly stable log-standard components mappping dominantly onto irreducible components 
of $C;$
\item components mapped to a point of $C$ via $\pi .$
\end{enumerate}

 The components mapped to a point, which we call {\bf pseudoelliptic}
 (see definitions \ref{I} and \ref{II}), 
are further devided in two types. For the 
moduli of {\it pairs} $(X,Q),$ we have: 
\begin{enumerate}

\item {\bf type I:} A pair $(Y, G+F)$ consisting of a surface $Y$
  endowed with a {\it structure morphism} $f:Y'\to Y ,$ which is regular and birational, from a
   {\it log-standard} surface $(Y' \to \PP ^1, Q+ G' + F'), $ with $1$ {\it
    scion} $F'$ and a number of {\it splices} $G' = \bigcup _i G_i '$
  and such that $F= f(F')$, $G_i = f({G _i}')$ and $G= \bigcup G_i.$ 
  The surface $Y$ is attached to the rest of $X$ along $F$ and
  possibly along some (or all) of the $G_i$'s. Furthermore, the
  morphism $Y' \to Y$ is obtained by
 torically (and explicitly) blowing-down the zero section $Q$ of $Y' ;$

\item {\bf type II:} A pair $(Y ,G + F_1  +F_2 ) $
  consisting of a surface $Y$ endowed with a birational morphism $f:
  Y' \to Y$ (the {\it structure morphism})
  from a {\it log-standard} surface $(Y', Q+ G' +{F _1 } '+ {F_2 }')$
  with $2$ {\it scions} ${F _1 } ' $ and $ {F_2 }' $ and a number of
  {\it splices } $G' = \bigcup _i {G_i} ';$ also $ F_i = f({F _i}')$,
  $G_i = f({G _i}')$ and $G= \bigcup G_i.$  
 Moreover $Y$ is attached to the rest of $X$ along the $F_i$'s and possibly along some (or all) of the $G_i$'s. 
Furthermore, the
  morphism $Y' \to Y$ is obtained by
 torically (and explicitly) blowing-down the zero section $Q$ of $Y' ;$

\end{enumerate}
We still call {\it splice } and {\it scion} the images in a {\it
  pseudoellitpic surface}, either of type I or of type II,  of the
corresponding curves via  the structure morphism.   
The attaching of a scion to a splice is
\'etale locally described by the fan given in Theorem \ref{log-flip}.

The attaching of a {\it pseudoellitpic surface 
of type II} along a marked fibre of $X$ is described 
\'etale locally by the cone given in Theorem \ref{logcancont}.

For the {\bf moduli of triples} $(X,Q, f:X \to \bM _{1,1})$
 we have the same situation except for having the {\it isotrivial}
analogues of {\bf type I} and {\bf type II}, and we ask that the map 
$C\to \bM _{1,1}$ be Kontsevich stable.

The {\bf type I} surfaces arise from log-flips, and the {\bf type II} ones 
from the log-canonical contraction of the zero section.

\begin{picture}(350,320)(0,0)
\put (180,50){\makebox(90,200) {\includegraphics{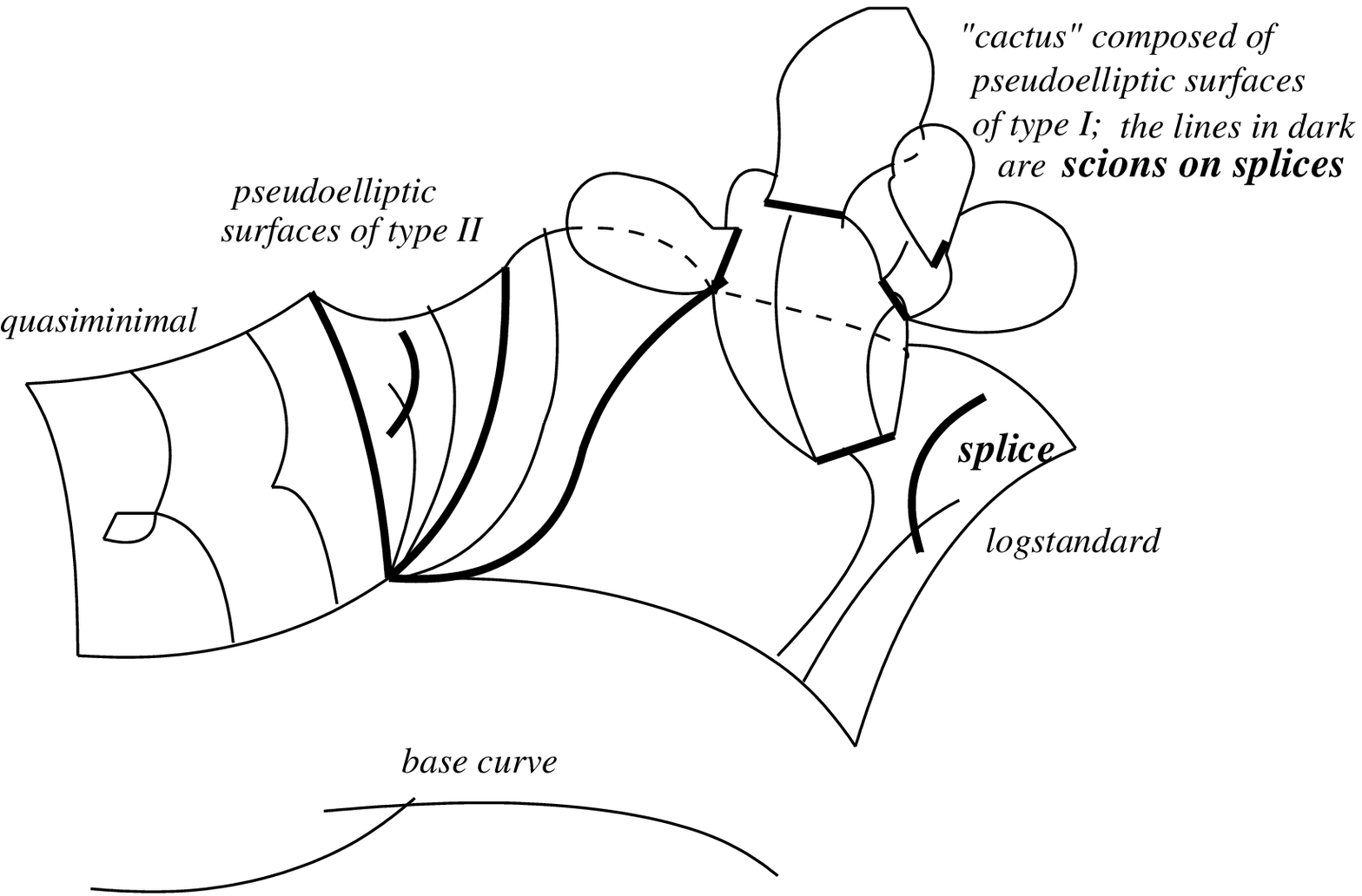}}}
\put(200,10){{fig.1}}
\end{picture}

\vskip 30pt

In the rational base case (see Theorem \ref{rationalbase}), 
there is no base curve at all. At the boundary here we get a union 
$X= \cup X_i$  where $X_i$ is a  
{\bf pseudoelliptic surface of type N} with $N=0$ or $N=1.$ 
 Furthermore the type $1$ pseudoelliptic surfaces come in pair, and
 they are attached one to the other (but not along as a scion).

Here, by {\bf pseudoelliptic surface of type N} we mean, loosely
speaking, a surface $S$
and a map $g: S' \to S$ from a log-standard surface $(S\to \PP ^1,
Q+\sum _{i=1} ^N F_i +\sum _{i=1} ^n G_i)$ with {\bf N } marked stable
or twisted fibres
$F_i$ and $n$ {\it splices} $G_i, $ furthermore, $g$ is obtained by an
explicit toric blow-down of the zero section $Q$ of $Y.$

In the elliptic base case (see theorem \ref{ellipticbase}) we have two
more types of surfaces: type $E_0$ and type $E_{I_N}.$

Loosely speaking, a {\bf type $E_0$ pseudoelliptic surface} (resp.
{\bf type $E_{I_N}$ pseudoelliptic surface}) is a surface $S$ endowed
with a map $g: S' \to S$ from a log-standard elliptic surface $(S\to
E, Q+F)$ with one marked fibre, and with base $E$ an elliptic curve
(resp. a closed chain of $\PP ^1$'s), and such that the exceptional
set of $g$ is the zero-section $Q.$ A type $E_0$ (resp. type
$E_{I_N}$) pseudoelliptic surface has an elliptic (resp. degenerate
cusp) singularity at $g(Q).$

\subsection{Aknowledgements} I would like to thank my advisor Dan Abramovich, 
whose invaluable and essential incouragement and teachings have been the sole 
reason for my accomplishing these results. The clarity of the
exposition has benefited tremendously from his constant reading.
Of course, any mistake of this work 
is to be blamed on me. Thanks are also due to Kenji Matsuki, for suggesting 
that the log-flips should be toric. On a personal level, I must
express my gratitude to my parents, Angelo and Milva, and my wife Verity.

\section{Preliminaries} \label{prelim}

In this section we will give a brief list of results we need from 
toric geometry and the theory of log-canonical surfaces.
 Varieties are always integral, reduced schemes of finite type over an 
algebraically closed field $k$ of characteristic zero. Let $S$ be a scheme 
of finite type over $k .$ If $X$ is proper over $S$ we set:
$$N_1 (X/S) = \{ \text{ ~1-cycles of } X/S \text{ modulo numerical 
equivalence} \} .$$
We then have a bilinear pairing:
$$Pic(X) \times N_1 (X/S) \to \QQ$$ 
defined by extending by linearity the map that associates $deg_C(\mathcal L)$ 
to the pair $(\cL , C),$ with $\cL$ a line bundle on $X$ and $C$ an effective 
irreducible curve.
We also set:
$$NE(X/S)=\{ B \in N_1 (X/S) \text{: } Z \equiv \sum a_i C_i \text{, } 
a_1 \in \QQ ^+ \cup {0}\}$$
and we denote by ${\overline {NE}}(X/S)$ its closure with respect to the 
euclidean topology in $H^2(X,\RR .)$
 {\bf Kleiman's Criterion of ampleness} states that, if 
$X$ is $S \text{-proper} ,$ 
 a divisor $A$ is $S$-ample on $X$ if and only if $A \cdot x > 0$ for each 
$x \in {\overline {NE}}(X/S) \setminus \{0 \}.$
An {\bf extremal ray} is a ray $R \subset {\overline {NE}}(X/S)$ such that 
if $x_1 + x_2 \in R$ then $x_1 \text{, } x_2 \in R ,$ for each $x_1 , x_2 
\in {\overline {NE}}(X/S).$
 
By the {\bf cone theorem} (see \cite{Mo} for the smooth case and \cite{KMM} for
the general case), given an extremal ray $R$ of
$NE(X)$ there exists an {\bf extremal contraction}, namely a morphism 
$\phi _R : X \to Y$ such that:

\begin{enumerate}
\item ${\phi _R} _* \cO _X \simeq \cO _Y ,$ i.e., $\phi _R$ 
has connected fibres;
\item a curve $C \subset X$ is contracted by $\phi _R$ if and only if its class
$[C]$ in ${\overline {NE}}(X/S)$ is such that 
$[C] \in R .$
\end{enumerate}

\begin{definition}
Let $\phi : X \to Y$ be an extremal contraction such that the codimension of 
the  exceptional set 
$E \subset X$ is $\geq 2,$ and let $D$ be a $\QQ-$Cartier divisor of
$X.$ A variety $X ^+ ,$ together 
with a birational morphism:

$$\phi ^+ : X^+ \to Y$$

is called a {\bf flip } (resp. {\bf log-flip}) of $\phi$ if:
\begin{enumerate}
\item
 $X^+$ has only {\it log-canonical singularities}
(see section \ref{slc} below) 
\item  $K _{X^+}$ 
  (resp. $K _{X^+} +D^+ $ where $D^+$ is the closure of
 $(f^+ \circ f^{-1})(D)$) is $f^+ $-ample.

\item the exceptional set of $\phi ^+$ has codimension  $\geq 2$ in $ X^+ .$
\end{enumerate}

\end{definition}

\subsection{Semi-log-canonical singularities} \label{slc}

By definition, a reduced scheme of finite type $X$ is said to be 
{\bf $\QQ$-Gorenstein} if
$\omega_X ^{[n]}$ is locally free for some $n$. Here
$\omega_X ^{[n]} = ((\omega_X ^{\otimes n})^\vee )^\vee .$ 
  For a $\QQ$-Gorenstein variety, the smallest such $n$
is called the {\it index}.
Following \cite{K-S}

\begin{definition}
A surface $X$ is {\bf semi-smooth} if it has only the following singularities:
\begin{enumerate}
\item  $2$-fold normal crossings with local equation $x^2=y^2$
\item  pinch points with local equation $x^2=zy^2$ 
\end{enumerate}
\end{definition}

and 

\begin{definition}
A {\bf good semi-resolution} of a surface $X$ is a proper map
$g:Y \longrightarrow X$ satisfying the following properties

\begin{enumerate}
\item { $Y$ is semi-smooth}
\item { $g$ is an isomorphism in the complement of a codimension
two subscheme of $Y$}
\item { No component of the double curve $D$ of $Y$ is exceptional for $g$.}
\item {The components of $D$ and the exceptional locus
of $X$ are smooth, and meet transversally.}
\end{enumerate}

\end{definition}

Finally, given a birational map $g:Y \to X, $ we call {\bf discrepancies  
  of} $K_X + D$ {\bf associated to} $g$ those integers $a_i$ such that:
 $\omega_Y^n ( \overline D ) = g^*\omega^{[n]}_X(D +ka_1 E_1 +...+ ka_n E_N),$
where $\overline D$ is the pushforward of $D$ via $g ^{-1} .$
When $D$ is the empty set, we call the $a_i$ just {\bf discrepancies}.

\begin{definition}
A surface $X$ is said to have {\bf semi-log-canonical} singularities if
\begin{enumerate}
\item $X$ is Cohen-Macaulay and $\QQ$-Gorenstein 
\item $X$ is semi-smooth in codimension one
\item The discrepancies $a_i$ of a good semi-smooth resolution of $g:Y \to X$ 
are all
greater than or equal to $-1 .$

\end{enumerate}
\end{definition}

It is not difficult to show that a surface $X$ is {\bf semi-log-canonical}
if and only if it is $S2$ and its normalization $\nu :X ^{\nu} \to X$ is such that 
the pair $(X, K_X + D)$ where $D \subset X ^{\nu}$ is the double curve of 
$\nu ,$ is log-canonical.

\subsection{Toric varieties} 

Toric varities are obtained by suitably patching affine toric varieties, which 
are, roughly speaking, {\it normal } zero sets of binomials. 
Given a lattice $N\cong\ \ZZ ^n$
and a {\it strictly convex rational polyhedral cone}
$\sigma\subset N_\RR = N\otimes _\ZZ \RR$ 
we will denote by $X_\sigma $ tghe affine toric variety associated
with $\sigma$ (see \cite{F})


We denote by $\sigma^{(1)}$ the 1-dimensional edges of $\sigma$.
 The variety $X_\sigma$ is nonsingular if and only if the
primitive points of $\sigma$ form a part of a basis of $N$. 

The toric variety $X_\sigma$ contains an n-dimensional algebraic torus
$T={\mathbb G }^n_m$ as an open dense subset, and the action of $T$ on itself
extends to a linear action on $X_\sigma$ ( hence the alternate name 
{\it torus embedding}). Thus, $X_\sigma$ is a disjoint
union of orbits of this action. There is a one-to-one correspondence
between the orbits and the faces of $\sigma$. In particular,
1-dimensional faces $\RR_{+}v_i$ correspond to codimension~1 orbits
${\mathbb O}_{v_i}$.

 A {\bf fan} $\Delta \subset N$ is a collection of strongly convex rational 
polyhedral cones such that: ~1) each face of a cone in $\Delta $ is also 
a cone in $\Delta$; ~2) the intersection of two cones in $\Delta$ is a 
common face of the two. To fans one associates {\bf toric varieties}, that 
are obtained by suitably patching the affine toric varieties corresponding
to the cones of the fan.

A {\bf toric morphism} $f:X_\sigma\to X_\tau$ is a dominant equivariant
morphism of toric varieties corresponding to a linear map $f_\Delta:
(N_\sigma,\sigma) \to (N_\tau,\tau)$.

In this paper we will write $\langle  f_1,...,f_2 \rangle $ for the 
cone in $N _{\RR}$ generated by the lattice vectors $f_1,...,f_n \in N$ 
for some lattice $N.$

For a toric variety $ {NE}(X/S)$ is a closed cone, hence 
${\overline {NE}}(X/S)=  NE(X/S) .$ Indeed M. Reid in 
\cite{R1} proves:

\begin{lemma} 
If $f: X \to S$ is a proper toric morphism, and if X is proper, then 
 $$ NE(X/S) = \sum \QQ ^+ \OO _w$$
where $\OO _w$ runs through the ~1-dimensional strata of $X$ in fibre of $f.$
Futhermore, if $X$ is projective, $ NE(X/S) $ is spanned by a 
finite number of {\it extremal rays}.
\end{lemma}

 Set $\Delta ^{k} = \{ \text{  k-dimensional cones of }\Delta \} .$

If $\sigma _1 = \langle e_1, ...,e_{n-1}, e_n \rangle $ and 
$\sigma _2 = \langle e_1, ...,e_{n-1}, e_{n+1} \rangle $ and $w$ is the face 
$\langle e_1, ...,e_{n-1} \rangle ,$ we write $ \sigma (w) =
\sigma _1 + \sigma _2$ for 
the cone: $\langle e_1, ...,e_{n-1},e_n ,  e_{n+1} \rangle .$

 Since $\{ e_1, ...,e_{n-1}, e_n \}$ is a $\QQ \text{-basis}$ for the lattice 
$N,$ there is a relation $\sum _{i=1} ^{n+1} a_i e_i =0 .$ Let $I_1 = \{ a_i
\text{ such that } a_i <0 \}.$

If $R$ is an extremal ray (i.e., $R=\QQ ^+ \OO _w$ with  $\OO _w \in NE(X/S) ,$ 
an extremal ray)
in a fan $F,$ write $F_R$ for the fan whose walls 
are $\Delta _R ^{n-1} = \Delta  ^{n-1} \setminus R .$ The corresponding 
toric variety $Y= X(F_R)$ is the contraction of $R .$  

 Define a simplicial subdivision 
$\Delta ^+$ of $\Delta _R$ by defining 
$${\Delta ^+ }^n = 
\Delta _R ^n \setminus \big(\bigcup _{w \in R} \sigma(w) \cup \bigcup _{w\in R 
, i\in I_1} \sigma_i(w)\big),$$  
 
where $\sigma _i(w):=\langle e_1,...,\hat {e_i}, ..., e_{n-1}, e_n,
e_{n+1} \rangle$
M. Reid in \cite{R1} proves: 
\begin{theorem}
 The toric morphism $\phi _1 : X^+ =X(\Delta ^+) \to Y$ corresponding to the 
simplicial subdivision $\Delta ^+$ of $\Delta _R$ is projective and an 
isomorphism in 
codimension ~1. $-R$ is identified with an extremal ray of $X^+$ and 
$\phi _1 = \phi _{-R}$ is the contraction of $-R .$ 
\end{theorem}



\subsection{Toric ~2-dimensional isolated singularities}

Let $\sigma $ be the cone generated by the vectors $f_2$ and $v=kf_1 - nf_2$ 
in the lattice $N=f_1 \ZZ \oplus f_2 \ZZ .$ Without loss of generality, we may
 assume that $k$ and $n$ are 
coprime. Then we can choose a unique $n' \in \ZZ$ such that $0 \leq n' <k$ and 
$n n' \equiv 1 (mod \text{ } k) .$
 Therefore, if $nn'=1-kb$ we can map $N$ isomorphically into itself and 
$\sigma$ to the cone $\sigma ' $ generated by $f_2$ and $v' = kf_1 - n'f_2$
 by means of the matrix:

$$ \bordermatrix {&& \cr
& n &k \cr
 & b &-n'  }  \in SL_2 (\ZZ),$$

thus inducing an isomorphism of toric varieties: 
$X(\sigma) \simeq X(\sigma ').$
 The surface $ X(\sigma ')$ is the normalization of:
$$W =\{ (x,y,z) \in \CC ^3 \text{: } x^k=yz^{k-n'} \}$$ 
and it is isomorphic to the quotient of $\CC ^2$ by $\bmu _k$ via the action:
$$\epsilon (x,y) = (\epsilon ^{n'} x, \epsilon y) $$
for a primitive k-th root of unity $\epsilon .$ Its minimal desingularization 
has as exceptional divisor a chain of rational curves $E_i$ with self 
intersections $E_i ^2 =-a_i\leq -2$ determined by the continued fraction:
$$\frac{k}{n'}=a_1 - \frac{1}{ a_2 -\frac{1}{...-\frac{1}{a_r}} }.$$
 Following \cite{B}, we call such a singularity an  
$A_{n',k} \text{-singularity} $ or a $\frac{1}{k} (n',1)$ singularity.

\subsection{Some toric ~3-dimensional isolated singularities}
 Let $n_1,$ $n_2$ and $n_3$ be generators of a ~3-dimensional lattice $N$ 
and let 
$\sigma$ be the cone $\langle n_1 , n_2 , an_1 + bn_2 +rn_3\rangle .$ 
 Then the affine toric variety $X(\sigma)$ is 
isomorphic to the quotient of $\CC ^3$ by $\bmu _r$ acting as:
$$ \epsilon (x,y,z)  = (\epsilon ^ {a} x, \epsilon ^ {b}y, \epsilon z) ,$$
where $\epsilon \in \bmu _r$ is a primitive ~r-th root of unity.
 Following M. Reid (\cite{R3}), we shall refer to such a ~3-dimenional 
isolated singularity as a 
$\frac{1}{r} (a,b,1)$ 
singularity, thus indicating the order of the cyclic group and the 
weights with which 
it acts.

\section{ Weierstrass Forms} \label{wform}

\subsection{Ampleness of the log-canonical divisor}

Let $\psi: Y \to C$ be a flat family of generically smooth stable
curves of genus $1$ over a smooth curve $C,$ with zero section $Q,$ and
let $\pi : X \to C$ be the surface obtained by  
contracting all the components of the possible singular fibres disjoint from 
$Q$. Then we can express $X$ in a Weierstrass form in the following
way.  Since $\pi$ is proper and flat and since $H^2(X_y, {\cO}_{X,y})=0$, by the 
theorem of base change in cohomology (\cite{H} theorem 12.11) 
$R^1 \pi _* {\cO}_X$
 is locally free, and since its rank is one, it is an invertible sheaf on $C$.
 Let ${\mathcal L} :=  (R^1 \pi _* {\cO}_X  )^{\vee}$ its inverse,
 $\cE := \cO _C \oplus \cL ^2 \oplus \cL ^3 $ and 
${\mathcal P}:={\mathbb P}({\mathcal O}_C \oplus {\mathcal L}^2 \oplus { \mathcal L}^3)$ and 
let:
$$x: \mathcal E \to {\mathcal L}^2,$$
 
$$y:\mathcal E \to {\mathcal L}^3$$ 
and 
$$z:\mathcal E \to {\mathcal O}_C$$
 be the canonical projections onto the given factor.
We have:

\begin{theorem} \label{weirstrassform}
There exist two sections: $g_2 \in H^0(C,\cL ^4)$ and $g_3 \in H^0(C,\cL ^6)$
such that $X$ is isomorphic to the Cartier divisor in $\mathcal P$ 
given by the equation:
$$y^2z=x^3 -g_2 xz^2 -g_3 z^3 .$$
Moreover:

\begin{enumerate}
\item $\Delta =4 g_2 ^3 - 27 g_3 ^2 \in H^0(C, {\mathcal L}^{12})$ is non-zero;

\item 
the sections $g_2 ^3$ and $g_3 ^2$ of ${\mathcal L}^{12}$ do not vanish to  
order $\geq 12$ at any point of $C;$

\item the zero section $Q$ of $X\to C$ corresponds to the section at infinity 
$(x,y,z)=(0,1,0) .$

\end{enumerate}

 \end{theorem}
\begin{proof} The proof is the same as the one of Theorem 1' in \cite{Mu}, 
although they only state the theorem in the case the elliptic fibration has 
no singular fibres.
\end{proof}
  Note that the  equation makes sense since all the monomials that appear 
in it are sections of the vector bundle 
${\mathcal Sym}^3 {\mathcal E}^* \otimes {\mathcal L}^6.$ 
The fact that we assumed that there is no component 
of the singular fibres that is disjoint from $Q$ implies that our $X$ is of 
this form. 
 Notice that by the Leray spectral sequence we get:
$$ \chi ({\mathcal O}_S) = \chi (\pi _* {\mathcal O}_S) - \chi( R^1 \pi _* 
{\mathcal O}_S)= \chi (\cO _C) -\chi (\cL ^{\vee}) =d,
$$
 where $d=c_1(\mathcal L)$. 

We have the important:
\begin{theorem}\label{canbundleform} Let $\psi : Y \to C$ be a
  relatively  minimal elliptic surface,
 with zero section $Q$ (in particular it has no multiple fibres). 
 Then 
$$ \omega _Y \simeq \pi ^* (\omega _C \otimes {\mathcal L} )$$ 

\end{theorem}

\proof See \cite{ B} theorem 12.1.
\qed

Two important consequences are the following:

\begin{corollary} \label{weierstrasscase} For a relatively minimal non-isotrivial 
elliptic surface $\pi: X \to C$ with 
smooth base curve $C$ and zero section $Q$ the divisor:
$$L_X= K_X  + Q + c_1(\pi^* j^* \cO _{\bM _{1,1}} (3))$$
is semi-ample and positive on every curve except for those 
($-2$)-curves of the fibres that do not meet $Q.$
More generally, for the $\QQ$-Cartier divisor:
$$L:=a Q + \pi ^* (\beta +\lambda)$$
for $a\in \QQ _+$ a number with $0<a \geq 1$ and $\beta $ a divisor on $C$, we have:

\begin{enumerate}

\item if  $c_1 (\beta)>0,$ then $L$ is ample;
\item 
\begin{enumerate}
\item if $0<a<1,$  $c_1 (\beta )=0$ and $c_1 (\lambda) >0,$ then $L$ is ample;
\item if $a=1,$ $c_1 (\beta ) = 0, $ and $c_1 (\lambda) >0,$
  the line bundle $L$ is
semiample, and for any irreducible curve $D\subset X$: 
$$L\cdot D=0 \text{  if and only if } D=Q;$$
\end{enumerate}
\item if $c_1 (\beta)<0$ and $c_1 (\beta +\lambda)>0$ then $Q$ is an
  extremal ray.
 
\end{enumerate}
 Note that this includes $K_X +a Q+\pi^*(\alpha)$ for any $\alpha \in
 Div(C)$ with $c_1(\alpha )+2g-2 \geq 0 ,$ if $g$ is the
 genus of $C.$
\end{corollary}
\proof 

In what follows we will use freely that $c_ (\lambda )= -Q ^2 .$

We will first show the general part and then we will reduce to it the statement involving
the canonical bundle.

First one checks easely that:

 $$\begin{aligned} (aQ +\pi ^* (\beta + \lambda))^2 &= 2a(\beta
   +\lambda)\cdot Q +a^2 Q^2\\ & =
 2a[c_1(\beta) + c_1(\lambda)] -a^2 c_1(\lambda) \\ &=2a \big( 2c_1(\beta)
 + (2-a) c_1(\lambda)\big) \geq 2a( c_1(\beta) + c_1(\lambda)) .\end{aligned}$$

Note that this is positive in all three cases, showing that $L$ is big
($L$ is clearly effective).

Also, if $F$ is a fibral divisor (i.e., if it is the class of a
fibre), since $X$ is relatively minimal, we always have:
$$L \cdot F= a >0 .$$

\noindent
 If $D \subset X$ is any irreducible curve, there are two possibilities: 
it either maps 
dominantly onto $C$ or it is a fibral divisor. Thus, we are left to deal
with the former case.

So we may assume that $D$ is an irreducible curve dominating $C .$
If $D \neq Q,$ since $Q$ is effective, we have: $D \cdot Q \geq 0 ,$
and then: 
$$L\cdot D \geq c_1 (\beta +\lambda) +
aD\cdot Q \geq c_1 (\beta +\lambda)>0 $$  

\noindent
in all the cases.

Therefore, all the cases will be differentiated according to the behavior of $L$
on $Q.$

We have: 

$$ L \cdot Q = a Q^2 +  c_1 (\lambda )+ c_1 (\beta )= (1-a) c_1
(\lambda ) + c_1 (\beta )$$

{\bf Case (1) and Case (2) (i):} In this case we have that: $L \cdot Q >0$ and we can
conclude by means of the Nakai-Mosheizon criterion.

{\bf Case (2) (ii):} in this case $L\cdot D \geq c_1 (\lambda )>0 $
for any irreducible (non-fibral) curve $D \neq Q$ and $L \cdot Q = 0$ 

{\bf Case (3):} in this case $Q$ is the only curve on which $L$ is
negative. This entails that $Q$ generates an extremal ray. In fact,
observe that  if $c_1 + c_2 \in [Q] \RR _+ \subset {\overline
  {NE}}(X)$ then $L \cdot (c_1 + c_2) <0.$ But if $c_1$ and $c_2$ are
limits of sequences $C_{1i}$ and $C_{2i}$ of curves that differ from
$Q,$ then $L \cdot C_{ji}>0$ for each $i$ and $j$, and thus $L\cdot
c_j \geq 0$ for each $j,$ which is a contradiction, since it would
imply that $L \cdot (c_1 + c_2) \geq 0.$ So at least one of the $c_i$
is in $[Q] \RR _+ ,$ and again, since $Q$ is the only irreducible curve
on which $L$ is negative, both $c_1$ and $c_2$ must be in $[Q] \RR _+ .$

Now, for the statement about $L= K_X + Q + c_1(\pi^* j^* \cO _{\bM
  _{1,1}} (3))$, note that the by the canonical bundle formula
(theorem \ref{canbundleform}) we have that:$ K_X =\pi ^* ( K_C + \lambda) ,$ where 
$\lambda$ is a divisor 
class associated to $ \mathcal L .$ 
This concludes the proof.

\qed

The following corollary is a special case of the previous one, but it
is worth stating separately, for heuristic reasons. 

\begin{corollary} \label{semiampl} 
Same hypoteses on $\pi: X \to C$ and $C.$
If furthermore the curve $C$ is not rational, then 
the divisor:
$$L_X= K_X  + Q$$
is semi-ample and positive on every curve, except for those 
($-2$)-curves of the fibres that do not meet $Q.$

\end{corollary}

What this means is that, if there are such $-2$-curves in the fibres, in order 
to make the log-canonical divisor ample we have to contract them. This 
explains the necessity of considering {\it Weierstrass forms} for the purpose 
of studying our moduli problems.

\subsection{Log-canonical singularities of Weierstrass equations} \label{LCWE}
 In this section, given a Weierstrass equation $y^2= x^3 + a (t) x + b(t)$ 
with 
$a(t), b(t) \in k[[t]]$ we will only write the low degree terms of 
$a(t), b(t).$ 

After having dealt with the ampleness of the log-canonical divisor,
the second question one needs to address in order to understand when an 
elliptic 
surface is stable is what kind of Weierstrass equations give rise to a 
log-canonical singularity. The following lemma answers this question:

\begin{lemma} \label{localcan}
Let $\pi:S \to \Spec k[[t]]$ be an elliptic surface with zero-section $Q$,
 given in 
Weierstrass form: $ y^2= x^3 + a t^ n x + b t^m,$ with $a,b$ units in $k[[t]]$.
 Assume furthermore that $j\neq \infty .$  
Then the pair $(S,Q)$ is log-canonical if and only if $min( 3n, 2m) 
\leq 12 .$ 
 Furthermore, if 
$F$ is a smooth fibre of $\pi$, the pair $(S,F+Q)$ is log-terminal  
if and only if S 
is,  and if $F$ is a cuspidal fibre, $(S,Q+F)$ is 
never log-canonical.   
\end{lemma} 

\proof 
 The proof of this lemma makes use of the list of Alexeev's of dual graphs 
of log-canonical surface singularities (cf. \cite {K} ch.3) in the sufficient
 direction, and by calculating the discrepancies in the necessary direction.

 If $min (3n,2m) <12$ then the minimal resolution has as dual graph the 
Kodaira graphs: $I_0$, $I_1$, $I_n$, $II$, $III$, $IV$, 
$I_0 ^*$, $I_n ^*$, $II^*$, $III^*$, $IV^*$( see 
\cite {S}, Ch. $IV$ pg.$354$). 

We will refer to table $(IV.3.1)$ pag.$41$ of \cite{M2}
for the singularities corresponding to the various 
Kodaira types.



 In case $I_0,$ and $II$ the surface is smooth, and in case $I_n ^*$
 it has a rational double point singularity, so there is nothing to prove.
 In case $III$ there is only one exceptional curve in the dual graph, and so 
we are in case $(1)$ of Alexeev's list . 
In the case of Kodaira type $IV$, then graph of the resolution 
is of type $A_2$ and 
for the 
Kodaira type $I_n$ the dual graph is of type $A_n$ and again this is in case 
$(1)$ of Alexeev's list. 

In the $*$ cases things become more interesting. 
 $I_0 ^*$ corresponds to a $D_4$ diagram, and this is in case $(2)$ again with 
$(\Delta _1, \Delta _2, \Delta _3)=(2,2,2)$; $I_n ^*$ corresponds to a 
$D_{n+4}$ graph and this is still in case $(2)$ with 
$(\Delta _1, \Delta _2, \Delta _3)=(2,2,N)$; $IV^*$ to an $E_6$ and again this
is in $(2)$ with $(\Delta _1, \Delta _2, \Delta _3)=(2,3,3)$; $III^*$ 
corresponds to an $E_7$ and we are again in case $(2)$ with 
$(\Delta _1, \Delta _2, \Delta _3)=(2,3,4)$; $II^*$ corresponds to $E_8$ and 
one gets case $(2)$ with $(\Delta _1, \Delta _2, \Delta _3)=(2,2,3)$. So these 
 are all even canonical (they are all Du Val singularities).

Now, if $min(3n,2m) \geq 12$, then one can write 
$(n,m)=(4k+n', 6k+m')$ with $n',m' \geq 0$ such that either 
$0 \leq n' \leq 3$ or  $0 \leq m' \leq 5$. 
 One can then consider the rational map:
 $ S' \to S$ given by $(x,y,t)= (x' t^{2k},y' t^{3k},t) $, where $S' \text{: } 
{y'} ^2 = {x'} ^3 + a t^{n'} x + b t^{m'}$. 
 Now:
$$\omega = \frac { dx \wedge dt} {2y} \in \Omega_{k(S)} ^2 $$ 
is a basis for the $k(S)$-module $\Omega_{k(S)} ^2 $, and so is 
 $$\omega ' = \frac { dx' \wedge dt} {2y'} \in \Omega_{k(S')} ^2 $$
for the $k(S')$-module $\Omega_{k(S')} ^2 $, moreover:
$$\pi ^* \omega = \frac { dx' t^{2k} \wedge dt} {2y' t^{3k}} =
\frac { dx'  \wedge dt} {2y' t^{k}}= \frac {\omega '} {t^{k}} \in 
\Omega_{k(S')} ^2 (kE) .$$
 The surface $S'$ is canonical for what has just been showed 
above since 
$min (3n',2m') < 12.$ Therefore $S$ is not log-canonical if  $k>1$, since the 
discrepancies 
would be at least $k$. One would be tempted to say that the $k=1$ case is 
settled and 
thereby claiming the log-canonicity in this case, but indeed the previous 
argument fails 
in this case, since the map $S' \to S$ exhibited above is not proper; 
there might still be exceptional divisors on a completion of $S' \to S$  
with discrepancies $>1.$

 Therefore a detailed analysis of the case $k=1$ is needed.  In order to do 
that one has to 
find a partial resolution and compute discrepancies. To achieve this goal, since 
$y^2= x^3 + a t^{4} x + b t^{6} ,$ with $a$ and $b$ units in $k[t],$ 
can be thought of as the double 
cover of the 
affine plane ramified along the curve 
$R \text{: }x^3 + a t^{4} x + b t^{6}=0$, we can 
simply take an embedded resolution $D\subset M$ of this curve, 
blow-up again at the points in which components of $D$ with odd multiplicity 
meet to obtain $D' \subset M'$, and then take the double covering 
of the resulting surface along the total transform of $R$. This double
covering has canonical singularities.

 The exceptional curve $E$ of a resolution of 
$y^2= x^3 + a t^4 x + b t^ 6 $ is a stable  
curve of arithemetic genus $1$ of self-intersection $-1$ attached to a 
rational curve of 
self-intersection $-1 .$  
Hence the singularity is log-canonical (but not log-terminal).

\qed

For the isotrivial $j= \infty $ case we have:

\begin{lemma} \label{infinity}
The local equation $y^2 = x^2 (x - \lambda t^k )$ with $\lambda$ a
unit in $k[[t]],$ is 
semi-log-canonical if and only if $k\leq 2 .$

\end{lemma}

\begin{proof}

The proof is very similar to the one of lemma \ref{localcan}.
Let us denote by $X$ the surface defined by $y^2 = x^2 (x - \lambda t^k ).$
Let us show first that if $k \leq 2,$ then the singularity is 
semi-log-canonical.
 If $k=0$ then the equation is locally around the $(0,0,t)$ isomorphic to 
$u^2=w^2, $ and it is thus semi-smooth, hence semi-log-canonical.
For $k=1,$ then seeting $z=x-\lambda t$ shows that this singularity is 
isomorphic to $y^2 = x^2 z$ which is again semi-smooth.
 For $k=2,$ after one blow-up we obtain a surface $X'$ and 
as exceptional divisor a nodal curve 
of the same type as the general fibre of $y^2 = x^2 (x-\lambda t^2).$ 
 It is now easy to see that $X'$ is semismooth and that $X$ is 
semi-log-canonical.
 If $k >2 ,$ then we can conclude as in lemma \ref{localcan} that $X$ is not 
semi-log-canonical.

\end{proof}

\begin{definition} \label{locstand} 

We call {\bf standard Weierstrass equation} an equation 
$ y^2= x^3 + a (t) x + b (t),$ which satisfies the 
condition $min( 3 n, 2m) \leq 12 ,$ where $a \text{, } b \in k[[t]]$ and 
$n$ and $m$ are respectively 
the order of vanishing of $a$ and $b$ at $t=0,$ 
or one of the form $y^2 = x^2 (x - \lambda t^k )$ with $k\leq 2
.$      
\end{definition}

 So, we can directly infer from lemma \ref{localcan} and \ref{infinity} 
the important:
\begin{corollary}
A generically stable elliptic surface $X \to C$ mapping to a smooth 
curve $C$ is 
semi-log-canonical if and only
if its local equation around each cusp is a standard Weierstrass equation.  

\end{corollary} 

\subsection{ Types of cuspidal fibres} \label{types}

 In this section we will introduce a bit of terminology and notation we 
will be using in 
the discussion of the special cases in section \ref{special cases}.
 Let $y^2 = x^3 + a x + b,$ with $a \text{, } b \in k[[t]] ,$ 
 be a {\it standard Weierstrass equation} with $j \neq \infty$ (see
 below for the case $j=\infty$).  
Let us denote with $\nu _0 (f) $ the order of vanishing at $t=0$ of a power 
series $f \in 
k[[t]].$ In the following table we set $N := min (3 \nu _0 (a) , 2\nu _0 (b))$
 and the 
first row will give a condition on $N,$ the second will give the Euler 
characteristic of the corresponding Kodaira fibre (see \cite {M2} table 
(IV.3.1) page 41) and the third 
will have a symbol we
 associate 
to the corresponding singularity.
\begin{center}
\begin{tabular}{c|cccccccccc|}
  $N=$& $0$& $6$&  $2$& $10$& $3$&
       $9$ & $4$ & $8$ &$12$\\
\hline
$\chi$ & $n$ & $n+6$ & $2$ & $10$ & $3$ & $9$ & $4$ & $8$ & &\\

\hline 
  type & $I$& $I^*$& $II$& $II ^*$& $III$& $III^*$& $IV$ & $IV^*$ & $L$\\ 
\end{tabular}
\end{center}

here $n$ is related to the order of vanishing of $\Delta$ (see ramark below).

\begin{remark}
 Our notation differs from Kodaira's in that we identify all the $I_n$ and 
$I_n ^*$ to two 
categories, namely: $I$ and $I^*.$
 In Kodaira's notation the former are characterized by the order of 
vanishing of the discriminant 
$\Delta = 4 a^3 + 27 b^2$ (in the $I_n$ case $\nu _0 (\Delta)=n$ and in the 
$I_n^*$ case $\nu _0 (\Delta)=n+6$). The reason for our notation stems
from the fact that 
 in our analysis we need not distinguish 
among them. Also the $L$ case does not appear in the classical litterature, 
because 
it is not a rational double point singularity. It is indeed elliptic, and 
hence it has moduli 
(for instance the $j$ invariant of the exceptional curve)
 as oppossed to the classical cases that do not.

 Note also that if a Weierstrass equation is isotrivial with $j=\infty$ then it is of 
the form $y^2 = x^3 -3 \lambda ^2 (t) x 
+ 2 \lambda ^3(t) $ and it can be transformed to ${y'}^2 = {x'}^2 (x' + 
2 \lambda (t))$ (and viceversa). This shows that the last equation is minimal 
if and only if $\nu _o (\lambda (t)) <2 ,$ which could also be argued by directly 
 computing the discriminant.
 \end{remark}

\section{ Abramovich-Vistoli's fibred surfaces and prestable reduction}

\subsection{Abramovich-Vistoli fibred surfaces} \label{abr}

  D. Abramovich and A. Vistoli in \cite {A-V} define 
{\it famillies of fibered surfaces} in order to compactify the moduli
space of {\it fibered surfaces}, that is to say surfaces $X \to C$
mapping flatly and properly to a curve with stable fibers and with a
number of sections $\sigma _1, ..., \sigma _n$ (or equivalently, to
compactify the moduli space of Kontsevich stable curves into the
Deligne-Mumford stack of stable $n$-pointed curves with fikxed genus).

We will need to borrow the following definitions and results from
their work:

Let $\Gamma$ be a finite group acting on a family of Delgne-Mumford
stable curves $Y\to V,$ over some scheme $V.$ Abramovich and Vistoli (in \cite{A-V})
give the following:

\begin{definition}

 This action is 
{\bf essential} if each $\gamma \in Stab(v),$ for some geometric point
$v\in V,$ acts nontrivially on the fibre $Y_s$ over $s.$  

\end{definition}

\begin{definition} (see {\bf Def. 4.1} in \cite{A-V}) Let $C\to S$ be a
  flat (not necessarily proper) 
family of nodal
curves,
$X\to C$ a proper morphism with one dimensional fibers, and
$\sigma_1,\ldots,\sigma_\nu\colon C\to X$
sections of $\rho$. We will say that
$X \to C \to  S $ is a {\bf
family of generically fibered surfaces} if $X$ is flat over $S$, and the
restriction of $\rho$ to $C_{sm}$ is a flat family of stable pointed curves.
If $S$ is the spectrum of a field we will refer to $X\to C$ as a
generically fibered surface. \end{definition}

\begin{definition} (see {\bf Def. 4.3} in \cite{A-V})
A triple $(U,Y\to V \to S, \Gamma)$ is called a {\bf chart} for a  
family of {\it generically fibred surface} $X\to C \to S$ if there 
is a diagram:
$$\begin{array}{ccccc}
Y&\to& X\times_{C} U & \to & X  \\
\down & & \down          &      &\down \\
V &\to& U  & \to & C \\
\down & & \down &       &\down \\
S       &=&  S  & = &S
\end{array}
$$ 
together with a group action $\Gamma\subset Aut_S(Y\to V)$ satisfying:
\begin{enumerate}
\item The morphism $U\to C$ is \'etale;
\item $V\to S$ is a flat (but not necessarily proper) family of nodal curves;
\item $\rho\colon Y\to V$ is a flat family of stable
$\nu$-pointed curves of genus $\gamma$,
\item the action of $\Gamma$ on $\rho$ is essential; \item we have
isomorphisms of
$S$-schemes $V/\Gamma\simeq U$ and $Y/\Gamma\simeq U\times_C X$ compatible
with the projections $Y/\Gamma \to V/\Gamma$ and $U\times_C X\to U$, such
that the sections $U\to U\times_C X$ induced by the $\sigma_i$ correspond
to the sections $V/\Gamma\to Y/\Gamma$.
\end{enumerate}
The fibre above $p$ is called a {\bf twisted fibre}.

\end{definition}

\begin{remark}
For our purposes, we do not need the chart $(U,Y\to V, \Gamma)$ to be 
{\it minimal}, i.e., we will not need it to satisfy property ~4.
\end{remark}

Let $(U,Y\to V \to S, \Gamma)$ be a chart for  $X\to C \to S,$ then:

\begin{proposition}
Let $\Gamma ' =Stab(v)$ be the stabilizer at the nodal point $v$ of a
gemetric fibre  $V _t$ of $V \to S$ and let $T_1$ and $T_2$ be the
tangent spaces of each branch at the node. Then:
\begin{enumerate}

\item $\Gamma '$ is cyclic and it sends each branch of $V _t$ to
  itself;
\item the generator $\gamma$ of  $\Gamma '$ acts on $T_1 $ 
and $T_2$ by multiplication with a primitive root of unity (of the
order of $\Gamma '$).  

\end{enumerate}

\end{proposition}

\begin{proof}
See \cite{A-V} Proposition 4.5.

\end{proof}

\noindent
In the same situation, Abramovich and Vistoli (in \cite{A-V}) set the following:

\begin{definition}
A chart $(U,Y\to V, \Gamma)$ is called balanced if for any nodal point
of any gemetric fiber of $V,$ the action of the two roots of unity
describing the action of a generator of the stabilizer on the tangent
spaces $T_1$ and $T_2$ of the branches are inverse to each other.

\end{definition}

Let $X\to C\to S$ be a family of generically fibered surfaces, 
$\alpha _1 =(U_1,Y_1\to
V_1,\Gamma_1)$, $\alpha _2 =(U_2,Y_2\to V_2,\Gamma_2)$ two charts; call $
p_i\colon V_1\times_C V_2\to V_i$ the $i^{\rm th}$ projection. Consider
the scheme
$$
I = { \mathop{\isom}\limits_{V_1\times_C V_2}}(p_1^*Y_1,p_2^*Y_2)
$$
over $V_1\times_C V_2$ representing the functor of isomorphisms of the two
families $p_1^*Y_1$ and $p_2^*Y_2.$  Abramovich and Vistoli call the 
{\bf transition scheme } from $\alpha _1$ to $\alpha _2$ the sheme theoretic 
closure of the section of $I$ over the inverse image $\widetilde V$ of
$C_{sm}$ in $V_1\times_C V_2$ corresponding to the isomorphism $p_1^*Y_1
\mid _{\widetilde V}\simeq p_2^*Y_2\mid _{\widetilde V}$ 

\begin{definition} Two charts $(U_1,Y_1\to
V_1,\Gamma_1)$ and $(U_2,Y_2\to V_2,\Gamma_2)$ are\/ {\rm compatible} if
their transition scheme $R$ is \'etale over $V_1$ and $ V_2$.
\end{definition}

\begin{definition} A\/ {\bf family of fibered surfaces} $$
\cX \to \cC  \to  S$$
is a
family of generically fibered surfaces $X\to C\to S$ such that $C\to S$ is
proper, together with a collection $\{(U_\alpha, Y_\alpha\to
V_\alpha,\Gamma_\alpha)\}$ of mutually compatible charts, such that the
images of the $U_\alpha$ cover $C$.
Such a collection of charts is called an\/ {\em atlas}.

A family of fibered surfaces is called\/ {\em balanced} if each chart in its
atlas is balanced.

The family of generically fibered surfaces $X\to C\to S$ supporting the
family of fibered surfaces $\cX\to \cC\to S$ will be called a family of {\em
coarse fibered surfaces}.
\end{definition}

We can now state the theorem of theirs that we are going to be using here:

\begin{theorem}\label{fibredsurfaces}

 Let $\cX _\eta \stackrel{\pi}{\to} \cC _\eta\to \eta$ be a {\em balanced 
stable} fibered
surface, with induced map $f_\eta: \cC _\eta \to \bM _{g,n},$ 
 with sections ${\sigma _i }_{\eta}\subset \cX _\eta$ and with sections 
${D_i}_\eta \subset \cC _\eta .$
 Then there is a finite extension of discrete valuation rings
$R\subset R_1$ and an extension
$$
\begin{array}{ccc}
\cX _\eta\times_\Delta \Delta_1 &\subset & \cX _1 \\
\down            &       & \down \\
\cC _\eta\times_\Delta \Delta_1 &\subset & \cC _1 \\
\down            &       & \down \\
\{\eta_1\} & \subset & \Delta_1, \end{array}
$$

with $\Delta =\Spec R$ and $\Delta _1 = \Spec R_1 ,$
such that:
\begin{enumerate}

\item $\cX _1\to \cC _1\to \Delta_1$ is a balanced stable family of fibered
surfaces with sections $\sigma _i ;$

\item there is a regular map $f_1 : \cC _1 \to \bM _{g,n}$ extending 
 $f_\eta \circ p_1,$ where $p_1 :\cC _\eta\times_\Delta \Delta_1 
\to \cC _\eta$ is the natural projection; 

\item $f_1$ is Kontsevich stable

\item $\omega _{\cX _1 / \cC _1} (\sum \sigma _i+ \sum \pi_1 ^* D_i)
\otimes f_1 ^* \cA$ 
is ample, for some ample line bundle $\cA$ on $\bM _{g,n} .$  

\end{enumerate}
 
The extension is unique up to a unique isomorphism, and its
formation commutes with further finite extensions of discrete valuation
rings.

\end{theorem}

\begin{proof}
See \cite{A-V} pg.20. prop. 2.1. and pg.28 prop. 8.13
\end{proof}

\subsection{ The extension lemma and Prestable reduction} \label{extlemma}

Suppose one is given an elliptic surface $\cX _{\eta} \to \cC _{\eta} \to \eta$
over the generic point $\eta$ of a DVR scheme $\Delta,$ and with induced
$j-$map $j _{\eta}: \cC _{\eta} \to \bM _{1,1} .$ If $\cX _{\eta} \to \cC
_{\eta}$ has no cuspidal fibres, then we can apply the theorem of
Abramovich and Vistoli (cf. Theorem \ref{fibredsurfaces} above), to
extend $\cX _{\eta} \to \cC _{\eta}$ and $j _{\eta}$ over the whole
$\Delta .$
 In case there are cuspidal fibres (and our surface is in Weierstrass
 form), the strategy we will adopt is to
 temporarely replace these with twisted fibres and mark the
 corresponding points on the base curve $\cC _{\eta}.$  

Once we do that, we want to be able to go back (at least generically) 
to Weiestrass forms.
 To this aim 
we prove the following:

\begin{lemma} \label{Extension Lemma} 
Let $S$ be the spectrum of a two-dimensional 
complete 
regular local ring, $S= \Spec k[[s,t]],$ let $W= S \setminus p,$ with 
$p=V(s,t)$ the closed point, set $U= S \setminus V(t),$  $\eta$ the generic
 point of $\Spec k[[s]]$ 
and let 
 $\pi :X _W \to W$ be a family of curves of genus $1$ with zero
section, 
such that:
\begin{enumerate}
\item  there is a map $j: S \to { \bM }_{1,1},$ the $j 
\text{-invariant};$

\item the family $X_W \mid _U ,$ is a family of stable elliptic 
curves;

\item $X_W \mid _{W_\eta}$ is an elliptic surface in minimal Weierstrass 
form.
\end{enumerate}

then $ X  _W$ extends over the whole $S,$ to a family $X\to S$ 
whose fibre over $s=0$ is an elliptic surface in minimal Weierstrass form.

\end{lemma}

\proof  Let us consider the line bundle ${\cL}'$ on $W$ whose dual is: 
$R^1 \pi _* {\cO}_W$. Since $S$ is non-singular, this line bundle 
extends to a unique 
line bundle $\cL$ over the whole $S$. By the Theorem of the base change in 
cohomology, 
${\cL}^* \otimes k (\eta) \simeq (R^1 \pi _* {\cO}_{X}) \mid
_{W_\eta}$, and thus the sections $g_2$ and $g_3$ of 
$H^0(W_\eta, (\cL \otimes k (\eta))^4)$ and of
 $H^0(W_\eta, (\cL \otimes k (\eta))^6)$ respectively, extend to two sections 
$a \text{ and } b$ of $H^0(S, \cL ^4)$ and of $H^0(S, \cL ^6)$ respectively.
So we can consider:

  $$X \text{: } \left\{ y^2 = x^3 + a x + b \right\} \subset \PP ( 
{\cO}_S \oplus 
\cL ^2 \oplus \cL ^3).$$
 
Note that since $a^3$ and $b^2$ are sections of the sixth power of the Hodge 
bundle 
$ \cL ^6$, $a$ and $b$ are defined up to the transformation:
 $(a,b) \to (\lambda ^4 a, \lambda ^6 b)$, where $\lambda \in \cO _S ^* .$

\begin{claim}\label{Claim 1} Let $V:=V(t)$, $A:= V(a)$ and 
$B:=V(b)$, and let  $A= \bigcup A_i \cup h V$ and 
$B= \bigcup B_i \cup kV$  
be respectively the decompositions in irreducible components of $A$ and 
$B$; then $V$ (if $h$ and $k$ are non-zero)
is the only irreducible component that $A \text{ and } B$ can share.

\end{claim}

  In fact, if they had another irreducible component in common,
 say $H,$ then , since in $ k[[s,t]] $  every prime ideal of height one is 
principal,
 there would be an $h \in k[[s,t]]$ such that $H=V(h)$ and so:

$$ X: y^2 = x^3 + h^{n} a' x + h^{m}b' ,$$

\noindent
where $a' {\text{ and }} b'$ are non zero on $H.$ 

The assumption that the curves in the family are stable away from $V=V(t)$, 
 is equivalent to the fact either one of $a' \text{ and } 
b'$ is nonzero on H and that $2n$  and $3m$ are divisible by $6$.
 In fact, if that were not the case, by pulling back our family $X$ via the 
map:
$$k[[s,t]] \to k[[s,t]] \otimes \kappa (h)$$
where $\kappa (h)$ is the residue field at the prime $(h)$, 
we would produce a family of curves that is not stable and that cannot be 
reduced to 
a stable one by changing $a$ and $b$ via the transformation 
$(a,b) \to ( g^{-2k} a, g^{-3k} b)$ for some $g \in k[[s,t]]$ and $k \in 
\ZZ $. 
 Therefore $2n$ and $3m$ are 
divisible by $6$ and we can apply the above mentioned transformation with 
$g=h,$ therefore $H$ does not exist. 

  Furtermore we may assume, according to the {\it purity lemma} of 
Abramovich-Vistoli (see {\bf Lemma 2.1}, pg.3 in \cite{A-V}), that 
the components $A_i$  do not intersect the $B_i$ 
away from $V.$ Indeed every point of $U$ satisfies the hypotheses of the 
purity lemma. If there were a point $p \in U$ such that $p \in A_i \cap B_i ,$ 
 since $j_U : U \to \bM _{1,1}$ is well defined and extends 
the map $j_{U \setminus \{p\}} :  U \setminus \{p\} \to \bM _{1,1}$ defined by 
the family over $U \setminus \{p\},$ according to the purity lemma
there would exist a 
stable family $\cX _U \to U$ extending $\cX _{U \setminus \{p\}} 
\to {U \setminus \{p\}} ,$ and thus $p$ could not be in $A_i \cap B_i
,$

So, next step is to show:

\begin{claim} $A_i  \cap B_i \cap V= \emptyset .$

\end{claim}

 Let us assume that some $A_i$ does meet some $B_j$ along $V .$
That means that there is some point $q \in V$ such that $a$ and 
$b$ are both zero at $q$. 

  We can write $a'=t^n a''$ and $b'=t^m b''$ where $a'' \text{, } b'' \in
 k[[s,t]] \setminus 
(t)$ do not vanish identically along $V$, and the non negative integers 
$n \text{ and } m$ 
can be zero if, respectively $A$ and $B$ don't contain $V$.
  Since we are assuming that some $A_i=V(a_i)$ does meet some $B_j=V(b_j)$
 along $V$, there must be some point $q \in V$ such that $a''$ and 
$b''$ are both zero at $q$.
 The $j \text{-map}$ of the statement is now a map:
$$j: \Spec k[[s,t]] \to {\bM}_{1,1}$$
that at this point has the form:

$$j(s,t) = \frac {1728 a^3 t^{3n}}{4 a^3 t^{3n} + 27 b^2 t^{2m}}.$$

This restricts to a map:
$$j_{G}: G \to {\bM}_{1,1}$$

 of the same form for any divisor $G$ of $\Spec k[[s,t]]$.
In particular we can restrict $j$ to $A_i$ and $B_j$, to get  two morphisms:
$$j_{A_i }:  A_i \to {\bM}_{1,1}$$

and 

$$j_{B_j}: B_j \to {\bM}_{1,1}$$
   
Since these two maps are induced via restriction by $j$, they have to 
coincide at $p$, 
 But from this we infer a contradiction, since $j_{A_i}\equiv 0$  
and $j_{B_j} \equiv 1720 .$  Therefore they cannot meet along 
$V$, because by hypothesis $j: \Spec [[s,t]] \to {\bM}_{1,1}$ is well defined.

\qed

Let $M \subset \PP ^r $ be a projective scheme. 
\begin{definition}
We call {\bf Kontsevich prestable} a flat family of maps over a scheme 
$S$ 
$f: \cC \to M $ if $\cC \to S$ is a flat family of nodal curves.

\end{definition}

and
\begin{definition}
We call {\bf quasiminimal} an elliptic surface 
$X \to  C$ such that $X \mid _{C_{sm}} \to C_{sm}$ is {\it in minimal
  Weiestrass form}, 
where $C_{sm} $ 
is the smooth locus of $C $
 and such that, for each point $p \in C _{sing} $  there is a chart $(U,Y\to V,
\Gamma)$ with $X \times _C V \simeq Y/\Gamma $
centered at $p.$

\end{definition}

The first step towards proving the stable reduction theorem is to  
show the following theorem

\begin{theorem} \label{fict}
Let $\Delta $ be the spectrum of a DVR, $\eta \in \Delta$ the 
generic point and 
$(\cX _{\eta} \to \cC _{\eta} , \cQ _{\eta} , f_{\eta} : \cC _{\eta} \to \bM _{1,1}),$ 
be a triple consisting of a relatively minimal elliptic surface 
$\cX _{\eta} \to \cC _{\eta}$ with section $\cQ _{\eta}$ and 
Kontsevich-stable map to moduli $f_{\eta}.$  Then we can find, after possibly 
a finite base change  $\Delta ' \to \Delta ,$  a map of 
$\Delta ' \text{-schemes}$ $\cX ' \to \cC '$ such that:
\begin{enumerate} 
\item the fibre over the 
special point  $0 \in \Delta ' ,$ $ \cX _0 \to \cC _0$ is the semi-log-canonical union 
of {\em relatively quasiminimal }elliptic surfaces with section $\cQ _0;$
\item the double curves of $\cX _0$ are either stable or twisted fibres;
\item
$f_{0} : \cC _{0} \to \bM _{1,1}$ is {\it Kontsevich prestable}. 
 
\end{enumerate}
\end{theorem}

 \begin{proof} 

We are given 
$(\cX _{\eta} \to \cC _{\eta} \to {\eta} ,f_{\eta} :\cC _{\eta} \to \bM _{1,1}) .$
 
As mensioned in the introduction to this section, we can mark the
points of $\cC _{\eta}$ corresponding to cuspidal fibres of 
$\cX_{\eta} \to \cC _{\eta}:$ let $\Sigma _{i \eta} $ be such divisor on
$\cC _{\eta}.$ 
 We can then  use the theorem of Abramovich-Vistoli (cf.{\bf Theorem}
 \ref{fibredsurfaces}), to
 get a triple $(\cX ' \to \cC ' \to \Delta ', \cQ ' ,f': \cC ' \to \bM
 _{1,1} )$ consisting of a family of fibred surfaces
$\cX ' \to \cC ' \to \Delta ',$ (with sections $s_i:\Delta ' \to \cC'$
extending $\Sigma _{i \eta}$) a $\QQ -\text{Cartier}$ divisor $\cQ'$
and a regular map $f': \cC ' \to \bM _{1,1} $ such that condition
~2 and the condition (stronger than ~3 above) that $f' _0 $ together
with the sections $s_j$ be Kontsevich stable hold.  

This family coincides with our family of elliptic curves in 
minimal Weierstrass form $\cX _{\eta} \to \cC _{\eta} \to \eta $ over $\cC ' _{\eta} \setminus 
\bigcup _i {\Sigma _i} _\eta .$  Let $\Sigma _i $ = $\overline {{\Sigma _i} _{\eta}},$ 
the closure in $\cC '$ of ${\Sigma _i} _\eta ,$ and $\Sigma = \bigcup
\Sigma _i .$

Hence, we can apply lemma \ref{Extension Lemma} to remove the twisted
fibres lying above $\Sigma$ and to replace them with the cuspidal curves induced
by the original ones lying above $\Sigma _{\eta},$ so that the generic
fiber $\cX ' _{\eta '}$ is now isomorphic to the original surface $\cX
_{\eta} \to {\eta}.$ According to lemma 
\ref{Extension Lemma} $\cX ' _0$ consists of {\it quasiminimal} elliptic 
surfaces. The map $f' _0$ may now be Kontsevich unstable (since
we have removed the sections $s_i$), but it is clearly {\it
  prestable}. This concludes the proof. 
\end{proof}

We shall refer to theorem \ref{fict} as the {\bf Prestable Reduction Theorem}.

\begin{remark}
 Lemma \ref{Extension Lemma} allowes us to remove the extra sections we 
added along the cuspidal 
fibres to use the stable reduction theorem of \cite {A-V} (see section $4$). 
 But this comes with a price: when we do so, the family of maps 
$j: \cC \to {\bM}_{1,1}$ may no longer be stable.

In fact,  if there is a rational component $C$ of the 
central fibre $\cC _0$ on which $j$ is constant, and which meets the rest of the 
central fibre in only one point ( e.g., if the surface 
$X=\cX _0 \mid _C\to C$ contains two 
cuspidal curves or more and meets the rest of $\cX _0$ transversally
along one fibre), then $j$ is no longer stable. Indeed we will see in proposition \ref {extr ray} 
that 
the components of $\cX$ that map onto such curves are not stable in
the sense of \bf Definition \ref{logstabsurf}.

\end{remark}

\subsection{The log-canonical divisor, Extremal rays and semiampleness } \label{iso}

As mentioned in the previous subsection, we have to deal with
components of the central fibre $\cX _0$ 
that lie above a rational component of $\cC _0$ contracted by the $j$-map. 
 These are then {\it isotrivial quasiminimal} elliptic 
surfaces, attached to at least one component in two possible ways: 
either along 
a smooth or a twisted fibre. By the very definition, for such a surface the 
$j$-map is constant, therefore the log-canonical divisor (for triples) does 
not have the contribution coming from $j^* \cO _{\bM _{1,1}} (1) .$
Hence the log-canoncial divisor for triples equals $ K _{\cX} + \cQ$ on 
such components.

  The zero section of such a component is always a log-flipping 
extremal ray, if this component is attached 
to the rest of the central fibre only along one fibre, and will be contracted 
by the log-canonical bundle if the isotrivial component is attached along two 
fibres.

Indeed we have, even more generally, that this happens, even if the 
component $X \to C$ is {\it not isotrivial}, in case the log-canonical 
divior is $L= K _{\cX} + \cQ$ (the one for pairs). This last divisor coincides 
with the log-canonical divisor for triples on the isotrivial components, as we 
remarked above.
This will turn out to be useful 
when dealing with pairs $(X,Q)$ only. 

\begin{lemma} \label{extr ray}
Let ${\mathcal X}_0 = X \cup X'$ be a decomposition of the central fibre 
where $\pi:X \to C $ is a {\em quasiminimal} elliptic surface,
 where $C \simeq \PP ^1$ a smooth rational curve. 
Then:
\begin{enumerate}

\item If $X$ is attached to $X'$ only along one fibre $G$, then:
$$L \cdot Q =-1$$
where $Q$ is the zero-section and 
$L= (K_{\cX} + \cQ) \mid _X=G+Q+K_X$ is the log-canonical divisor;

\item If $X$ is attached to $X'$ along two fibres $G_1 \text{ and } G_2$, then:
$$L \cdot Q =0$$
where $Q$ is the zero-section and 
$L= (K_{\cX} + \cQ) \mid _ X=G_1 + G_2+Q+K_X$ is the log-canonical divisor.

\end{enumerate}
\end{lemma}
\proof

 Let us assume first that the attaching fibres are stable. In this case, 
the zero section $Q$ goes entirely through the smooth locus of the morphism 
$\pi : \cX \to \cC$. Therefore, if we denote by $K_{X/C}$ the relative 
canonical divisor:
$$K_{X/C} \cdot Q + Q^2 =c_1(\omega _Q \otimes \pi ^*\omega _C)=0$$
moreover the attaching fibres are reduced and therefore linearly equivalent 
to the general fibre $F$. Hence, since $\pi^* K_C= -2F$:
$$L \cdot Q = (K_{X/C}+Q +G_1 + \pi ^* K_C)\cdot Q=(-2F +F)\cdot Q = -1$$
in the case of one attaching fibre, and:

$$L \cdot Q = (K_{X/C}+Q +G_1+G_2 + \pi ^* K_C)\cdot Q=(-2F +F+F)\cdot Q = 0$$

\noindent
in the case of two stable attaching fibres.

 If the attaching fibres are not stable, then there is a diagram:

$$\begin{array}{ccc} Y &\stackrel{f} {\rightarrow} & X \\
 \dar & &\dar \\
C' &\stackrel{\phi} {\rightarrow} &C
\end{array}$$

\noindent
such that $\phi :C' \to C$ is a finite morphism with Galois group 
$\Gamma$, $f: Y \to X $ is also finite and 
$Y \to C'$ is a relatively minimal elliptic surface. 
 Moreover $\phi :C' \to C$ is ramified at two points,
and:

\begin{enumerate}

\item  $\Gamma \simeq \bmu _k$ if there is only one twisted attaching fibre 
with monodromy of order $h$, 
in this case the branch points are $0$ (where $(X_0)_{red} =G_1$) and infinity;

\item  $\Gamma \simeq \bmu _k $ with $k = l.c.m.(k_1,k_2 )$ when there are two 
twisted attaching fibres of monodromy of orders $k_1$ and $k_2 ,$ respectively.

\end{enumerate}

Let $Q'=f^*Q$ and $F'=f^*G .$ 
 Then, in case (1), by the projection formula we have:
$$ L \cdot Q = \frac {1}{k} f^* L \cdot Q' = 
\frac {1}{k} [(K_{X'} - (k-1)F' + F' + Q') \cdot Q']=-1$$

since $f^* G_1 = F'$ and by Riemann-Hurwitz:
$$f^* K_X = K_{X'} - (k-1)F'.$$

In case (2), the Riemann-Hurwitz formula now 
reads:
$$f^* K_X = K_{X'} - (k_1-1)F'-(k_2-1)F'$$
and ones again, by means of the projection formula, we can conclude:
$$ L \cdot Q = \frac {1}{k} f^* L \cdot Q' = 
\frac {1}{k} [(K_{X'} - (k_1 -1)F' +(k_2 -1) F'+k_1 F'+k_2 F' + Q') \cdot Q']=0.$$
\qed

 It will turn out to be useful to generalize these computations for any
 {\it quasiminal elliptic surface} $X \to C$ with no limitation on
  the number of twisted fibres and the genus of the base curve. 

Let $\cX _0 = X \cup X_R$ a decomposition of the central fibre
such that $X$ meets $X_R$ in $r$ attaching fibres (stable and
twisted) $\bigcup _{i=1} ^r G_i.$ Let $Q$ be the zero-section of 
$\pi :X\to C$ (where $C$ is the base curve of $X$). 
Let $k_i$ be order of the monodromy 
group around $G_i ,$ and $g$ the arithmetic genus of $C.$

We have:

\begin{proposition} \label{genlogcan} 

With these hypotheses and notations, one has:

$$ L_X \cdot Q =  2g-2 +  r .$$

\end{proposition}

\begin{proof}

The proof is very similar to the proof of the previous lemma.
 Define $\overline r $ to be the integer that
equals $r$ if $r \equiv 0 \text{(mod 2)}$ and equals $r+1$ if $r
\equiv 1  \text{(mod2)}.$

Let $p_i := \pi (G_i),$ $p_{r+1}\in C$ some extra point,
 and let $I$ to be the set that equals 
$\{ p_1, p_2 , ..., p_ {\overline r} \}.$  Now let $g:C' \to C$ be the finite
ramified covering  with Galois group
 $\bmu _k, $ with $k=lcm(k_i)$ totally ramified at the $p_i 's ,$ with
 $i \in I.$

Because of the  underlying structure of quasiminimal elliptic
 surface,there is an  atlas $(U_i,Y_i\to V_i,\Gamma _i \simeq \bmu
 _{k_i})$ such that  $Y_i \to V_i$ is stable. 
 In particular, the normalization of the pull-back $X' \to C'$ of $X \to C$ to $C'$ will be
 a family of stable curves, and therefore, if we denote by $Q'$ the
 pull-back of the zero section $Q \subset X,$ $Q$ must be entirely
 contained in the smooth locus of $X'.$  Let $f: X' \to X$ be the
 natural morphism.

If we indicate by $g'$ the genus of $C',$
 it must then be the case that $(K_{X'} +Q')\cdot Q' =2g'-2.$ But by
 Riemann-Hurwitz applied to $C' \to C$ we obtain:

$$2g'-2 = k(2g-2) +  {\overline r} (k-1)$$

The Riemann-Hurwitz formula applied to $f$ reads:

$$ K_{X'} = f^* K_X + \sum _{i=1} ^ {\overline r} (\frac{k}{k_i}-1) G_i '$$ 
where $G_i ' =(f^* G_i)_{red} .$
 Therefore, if we let $Q' = f^* Q :$

$$ L_X \cdot Q = \frac {1}{k} \big( K_{X'} - \sum _{i=1} ^ {\overline r} (\frac{k}{k_i} -1) G_i '
+ 
\sum _{i=1} ^r \frac{k}{k_i}  G_i'  + Q' \big) \cdot Q' ,$$

where $ G_{\overline r}= \varnothing $ if $r \equiv 0 \text{ ( mod 2)} $ and
$  G_{\overline r}$ is the fiber over the point $p_{r+1}$ otherwise.

 Thus, if we let $\epsilon (N)$ be the function over the integers that
 equals $0$ if $N \equiv 0 \text{ ( mod 2)} $ and $1$ otherwise we have: 

$$ L_X \cdot Q =\frac {1}{k} (k(2g-2) + {\overline r} (k-1) -\epsilon
(r)k +{\overline r})= 2g-2 + r. $$

\end{proof}

\begin{remark} \label{self int}
Note that a similar computation can be carried out to measure 
the failure of $X$ to satisfy adjunction, i.e., to compute the {\it
  different} (see \cite{K}, page 175 for a definition).
 Another point worth observing is that if $X$ has only one attaching fibre,
 in $X \mid _U$:

 $$Q ^2 = \frac{1}{k} f^* Q \cdot Q' = \frac{1}{k} {Q'} ^2 $$
so that this is the contribution a twisted fibre with monodromy $k$ gives 
to the self intersection of $Q .$

\end{remark}

\begin{remark}\label{canint}

Note that we have showed, en passant, that:

$$ f^* L_{X}=L_{X'} + (1-k \epsilon (r)) G_{\overline r} .$$

\end{remark}

\begin{corollary}  \label{twistedampleness}

Same hypotheses as in Proposition \ref{genlogcan}, then:

\begin{enumerate}
\item if $2g-2+r > 0,$ then $L_X$ is ample;
\item if $2g-2+r = 0,$ and $X\to C$ is non-isotrivial, then  $L_X $ is
  semiample and for any irreducible curve $D$:
 $$L_X \cdot D =0 \text { if and only if } D=Q ;$$ 
\item if $2g-2+r <0$ then $Q$ is an extremal ray.

\end{enumerate}

In fact, more generally the same holds for the $\QQ$-Cartier divisor:

$$M := K_X + \sum G_i + a Q$$

with $ 0<a<1$ any rational number.

\end{corollary}

 \begin{proof}

Recall from Remark \ref{canint} that:

$$ f^* L_{X}=L_{X'} + (1-k \epsilon (r)) G '_{\overline r}, $$
    
\noindent
where $f:X' \to X$ is as constructed in the proof of Proposition \ref{genlogcan}.
Therefore: 
$$L_X ^2 = L_{X'}^2 + 2(1-k \epsilon (r)) L_{X'} \cdot  G '_{\overline r}$$

\noindent
since ${G '_{\overline r}} ^2=0.$ 
The surface $X'$ is now in minimal Weierstrass form, and therefore we can apply
the {\it canonical bundle formula} (cf. theorem \ref{canbundleform}),
and write: 

$$f^* L_{X}= ({\pi '}^* (K_{C'} + \lambda '+ \sum _{i=1} ^r p_i ' +(1-k
\epsilon (r)) p_ {\overline r} ')+Q',$$

\noindent
where $p_i' \in C'$ denotes the point whose fibre via $\pi '$ is $G_i
.$

Note that $c_1(K_{C'} +  \sum _{i=1} ^r p_i ' +(1-k \epsilon (r)) p_
{\overline r} '=k (2g-2+r).$

We can then conclude (1), (2) and (3) by means of corollary \ref{weierstrasscase}.

We can now conclude also the more general statement about $M$
appealing to the same corollary, and by observing that if one defines
the $\QQ$-Cartier divisor on $X'$:
$$M' := K_{X'} + \sum G_i ' + a Q' $$
one has that:

$$ f^* M= M' + (1- k \epsilon(r)) G_{\overline r} = ({\pi '}^* (K_{C'} + \lambda '+ \sum _{i=1} ^r p_i ' +(1-k
\epsilon (r)) p_ {\overline r} ')+a Q' .$$

\end{proof}

\section{The special cases and standard elliptic surfaces} \label{special cases}

\subsection{Special cases} \label{specas}

 


  Let $\cY \to \cC \to S$ be an elliptic threefold in Weierstrass form 
 $$\cY = V(y^2 - x^3 - ax -b) \subset \PP (\cO _{\cC} \oplus \cL ^2 \oplus 
 \cL ^3),$$ 
 where $\cL = \omega _{\cY \mid \cC} .$
For any fibre $\cC _t$ of the base $\cC$ one can define an 
integer valued  function on $\cC _t$ as: 
$$ N_{\cC _t} (q)= min (\nu _q (a ^3 \mid _{\cC _t}), 
\nu _q (b^2 \mid _{\cC _t}) ),$$ 
where $\nu _q(g)$ is the order of vanishing of any section 
$g\in \cO _{\cC _t}$ at $q \in \cC _t.$ In particular we will simply write 
$n(q)$ for this function defined for the central fibre $\cC _0.$
 
If $\cC$ contains a negative rational curve $E$ (not necessarely irreducible),
 one can contract it to obtain a morphism 
$\rho :\cC \to \cC ' .$
On $\cC '$ one can then construct an elliptic 
threefold 
$\cY ' \to \cC '$ in the following manner.
 Let $p'$ be the image of $E$ and $W'= \cC ' \setminus p' .$ 
Then, one naturally has a Weierstrass equation on $W'$ induced by 
$\cY \mid  _{W'}.$ In fact $\rho$ induces an isomorphism 
$\rho _W : W \to W'$ where $W= \cC \setminus E ,$ and, 
letting $\cL '_{W'} = \rho _* \omega _{\cY _W / W'}$ 
and  $a'_{W'}$ and $b'_{W'}$ respectively the push-forward of the 
section 
$a \mid _ W$ of $\cL ^2 \mid _W$ 
and $b\mid _W$ of $\cL ^3\mid _W,$ 
one can then define an ellitpic threefold in 
Weierstrass form: 
$${\cY}' = V(y^2 - x^3 - a'x -b') \subset \PP (\cO _{\cC '} 
\oplus {\cL '} ^2 \oplus 
 {\cL '}^3).$$ 
 Since $\cC '$ is smooth at $p'$ one can then extend the line bundle $\cL ' _{W'}$ 
uniquely to a line bundle $\cL '$ on $\cC ',$ and correspondingly uniquely 
extend the sections $a' _{W'}$ and $b' _{W'}$ to global sections $a'$ and 
$b'$ of 
${\cL ' }^2$ and ${\cL '}^3$ respectively.
 To put it differentely, we take the {\it saturation } 
$\cL ' = (({ \rho _* \cL}) ^{\vee})^{\vee}$ and the corresponding sections 
$a'$ and $b'$ to 
construct the Weiestrass equation presented above.

  In general, given an elliptic threefold $\cX \to \cC$ 
(not necessarely in Weierstrass 
form) with no multiple fibres, but possibly twisted fibres, 
and a contraction of a negative curve of the base 
$\rho : \cC \to \cC ',$ one can define 
an elliptic threefold  $\cX '$ over $\cC'$ as the elliptic threefold 
$\cX '$ extension of a local Weierstrass model around the image point. 
 That is to say we choose a Zariski neighborhood $Z$ of $E$ in $\cC$  
such that the only possible twisted fibres are above the nodes of $E.$ 
Let $W$ be a Zariski neighborhood of $E$ in $Z$ such that $\cX _W$ has a 
Weierstrass representations away from the nodes of $E.$ $\rho :\cC' \to \cC$ 
induces a map $\rho : W \to W'$ to a Zariski neighborhood of $p' \in \cC ' .$
 Now apply the construction above to $W'$ and $W.$ 
This construction may lead to a threefold which does not have 
log-canonical singularities, even if $\cX$ had log-canonical singularities to 
begin with.   

 Let $\cX \to \cC \to \Delta$ be a ~1-parameter family of elliptic surfaces with 
section $\cQ ,$ with general fibre $\cX _{\eta}$ in minimal Weierstrass 
form. Assume furthermore that there is a component $X_1 \to C_1$ of the 
central fibre $\cX _0$ that maps to a rational nodal not necessarely 
irreducible curve $C_1 \subset \cC ,$   
that $X_1$ is attached (transversally) to another 
component $X_2\to C_2$ of $\cX _0$ along only one fibre $G$ 
(twisted or stable). One can then  
contract $C_1$ in $\cC$ to a point $p\in \cC '$  and therefore emulate 
the construction of 
$\cX '$ above. The question 
one must ask oneself is: when does it happen that the surface $X_2 ' ,$ image 
of $X_2$ via the contraction of $X_1$ has 
log-canonical singularities? Proposition \ref{specialcas} will answer 
just that. We call the cases when this occur, the {\bf special cases}.

In this chapter we will use (mostly implicitely) the following well known result:

\begin{lemma} \label{nonsing}
The contraction of $C_1$ in $\cC$ produces a surface $\cC '$ that is smooth 
in a Zariski neighborhood of the image point.  

\end{lemma}

\subsection {The Idea} 

The idea of the construction of the log-canonical model, is based on the 
possibility (cf. section \ref{abr}) of taking an \'etale neighborhood $V$ of  
$ p\in \cC$ such that the pull-back of $\cX $ to
 $V$ has a nodal fibre at $p,$ 
i.e., of finding a chart (not necessarely minimal) 
$(Y\to V, U, \Gamma )$ at least of a Zariski neighborhood $W$ of $p \in \cC .$
 Therefore $Y \to V$ has a Weiestrass model and we can naively
 contract the rational components of the central fibre of $\cC \to S$
 and write ``explicit'' equetions for the analityc singularity we thus obtain.

\subsection{Standard elliptic surfaces}

Let us assume that $C_1$ is irreducible and at set $p =C_1 \cap C_2 .$ 
 Let moreover $\cC _R= \overline {\cC _0 \setminus C_1} .$
 It might happen that the fibre of $\cX$ over $p$ is {\it twisted}.
But we can find a chart $(V,Y \to U, 
\Gamma  )$ such that the fibre over the point lying above $p$ is stable.
 On $U$ we can define two divisors $A$ and $B$ as follows. 
Let $W$ and $W'$ as in section \ref{specas}, so that 
$\cX \mid _ {W' \setminus p}$ is 
in Weierstrass form ${y^2= x^3 + ax +b} \subset 
\PP (\cO _{W' \setminus p} \oplus \cL _{W' \setminus p}^2 \oplus 
\cL _{W' \setminus p})$ Since $\cC$ is normal at $p,$ the 
two sections $a$ and $b$ extend to sections all over $W.$ 
 Set $A= V(a^3)$ and $B=V(b^2).$ 

\begin{lemma}\label{nonvan}
We have: $p \notin A \cap B.$ 

\end{lemma}

\begin{proof}
Let $U$ be a Zariski neighborhood of $p$ such that the only twisted fibre 
of $\cX \mid _U$ is at $p$ and that there are no cuspidal fibres. 
Let $(V,Y \to V', \Gamma  )$ be a chart of $\cX _U \to U$ 
as in section \ref{abr}. Then the thesis holds, because if 
both $A$ and $B$ passed through $p,$ the family $Y \to V'$ 
could not consist of stable curves.
\end{proof}

Let $c: \cX \dashrightarrow \cX ' $ be the contraction 
constructed in the previous section and $ X_2 ' =c_* X_2 .$ 
Let $C_1$ be an irreducible $(-1)$-curve in $\cC$ meeting the rest of 
 $\cC$ transversally 
in only one point $p.$ Let $C_2$ be the irreducible component 
of $\cC$ meeting $C_1 .$

At this point we ask ourselves what happens to the self intersection 
of the zero section $Q_2$ of the component $X_2$ to which $X_1$ was attached. 
The answer is given by the following, where we let $Q_2 '$ be the
zero-section of $X_2 ' \to C_2 ':$

\begin{proposition} \label{self-int}
  The self intersection ${Q_2 ' }^2$ of $Q_2 '$ in $X_2 '$ is:
$${Q_2 ' }^2 = Q_2 ^2 +{Q_1  }^2 $$
 where  $Q_2 ^2$ 
is taken in $X_2$ and ${Q_1  }^2$ is taken in $X_1 .$
\end{proposition} 
\begin{proof}
 Let us consider how the number $Q_2 ^2 + Q_1 ^2$ changes 
after the contraction of $X_1.$ Note that it is the same
as $(Q_2+ Q_1)\cdot \cQ ,$ where $\cQ$ is the divisor of $\cX$ swept
out by the zero sections of the fibres of $\cC \to S.$ On the
other hand, $Q_1 + Q_2 \equiv \cX _0 \cdot \cQ  - \cX _R \cdot \cQ$ 
where $\cX _R = \overline {\cX _0 -(X_1 +X_2)}$ is the rest of
the central fibre $\cX _0 ,$ and $\equiv$ denotes numerical equivalence. 
Hence 
$$Q_2 ^2 + Q_1 ^2 =(\cX _0 \cdot \cQ  - \cX _R \cdot \cQ)\cdot \cQ .$$
This number does not change after performing the contraction, 
since this operation does not touch
$\cX _R$ and $\cX _t \cQ $ is constant as $t$ varies among the geometric 
points of $\Delta$ because of the flatness of $\cX \to \Delta .$ So 

        $$(Q_2')^2 = Q_2^2 + Q_1^2 $$

\end{proof}

Our goal is to relate the invariant 
$\sum _{q \in \cC \setminus \{p\}} N_C(q)$ with the self 
intersection of the zero section, therefore suggesting tha the latter 
should be treated as an 
{\it invariant of the geometry near the zero section}. This is taken care by:

\begin{proposition} \label{twelve}
Let $X \to \PP ^1$ an elliptic surface whose 
singular fibres are all standard Weierstrass equations, except for possibly 
having at most one twisted fibre over $0 \in \PP ^1$
 whose monodromy group is isomorphic to 
$\bmu _h .$ Assume also that away from $0\in \PP ^1,$ the surface $X$ is in Weierstrass 
form (not necessarely minimal).
Then $12 Q^2 = -[\sum _{q \in \PP ^1 \setminus {0}} n_C(q) + deg(j)]$ 
and in particular $\sum _{q \in \PP ^1 \setminus {0}} n_C(q) + deg(j) \leq 12$ 
if and only if 
$Q_1 ^2 \geq -1 .$  
\end{proposition}

\begin{proof}
Let us assume first that $X\to \PP ^1$ is not isotrivial.
  Noether's formula reads:

$$\chi (\cO _X) = \frac{ K_X ^2 + \chi (X)}{12}$$ 
and therefore, since $K_X ^2 =0$ and $\chi (\cO _X) = c_1 (\cL)$ 
(see section \ref{weirstrassform}) we obtain:
$$\chi(X) = 12  c_1 (\cL)= -12 Q^2 .$$ 
 
 But since $X\to \PP ^1$ is locally trivial in the euclidean (resp. \'etale) 
topology away from the cuspidal fibres, the (\'etale) 
Euler characteristic $\chi (X)$ 
must equal:

$$\chi (X) = \chi (X_{\eta}) \chi (f^{-1} (C \setminus \delta) )
+ \sum _{p\in \PP ^1} \chi (X _p), $$
where $X_{\eta}$ is the generic fibre and $\delta =\{ p \in \PP ^1 : 
X_p \text{ is singular } \} .$ Hence, since $\chi ( X_{\eta})=0, $

$$  Q^2 =-\frac{1}{12} \sum_{p \in \delta}  \chi (X_p) .$$

Performing a case 
by case analysis of  the table in 
section \ref{types} on sees that in all cases but for $I_n$ and $I_n ^*$ 
one has that $\chi (X_p) = N_p$ and in the latter cases $\chi (X_p)= N_p +
n .$  

 Therefore, since by Corollary $IV$.4.2 of \cite{M2} one has that 
$deg(j) = \sum _{n\geq 1} n(i_n + i_n ^*)$ where $i_n$ and $i_n ^*$ indicate 
respectively the number of fibres of type $I_n$ and $I_n ^* ,$ we conclude 
that 

$$  Q^2 =-\frac{1}{12} [\sum _{p \in \delta} N_p + deg(j)].$$

 If there is a twisted fibre $G = (X \otimes k(p))_{red}$ 
(for 
some point $p \in \PP ^1$), there is a ramified 
morphism of order $h$ (the monodromy of $G$) of the base 
$s:C\simeq \PP ^1 \to \PP ^1,$ totally 
ramified at $0$ and at infinity,
  and a diagram:
$$\begin{array}{ccc}
X' & \to &X \\
\downarrow & &\downarrow \\
C & \stackrel{s}{\to} &\PP ^1
\end{array}$$
such that $X' \to C$ is stable.
 
Thus we conclude by means of the same argument as above and the observations 
that $Q^2 =\frac { {Q' }^2}{h}$ (see remark \ref{self int}) and that 
if $q' \in C$ is a point mapping to $q$ via $s,$ $N_{q'} = \frac {N_q}{h}.$

If $X\to \PP ^1$ is isotrivial, there is nothing to prove if $j\neq
\infty,$ since in that case there are no $I$ or $I ^*$ fibres.

\end{proof}

\begin{proposition} \label{specialcas}

 The singularity of $X_2 '$ in a neighborhood of the fibre over 
the point $p'$ to which 
$C_1$ gets contracted is of type $y^2 = x^3 + a' x +b'$ with 
$min( \nu _{p'} ({a' \mid _{{C_2'}}}^3), \nu _{p'} 
({b' \mid _{{C_2'}}}^2) ) = 
\sum _{q\in \cC \setminus I} N_{C_1} (q)+deg(j \mid _{C_1}) .$
 In particular $X_2 '$ is log-canonical if and only if 
$ \sum _{q\in \cC \setminus I } N_{C_1} (q) +deg(j \mid _{C_1}) \leq 12 .$

\end{proposition}

\begin{proof}

 We can construct $\cX '$ as in section \ref{specas}.

 Let $C_2 ' = \rho _* C_2$ and $p= C_1 \cap C_2 .$

According to proposition \ref{twelve} we have that:

$$12 {Q_2 '}^2 = -[\sum _{q \in C_2 '} 
N_{C_2 '}(q) + deg(j _{C_2 '})]$$

We remark that the fibre over $p'$ will no longer be twisted.
Note also that even though $j'$ might only be a rational map on the whole 
$\cC ^c ,$ (the surface resulting from contracting $C_1$)
 it extends to a regular map when restricted to any smooth 
curve $B \subset \cC ^c .$

and:
$$ \begin{array}{cc}
12 (Q_1 ^2 +{Q_2 }^2) &= -[\sum _{q \in C_2  \setminus \{ p \}} 
N_{C_2 }(q) + deg(j_{C_2})] -[\sum  _{q \in C_1  \setminus \{ p \}} 
N_{C_1 }(q) + deg(j_{C_1})]\\
 &\\
&\\

\end{array}$$

On the other end, according to proposition \ref{self-int} one has:
$12 (Q_1 ^2 +{Q_2 }^2) =12 {Q_2 '}^2.$ Hence: 

$$ \begin{array}{cc}
 \sum _{q \in C_2 '} N_{C_2 '}(q) &= \sum _{q \in C_2  \setminus \{ p \}} 
N_{C_2 }(q)  +\sum  _{q \in C_1  \setminus \{ p \}} 
N_{C_1 }(q) + deg(j_{C_1})]\\
 &\\
&\\

\end{array}$$

since away from $p \in C_2$ we have not changed $X_2 \to C_2,$ and
thus that $deg(j_{C_2})= deg(j _{C_2 '}).$ 
 The proposition now follows from the observation that, since we have
 not changed $X_2 \to C_2$ away from $p,$ it must be that 
$N_{C_2 '}(q)= N_{C_2 }(q)$ for every $q\neq p .$

\end{proof}

What this proposition says is that one can contract all these components 
to begin with, before proceding with the stable reduction. It is thus 
convenient to introduce 
the following generalization of the concept of minimal Weierstrass form:

\begin{definition}\label{standard}
Let $X \to C$ a be an elliptic surface 
mapping to a nodal irreducible curve 
$C.$ Such a surface will be named 
{\bf standard elliptic surface} if the local 
equation around each fibre above $C_{sm}$ is a 
{\it standard Weierstrass equation} (see 
definition \ref{locstand}) and if, 
for each point $p \in C _{sing} $  there is a chart $(U,Y\to V,
\Gamma)$ centered at $p$ with $X \times _C U \simeq Y/\Gamma .$

\end{definition}

We have:

\begin{theorem} \label{divcont}
Same hypotheses as in theorem \ref{fict}; then there is a finite 
base change of DVR schemes $\Delta ' \to \Delta$ and a $\Delta '$-family 
of elliptic surfaces $(\cX ' \to \cC ', \cQ ' )$ 
with section $\cQ ' ,$  
such that:
\begin{enumerate}
\item if $\eta '$ is the generic point of $\Delta ' ,$ $\cX ' _{\eta '} \simeq 
\cX _{\eta}$;
\item the central fibre $\cX ' _0 $ is composed of 
{\it standard elliptic surfaces}; 
\item  if $X \to \PP ^1$ is a component of 
$\cX ' _0$  attached along only one fibre (twisted or stable), then 
$\cQ \mid _X ^2 < -1$ in $X.$ 
\item if we ask condition ~3 only of isotrivial components, then 
there is a well defined regular map $j' : \cC ' \to \bM _{1,1} $
and if we set $f'= \pi ' \circ j',$ then  
$\omega _{\cX ' / \Delta '} (Q) \otimes {f'}^* \cO _{\bM _{1,1}} (3)$
is ample away from those isotrivial components of the central fibre 
$\cX ' _0$ 
that meet transversally the rest of $\cX _0 '$ along one or two fibres 
of $\cX _0 ' \to \cC _0 '.$

\end{enumerate}

\end{theorem}
\begin{proof}

According to theorem \ref{fict} there exists a finite morphisms of DVR schemes 
$\overline {\Delta} \to \Delta $ and a triple 
$(\overline \cX \to \overline \cC ,\overline \cQ ,  \overline f : 
\overline \cX \to \bM _{1,1} )$ satisfying conditions ~1, ~2 and 
$\overline \cX _0 \to \overline \cC _0 $ is a union of 
{\it realively quasiminimal} ellitpic surfaces.  

 Let $\overline B$ be a connected component of $\overline \cC _0$ consisting 
of a tree of rational curves (we call such a component $\overline B$
a {\bf tree-like component}), meeting the rest of
 $\overline \cC _0$ trasversally in one point, and such that 
in $\overline X_ {\overline B}  = \overline \cX \mid _{\overline B} \to 
\overline B$  
we have ${\overline \cQ } _{\overline B} ^2 \geq -1.$

The proof will be by induction on the number $N$ of such components.

Let us asume first that $N=1.$
Let $\overline C_1,$ be a ``leaf''(i.e. an irreducible component 
furthest away from $\cC _0 \setminus B$) of $\overline B$ 
 Let $\rho :\overline \cC \to \overline \cC '$ be the contraction of 
$ \overline C_1.$ If $ \overline X_1 \to  \overline C_1$ is isotrivial, 
the map $ \overline j : \overline \cC \to \bM _{1,1}$ descends to a map 
$ \overline j ' : \overline \cC '\to \bM _{1,1}.$
 
 Let 
$\overline C_2 $ be the component of $\cC _0$ meeting 
$\overline C_1,$ let $\overline C_2 '$ be its image in $\overline \cC '$ 
via $\rho$
 and $p' = \rho ( \overline C_1) .$ We can apply proposition \ref{specialcas} 
to $\overline C_1$ and    
find a threefold $ \overline \cX ' \to \overline \cC '$ that 
satisfies ~1 and ~2 (and with a regular map  
$ \overline j ' : \overline \cC '\to \bM _{1,1}$
if $C_1$ was $j$-trivial) and such 
that the local equation of the fibre above $p'$ is in {\it standard Weiestrass
 form}. 
 The image of the curve $\overline C_2 $ is again a leaf itself,
 because there are no other tree-like components. 
According to Lemma \ref{twelve}, if $\overline \cQ '$ denotes the 
divisorial push-forward 
of $\overline \cQ :$ 
$${\overline \cQ ' \mid _{\overline X_2 '}} ^2 ={\overline \cQ  
\mid _{\overline X_2 }} ^2 
  + {\overline \cQ  \mid _{\overline X_1 }} ^2 .$$
If this number happens to be $\geq -1 ,$ then we can appy the procedure to 
$C_2 ' ,$
which is now attached to the rest of ${\overline \cC }' _0$ only at one point, 
since we have contracted $\overline C_1 .$ We can inductively iterate this 
procedure untill we get to a component for which the 
self intersection of the zero-section is less than $-1 .$ 

Let us now do the general case: we assume we know the result for any
number $k$ of tree-like components with $k<N .$ Let $\overline B$ be a
tree-like component which has the property that it is attached to all
the other tree-like components at one end only (such a component must
exist, even though it may only consist of one edge). Then we can apply
the previous argument to this tree-like component, and ``prune'' all
its edges. Now the number of tree-like components is $N-1 ,$ and we
can then conclude by induction.

 If we perform these operations only for tree-like components $\overline B$ 
for which $\overline X \to \overline B$ is isotrivial, then
 we get the 
desired triple $(\cX ' \to \cC ', \cQ ', 
f:\cX '\stackrel{\pi '}{\to }\cC' \stackrel{j'}{\to } \bM _{1,1})$ on 
$\Delta ' = \overline \Delta .$ In this case, the claim about the ampleness of 
 $\omega _{\cX ' / \Delta '} (Q) \otimes f^* \cO _{\bM _{1,1}} (3)$
 away from those isotrivial components of the central $\cX ' _0$ 
that meet transversally the rest of $\cX _0 '$ in one or two fibres, is a 
consequence of theorem \ref{fibredsurfaces}. 
In fact we have only perfomed birational operations to  
some surfaces with constant $j$-invariant (namely, those surfaces $X$ 
for which $\cQ \mid _ X ^2 \geq -1$)  
which met the central fibre ${\overline \cX }_0$ along just one 
fibre, and so $\cC ' \to \bM _{1,1}$ is still Kontsevich prestable, i.e., 
$\omega _{\cC'/ \Delta '} \otimes {j'}^* \cO _{\bM _{1,1}} (3)$ 
is semiample and ample away from those components $C \subset \cC ' _0 $ 
that meet the rest of $\cC ' _0 $ in one or two points and that are 
$j-$trivial; 
and since $\omega _{\cX '/\cC'}(Q)$ is relatively ample, we are done.

\end{proof}

\begin{definition}
A {\it pair} $(X \stackrel{\pi}{\to} C, 
Q)$ will be called {\bf strictly prestable} if:
\begin{enumerate}
\item 
$X \to C$ is a  {\it standard elliptic surface} 
and if for each rational component 
$B$ of $C$ that meets the rest of $C$ in only one point;
\item 
${\cQ \mid _ {X \mid _B}} ^2 < -1 .$
 \end{enumerate}
 Similarly, a {\it triple }
$(X \stackrel{\pi}{\to} C, 
Q, f: X \stackrel{\pi}\to C \stackrel {j}{\to} \bM _{1,1})$ 
will be called {\bf stricly prestable}
if condition $(1)$ above hold and if $(2)$ holds only for 
those components $B$ for which 
$X \mid _ B \to B$ is {\it isotrivial}.
\end{definition}

Therefore, theorem \ref{divcont} says that one can perform the {\bf strictly 
prestable reduction} of a family of minimal elliptic surfaces,
possinly after a finite base change.

\begin{remark}\label{sum}
From the formula 
$  Q^2 =-\frac{1}{12} \sum_{p \in \delta}  \chi (X_p) $
(see proposition \ref{twelve}) and from table IV.3.1 of \cite {M2}, we infer 
the following table for the contribution of a Kodaira fibre to $Q^2 :$

\begin{center}
\begin{tabular}{c|cccccccccc|}
 
fibre type & $I$& $I^*$& $II$& $II ^*$& $III$& $III^*$& $IV$ & $IV^*$ &L \\ 
\hline
cont. to $Q^2$ &$0$& $-\frac{1}{2}$& $-\frac{1}{6}$&$-\frac{5}{6}$&
$-\frac{1}{4}$& $-\frac{3}{4}$&$-\frac{1}{3}$&$-\frac{2}{3}$& -1\\
\end{tabular}
\end{center}

Note also that the contribution to the self intersection is exactly 
$\frac{-N}{12}$ (see the table in section \ref{types}).

\end{remark}

 \section{ The Toric Picture}

\subsection{One attaching fibre}

In this section we will show that we can perform the necessary log-flips and 
log-canonical contractions torically.
 The first step towards understanding the toric picture in the case of 
one attaching 
fibre, is 
to understand what it looks like on the base curve $\cC ,$ or 
equivalently on the zero-section $\cQ .$

 Let $R$ be a discrete valuation ring (DVR),  $\Delta = \Spec (R) ,$ $\eta$ 
its generic 
point and $0$ its special point.  

For a toric variety $Z$ with torus $T$ we 
write $D_Z$ 
for the complement of the torus $D_Z = Z \setminus T .$

Let $ \cC \to \Delta$ be  a family of nodal curves $\cC .$ 
Assume that the central fibre $ \cC _0 $ has 
a rational component $C_1$ meeting the rest of $\cC _0$ transversally 
in only one point $p , $ and assume that the 
singularity at $p$ is an $A_{k-1}\text{-singularity}. $
Let $C_2$ be the rest of $\cC _0$, $S$ a divisor of $\cC $ meeting $C_1$ 
transversally in only one point and $\cI \subset \cO _{\cC}$ the ideal sheaf 
of $C_1 \cup C_2 \cup S .$

%
%

%
%
%

\begin{lemma} \label{torbasonefib}
Let $\cC \to \Delta$ as before.
Then there is a Zariski open neighborood $U$ of $C_1$ in $\cC ,$ a 
~2-dimensional 
toric variety $Z$ and an \'etale map $t:U \to Z$ such that:
\begin{enumerate}

\item the fan of $Z$ is $F = \langle f_1, f_1 + k f_2\rangle  \cup  
\langle f_1 + k f_2 ,f_2 \rangle $
in the lattice $N = f_1 \ZZ \oplus f_2 \ZZ ;$

 \item the pull-back via 
$t$ of the ideal of $D_Z $  is the ideal $\cI (U). $ 
\end{enumerate}
\end{lemma}

\begin{proof}
Let $U \subset \cC$ be a Zariski open neighborhood of $C_1$ such that $U \cap 
(C_1 \cap S) =p , $ 
and such that $S$ is not entirely contained in $U .$
We can now contract $C_1$ to a smooth point $q.$ Let $c : \cC \to \cC ''$ 
be such contraction. 
Let $C_2 '' = c' _* C_2$ 
and 
$S'' = c'_* S$ and let $\cI ''$ be the ideal sheaf of $C_2'' \cup S'' .$
 Therefore, we can find an \'etale neighborhood $U''$ of $q\in \cC''$ 
and a map $t'' : U'' \to \bbA _{k} ^2$ to the toric variety $\bbA _{k}^2 .$
 We shrink $U$ and $U''$ if necessary so that $c^{-1}(U'') =U .$

Let $Z'$ be the toric variety whose 
fan $F'$ is the union of the cones:
$\sigma _1 =\langle f_1, f_1 + f_2 \rangle  ,$ $\sigma _2 =\langle f_1
+   f_2 , f_1 +2 f_2\rangle ,$ 
....$\sigma _k = \langle f_1 + (k-1) f_2 , f_1 + k f_2\rangle ,  
\sigma _{k+1}=\langle f_1 + k f_2 , f_2 \rangle $
 in the lattice 
$N= f_1 \ZZ \oplus f_2 \ZZ ;$ and let $U'$ be normalization of $Z'$ in the 
function field $k(U)$ of $U.$ Thus, by definition we have morphisms 
$b: U' \to U $ and 
$t':U' \to Z'$ such that ${t'}^* \cO (D_{Z'}) \simeq \cO (C_1' + C_2 ' +S')$ 
where $C_i ' = b^* C_i$ and $S' =b^* S.$ The surface $U'$ is in fact the 
minimal resolution of $U.$ 

 The morphism $t'$ induces a rational map $t :U \dashrightarrow Z .$ 
Let $W$ be its graph, then we have a commutative diagram:

$$\begin{array}{cccccc}  U'&&\stackrel{t'}{\longrightarrow} &&Z'  \\ 
&\searrow \alpha&&&\\
b\downarrow &&W&&\downarrow \beta\\
& p_1 \swarrow  &&\searrow p_2&&\\
 
  U &&\stackrel{t}{\dashrightarrow} && Z
\end{array}$$

Where $\alpha$ is the morphism whose existence is ensured by the minimality 
of $U' .$ 
If we show that $p_1$ and $p_2$ are finite, 
then $t$ is in fact a regular morphism 
and by construction it satisfies the properties of the thesis.

 But $\beta \circ t' = p_2 \circ \alpha ,$ and since ${t'} ^* F = E,$ if 
$E$ is the exceptional divisor of $b$ and $F$ the one of $\beta, $ we
have  that $p_2$ is finite ($t'$ is \'etale). Similarly one concludes that $p_1$ 
is finite and therefore $t$ is a morphism.

\end{proof}

Our goal is to look at a {\it strictly prestable} 
family $\cX \to \cC$ of elliptic 
surfaces, 
and in particular the base curve $\cC$ is either a family of {\it Kontsevich 
prestable} curves (in the case of moduli of triples) 
or a surface obtained from a 
Kontsevich prestable family by contracting ~{-1}-curves.
So it may happen  that $\cC$ is singular, but the singularities are of type 
$A_{k-1}.$ 
In the particular context of lemma \ref{torbasonefib}, the total space of the
 family of 
base curves $\cC$ might have a singularity of type $A_{k-1}$ at 
$p,$ in which case we 
want to be able to find a finite ``toric'' morphism from a smooth 
surface. We have:

\begin{lemma} \label{toricbasext}
In the hypothesis of lemma \ref{torbasonefib}, then there exist a 
ramified cyclic covering $f:V \to U$  of order $k$ ramified at 
 $p$ and along some section $S$ and an \'etale map 
to a ~2-dimensional toric variety $t' : V \to Z'$ such that:
\begin{enumerate}
\item the fan of $Z'$ is $F'=\langle  f_1 ' , f_1 ' + f_2 '\rangle  \cup 
\langle f_1 ' + f_2 ', f_2 '\rangle  ;$
\item the pull-back via $t'$ of the ideal of the toric divisor $D_{Z'} $ is 
the ideal of $C_1 ' \cup C_2 ' \cup 
S'$ where $C_2 ' = f^* C_2$ and $S' = f^* S ;$
\item $f$ induces a toric morphism $f : Z' \to Z $ given by $f_1 '=f_1$ and 
$f_2 ' = k f_2$

\end{enumerate}

\end{lemma}
\begin{proof}

 Choose a divisor $S$ that meets $C_1$ transversally in only one point, 
by shrinking $U$ 
if necessary, we can apply 
lemma \ref{torbasonefib} and find a toric 
variety $Z$ and an \'etale morphism $t:U\to Z$ such that $t^*\cO_Z (\OO _{f_1}
 + \OO _{f_1 + kf_2} + \OO _{f_2}) \simeq \cO _U 
(C_1 + C_2 + S) .$ Since 
$t^* \cO _Z (f_2) \simeq \cO _U ( S)$ the cyclic covering corresponds to 
taking the 
sub-lattice $N' =f_1 ' \ZZ \oplus f_2 ' \ZZ$ of $N=f_1  \ZZ \oplus f_2  \ZZ$ 
with 
$f_1'=f_1$ and 
$f_2'=kf_2 .$ $U'$ is the toric variety given by the fan 
$F' = \sigma _1 ' \cup \sigma _2 '$ union of the two cones generated by 
$\{ f_1', f_1 '+f_2 '\}$ and $\{ f_2', f_1' +f_2 '\}$ respectively in the
 lattice $N' .$ From the description of $Z' \to Z$ we can easely read off the 
ramification.
 
Take $V$ to be the normalization of $U$ in $k(Z')$ with the induced 
morphisms $f: V \to U$ and $t': V \to Z'.$  By construction, they satisfy 
the hypotheses of the lemma.
\end{proof}

Let $(\cX \to \cC \to \Delta, \cQ )$ be a family of
 {\it strictly prestable} elliptic surfaces 
with zero section $\cQ .$  Let $X_1 \subset {\mathcal X}_0$ be component of 
the central fibre 
$\cX _0$, mapping down to a rational curve $C_1 \subset 
{\mathcal C}_0;$ also, let $X_2 \subset {\mathcal X}_0$ be the 
rest of ${\mathcal X} _0$ to which $X_1$ is attached along one and only one 
fibre, which is either stable or twisted. 
We have:

\begin{proposition} \label{one fibre}
 
 Let $Q_1 = \cQ \mid _{C_1} $ and let $\cS$ be a divisor of $\cX$ meeting
 $X_1$ 
transversally in only one stable fibre. 
Then there is a Zariski open neighborhood $\cU$ of $Q_1$ in $\cX$ and a toric 
variety 
$\cY$ with an \'etale morphism $T:\cU \to \cY$ such that:

$T^*\cO _{\cY} (D_{\cY}) \simeq \cO _{\cU} ( X_1 + X_2 + \cS)$

\end{proposition} 

\begin{proof}
Let $U$ be a Zariski naighborhood of $C_1$ with an \'etale 
map from a toric variety $t:Z \to U$ as in lemma \ref{torbasonefib}. 

Let $f:V \to U$ be as in lemma \ref{toricbasext}. Since $V$ is smooth 
(at least in a neighborhood of $C_1 \cap C_2$), 
 by the purity lemma of
Abramovich and Vistoli (cf. {\bf Lemma 2.1}, \cite{A-V} ) the
attaching fibre $F := (X_1 \cap X_2)_{red}$ is a quotient of a stable 
curve by a 
cyclic group of order $h$ dividing $k$ and the normalization 
of the pull-back of $\cX \mid _U$ to $V$ is a family of stable curves.

Therefore we have a 
diagram:

$$\begin{array}{ccl} {\cV} & \stackrel{s}{\to} & \cX \mid _U  \\ 
\pi ' {\dar} & & {\dar}\pi\\
  V &\stackrel{f}{\to}
 & U
\end{array}$$

such that $ \pi' :\cV \to V$ is a family of stable curves and the action of 
the Galois group $\Gamma$ of $f: V \to U$ extends to an action 
on $\cV$ such that $ \cV /\Gamma \simeq \cX \mid _U.$
  
 According to lemma and \ref{toricbasext} we have  
Zariski neighborhoods $U$ of $C_1$ in $\cC$ and  $U'$ of $C_1 '$ in $\cC'$ and 
toric varieties with \'etale maps $U \to Z$ and $U' \to Z'$ such that the 
induced map $Z' \to Z$ is toric. 

It is thus enough to show that the pull-back 
$\cY'$ of $\cV $ to $Z'$ is toric and that the induced action of the Galois 
group 
 $\Gamma$ on $\cY'$ is an action by a subgroup of the torus. In this case 
the quotient $\cY / \Gamma \to \cX \mid _U$ would be the \'etale map in the 
statement of the proposition.

Indeed, the total space of 
the normal bundle $\cN _{\cQ /\cX '} \mid _U$ of the zero section
 $\cQ$ in $\cX$ is such a space, since every line bundle on a 
toric variety 
is a toric bundle, according to the proposition of page ~{63} in \cite{F}. 
Thus if we can show that the action of $\Gamma$ on $\cC '$ lifts to an action 
of 
$\Gamma$ on  $\cY$ and that $\Gamma$ acts as a subgroup of the torus 
on $\cY$. 

But $\Gamma$ is a finite group and for each point $(u,v) \in \cY$
it acts on the element $v$ of finite dimensional vector space $\cN _{Q/\cX }
 \otimes \kappa (u),$ where $\kappa (u)$ is the residue field of $u \in V.$
  Therefore $\Gamma$ must act linearly on $\cN _{Q/\cX }$, i.e.,
 via multiplication by a character, and since $\Gamma$ acts as a subgroup 
of ${\CC ^*} ^2$ on $U$ we are done.  

\end{proof}

In the following lemma, $Q_1 ^2$ denotes the self-intersection 
as $\QQ$-divisor of $Q_1$ in $X_1.$

\begin{lemma} \label {toric one fibre}
In the same hypotheses as proposition \ref{one fibre}, the fan $F$ of the 
toric variety $\cU,$  is $F=\sigma _1= \langle e_1 , e_1 +e_2 , e_3 \rangle  
\cup 
\sigma _2 =\langle e_3 ,e_1 +e_2,w\rangle $ in the lattice 
$N :=  e_1 \ZZ \oplus  w \ZZ \oplus e_3 \ZZ$ where 
 $\{e_1 ,e_2 ,e_3 \}$ is the standard basis of $\RR ^3$
and $w =\frac{1}{k} (  e_2 + n e_3).$ 
 Here $-n=k Q_1 ^2$ 
(see remark \ref{self int}) and 
$\cC$ has an 
$A_{k-1}$ singularity around $C_1 \cap C_2 .$
\end{lemma} 

\begin{proof} 
  Let 

$$\begin{array}{ccc} {\cX}' & \stackrel{f}{\to} & \cX  \\ 
\pi ' {\dar} & &{\dar}\pi\\
{\mathcal C}' &\stackrel{s}{\to}
 & \mathcal C
\end{array}$$
be the diagram as in \ref{one fibre}. In an \'etale neighborhood of $C_1$ 
the surface $\cC$ is described by the quasi-projective toric surface defined 
by the fan $\Delta$ obtained by the two cones generated by 
$\{f_1 , v=f_1 + k f_2 \}$ and  $\{f_2 , v \}$ respectively in the 
lattice $L= f_1 \ZZ \oplus f_2 \ZZ +  (k f_2 ) \ZZ ,$ where 
 $\{ f_1 , f_2 \}$ is the standard basis for $\RR ^2 .$ This is the surface 
we refer to in \ref{one fibre}.

The fan of $V$, i.e., the toric \'etale neighborhood of $C_1 '$ in $\cC '$ 
as in \ref{one fibre}, is obtained by taking the the two cones generated by
$\{ f_1 ' , f_1 ' + f_2 ' \}$ and $\{ f_2 ' , f_1 ' + f_2 ' \}$ in 
the sub-lattice 
$L' = L + \frac {1}{k} ( f_2 ') \subset L,$ where $f_1 '=f_1$ and 
$f_2 ' = k f_2.$ Set $\Lambda = Hom _{\ZZ} (L, \ZZ)$ and 
  $\Lambda ' = Hom _{\ZZ} (L', \ZZ),$ and let $\langle ,\rangle  : 
\Lambda '/\Lambda 
\times L/L' \to \QQ /\ZZ$ the natural pairing.

The Galois group 
$\Gamma$ as in \ref{one fibre} is isomorphic to $L /L' \simeq \ZZ / k \ZZ$ 
and its action on 
 the ring of functions of $V$ is given by:
$$ \gamma \cdot \chi ^{n'} = e^{2 \pi i \langle \gamma , n'\rangle } \cdot \chi ^{n'},$$
for $\gamma \in \Gamma$ and $n' \in L';$ i.e., by:
$$\gamma u = u \text{ and } \gamma uv = \zeta uv, $$
where $u= \chi ^{f_1 '} \text{, } uv=\chi ^ {f_1 ' + f_2 '}$ and 
$\zeta = e^{2 \pi i \frac{1}{k}}$ is a primitiv $k\text{-th}$ root of unity.
We refer to \cite{F} for the notation. Let $\{ e_1 ,e_2 , e_3 \}$ be the 
standard basis of $\RR ^3 .$ 
 
 The fan of $\cY$ is the union of the two cones $\sigma _1  '$ 
and $\sigma _2  '$  generated by $\{ e_1 , e_1+ e_2 , e_3 \}$ and 
  $\{ e_1+e_2 , w' \}$ respectively, 
in the lattice $N' = e_1 \ZZ \oplus e_2 \ZZ \oplus e_3 \ZZ + w' \ZZ
 + (e_1 + e_2 ) \ZZ,$
  for some vector $w'$. 

We want to find 
$w'.$ Well, we know that the projection onto $e_2 \RR \oplus e_3 \RR$ 
along $e_1 + e_2$ is the fan of a toric \'etale neighborhood of $Q_1 '$ in
$X_1 ',$ where $X_1 ' = \cX ' \times _{C_1} C_1 '$ and 
$Q_1 ' =\cQ ' \cap X_1 '$ is the corresponding zero section.
 Let $\pi _{e_1 + e_2} : \RR ^3 \to e_2 \RR \oplus e_3 \RR \simeq \RR ^2$ 
be such projection. Then $\pi _{e_1 + e_2} (x,y,z) = (y-x,z) $ and 
$ \pi _{e_1 + e_2} (w') = (-1, n).$ 
But $w'$ must project onto the vector $e_2$ via the projection 
$\pi _{e_3} : \RR ^3 \to  e_1 \RR \oplus e_2 \RR \simeq \RR ^2,$ since the 
latter maps the cone of $\cY$ onto the cone of $V,$ as the ray $e_1 \RR _{+}$ 
represents the zero section of $\cY.$ 
 Hence $w'= e_2 + n e_3.$ 

In order to find the fan of $\cY / \Gamma$ as in 
\ref{one fibre} we need to identify a lifting in $N'$ of the sublattice $L$ 
of $L'.$ So $N = e_1 \ZZ \oplus e_2 \ZZ \oplus e_3 \ZZ + w' \ZZ
 + (e_1 + e_2 ) \ZZ,$ and we want to find $w.$ 
Let $ \pi _{e_3} : \RR ^3 \to e_1 \RR \oplus e_2 \RR$ be the projection.  
 Since the divisor corresponding to $e_3$ is the zero section $\cQ ',$ 
$ \pi _{e_3} (w)= \frac{1}{k} e_2.$  
 Also, the action of $\ZZ / k \ZZ $ on the divisor $\cS$ 
corresponding to $w$ (i.e., 
the pull-back to $\cY$ of the divisor in $\cX$ corresponding to 
$S \subset \cC,$ with the notation as in lemma \ref{one fibre} ), is trivial, 
therefore it must be trivial on $X_1 ' \cap \cY.$ Thus $\pi _{e_1 + e_2} (w)$ 
must lie on $w'=e_2 + n e_3;$ from this and from 
$\pi _{e_3} (w)= \frac{1}{k} e_2,$  we infer that $w=\frac{1}{k}(e_2 + n e_3).$
\end{proof}

\subsection{Two attaching fibres}

The analogous lemmas and proposition hold for the case of a chain of rational 
curves joining two curves in the central fibre:

Let $ \cC \to \Delta$ be  a family of nodal curves $\cC .$ 
Assume that the central fibre $ \cC _0 $ has 
a chain of rational components $C=\bigcup _{i=1}^N C_i$ meeting the rest of
 $\cC _0$ 
transversally 
in only two points $q_1 \text{ and } q_2 .$ The singularities of 
$\cC$ around $p_i=C_i \cap C_{i+1}$
 and the $q_i$ are at worst $A_{k_{i}-1}$ singularities, $i=1,...,N+1.$

We have:

\begin{lemma} \label{torbastwofib}

 Let $B$ be the rest of $\cC _0$ 
and $\cI \subset \cO _{\cC}$ the ideal sheaf of $C \cup B  .$
Then there is an \'etale neighborood $U$ of $C$ in $\cC ,$ a
 ~2-dimensional 
toric variety $Z$ and an \'etale map $t:U \to Z $ such that:

\begin{enumerate}
\item the fan of $Z$ is $F= \bigcup _i \sigma_i$ 
where $\sigma_1=\langle  f_1, f_1 + k_{1} f_2\rangle,$
$\sigma_i=\langle   f_1+\sum _{j=1} ^{i-1} k_{j}f_2, f_1 +
 \sum _{j=1} ^{i} k_{j} f_2\rangle $ and $\sigma_{r+1}= \langle f_1+
 \sum _{j=1} ^{r} k_{j}f_2, f_2\rangle$
\item the pull-back 
via $t$ of the ideal of $D_Z $ is the
 ideal $\cI \otimes \cO _U .$ 
\end{enumerate}
\end{lemma}

\begin{proof}
The proof is analogous to the one of lemma \ref{torbasonefib}, except for a
 few 
variations. 
 Let $W$ be a Zariski neighborhood of $C$ such that 
$W \cap (C \cap B) =
\{p , q\}.$

We can now contract $C$ to a rational double point of type 
$A_{r}$ where 
 $r+1=\sum _i k_i .$  
Let $\rho :\cC \to Y$  be the contraction map and let $W' = \rho (W)$ and 
$p' = \rho (C) .$
 Choose $U'$ to be an \'etale neighborhood of $W'$ that 
splits the node.
 Hence, there is an  isomorphism 
$\phi ^{\sharp}: R=k[[x,y,t]]/(xy-t^{r+1}) \to \hat {\cO} _ {W' , p'},$ 
such that the 
pull-back of the maximal ideal 
$\mathfrak m _{W' , p'} \subset \hat {\cO}_{W',p'} $ is the maximal 
ideal $(x,y,t) \subset R.$

The homomorphism of complete local rings $\phi ^{\sharp}$ produces an 
\'etale map $\phi :U' \to \Spec R$ for some \'etale neighborhood 
$U'$ of $p'$ in $Y .$ Note that the maximal ideal 
$\mathfrak m _{\cC , p}$ is generated by local equations of the branches 
of $B\cap U.$

Let $Z'$ be the toric 
variety whose 
fan $F'$ is the union of the cones 
$\sigma _1,$ ..., $\sigma _r ,$ $\sigma _{r+1}$ respectively generated 
$\{f_1, f_1 + f_2 \} ,$ ..., $\{f_1, f_1 + (r+1) f_2 \} ,$ 
$\{f_1+(r+1) f_2, f_2 \} $ in the lattice 
$N= f_1 \ZZ \oplus f_2 \ZZ .$ 

Let $U' \to Z'$ be the normalization of $U$ in the function field of $Z'.$
As in lemma \ref{torbasonefib} we construct a commutative diagram:

$$\begin{array}{cccccc}  U'&&\stackrel{t'}{\longrightarrow} &&Z'  \\ 
&\searrow \alpha&&&\\
b\downarrow &&W&&\downarrow \beta\\
& p_1 \swarrow  &&\searrow p_2&&\\
  U &&\stackrel{t}{\dashrightarrow} && Z
\end{array}$$

where $Z$ is the toric variety whose fan is $F =\bigcup \sigma_i$ 
where $\sigma_i=\langle  f_1+ \sum _{j=1} ^{i-1} k_{j}  f_2, 
f_1 +  \sum _{j=1} ^{i} k_{j}  f_2\rangle $ in 
$f_1 \ZZ \oplus  f_2 \ZZ,$ for $i\neq 1, r+1$ and where
$\sigma_1=\langle  f_1, f_1 + k_{1} f_2\rangle,$ and $\sigma_{r+1}= \langle f_1+
 \sum _{j=1} ^{r} k_{j}f_2, f_2\rangle .$

 An argument similar to the one given there, shows that the rational map 
$U \dasharrow Z$ is in fact regular.

Let $B '' = c' _* B$
 and let $\cI "$ be the ideal sheaf of $B .$
  By construction the pull-back of 
the ideal sheaf of $D_Z$ is the ideal sheaf 
$\cI \otimes \cO _U,$ by construction.

\end{proof}

As in the case of one attaching fibre, we want to be able to find a toric 
finite covering that ``untwists'' all the fibres of $\cX \to \cC$ 
above the $p_i$'s.
In the hypothesis of lemma \ref{torbastwofib}, let the possible 
$A$-singularity of $\cC$ at the ponts $p_i$ be of type $A_{k_i-1} .$  
 We have:

\begin{lemma} \label{toricbasexttwo}
 There is a $\bmu _n \times 
\bmu _n $-covering 
$f: V \to U$ ramified along $C$ and along the 
two branches 
of $B$ coming off $C$ 
 and an 
\'etale map 
to a ~2-dimensional toric variety $t' :V \to Z'$ such that:
\begin{enumerate}
\item the fan of $Z'$ is $F' = \sigma '_1 = \langle f_1', 
(f_1 '+ k_1 f_2)\rangle \cup ... \cup \sigma _i '=  \langle f_1' + 
\sum _{j=1} ^{i-1} k_{j}  f_2 ', f_1 ' +\sum _{j=1} ^{i} k_{j} 
 f_2 ' \rangle \cup ...\cup  
\sigma _{r+1} '= \langle  f_2',
f_1 ' + \sum _{j=1} ^{r+1} k_{j}  f_2 ', \rangle $ in the 
lattice $N' =f_1 ' \ZZ \oplus f_2 ' \ZZ.$
\item the pull-back via $t'$ of the toric ideal sheaf $D_Z $
is the ideal sheaf $\cI _{C ' \cup B ' \cup S'} \mid _V$
of $C ' \cup B ' \cup 
S'$ restricted on $V,$ where $B ' = f^* B;$ 
\item $f$ induces a toric morphism $f : Z' \to Z $ given by $f_1 ' = n f_1$ 
and $f_2 ' = n f_2 .$

\end{enumerate}

\end{lemma}
\begin{proof}

 We can apply 

lemma \ref{torbasonefib} and find a toric 
variety $Z$ and an \'etale morphism $t: Z \to U$ such that 
$t^*\cO_Z (\OO _{f_1} + \sum _{i=1} ^N \OO _{f_1 + k_if_2} + 
\OO _{f_2}) \simeq \cO _U (B+C).$

 Let $N'$ the 
sub-lattice of $N$ 
$N' =f_1 ' \ZZ \oplus f_2 ' \ZZ$ of $f_1  \ZZ \oplus f_2  \ZZ$ 
with 
$f_1'=nf_1$ and $f_2'=nf_2 $ and let $Z'$ be the toric variety given by 
the fan 
$F' =  \sigma _1 ' \cup \bigcup _i \sigma _i' \cup \sigma _{r+1} '$ 
union of the cones $\sigma '_1 = \langle f_1', 
(f_1 '+ k_1 f_2)\rangle , $ $\sigma _i '=  \langle f_1' + 
 \sum _{j=1} ^{i-1} k_{j}f _2 ', f_1 ' +\sum _{j=1} ^{i} k_{j} 
 f_2 ' \rangle $ and 
$\sigma _{r+1} '= \langle  f_2',
f_1 ' + \sum _{j=1} ^{r+1} k_{j} f_2 ', \rangle $ in the 
lattice $N' .$

We obtain the desired ramified finite covering $V \to U$ by
letting $V$ be the normalization of $U$ in the function field of $Z'.$

\end{proof}

What we have in mind is to look at those families of elliptic surfaces 
with section $(\cX \to \cC \to \Delta , \cQ \to \Delta)$ such that the 
central fibre contains a chain of rational components meeting the rest of the 
central fibre in only two fibres, and show that the picture is toric in this 
case too. This will allow us to perform the small contractions torically.
 
 So, let $X=\bigcup _{i=1} ^n X_i \subset {\mathcal X}_0$ be a 
chain of components of the 
central fibre, mapping down to a chain of rational curves $C=
\bigcup _{i=1} ^n C_i \subset 
{\mathcal C}_0.$ Assume that $C$ meets $B,$ the rest of $\cC _0 ,$ 
transversally 
in only two points, and corrispondingly $X$ is attached to the rest of 
$\cX _0$ along two fibres that we assume to be either stable or twisted. 
Let $Q=\bigcup _i Q_i= \cQ \mid _C$

 We have:

\begin{proposition} \label{many fibres}
 There is an \'etale  
neighborhood $\cU$ of $Q$ in $\cX$ and an \'etale morphism to a toric 
variety $T:\cU \to \cY $ such that $T^*\cO _{\cY} (D_{\cY})\simeq
 \cO _{\cU}(X + \cX \mid _B ).$   

\end{proposition} 

\begin{proof} As in proposition \ref{one fibre}, we can find a diagram:

$$\begin{array}{ccc} {\cV} & \stackrel{s}{\to} & U  \\ 
\pi ' {\dar} & &{\dar} \pi\\
  V &\stackrel{f}{\to}
 & U
\end{array}  $$

where $ \pi ':\cV \to V$ is a family of stable curves and the action of 
the Galois group $\Gamma \simeq \bmu _n \times \bmu _n$ of 
$f: V \to U$ extends to an action 
on $\cV$ such that $\Gamma \backslash \cV \simeq \cX.$ 
  Let $k_i$ with $i=2,...,N$   
and $k_1$ and $k_{N+1}$ be the orders of singularities around 
$ p_i:=C_i \cap C_{i+1}$ and 
$p_j:=C \cap B_j$ respectively ($i=2,...,N$ and $j=1,2$); and let $h_i$ 
and $n_j$ be such that the group acting non-trivially on the fibre above $p_i$
is $G_i \simeq  \ZZ / h_i\ZZ $ and on the one above $q_j$ is $\Gamma _j 
\simeq \ZZ /n_j \ZZ .$ 

According to lemmas \ref{torbastwofib} and \ref{toricbasexttwo} 
we can find toric varieties and \'etale morphisms:
$t:U \to Z $ and $t':U' \to Z'$ such that $f:Z' \to Z$ is toric. 
 
We can choose toric cyclic coverings:
 $$t_i :V_i \to U$$
 with Galois group $\bmu _{k_i}$ and that desingularize the 
singularities around points $p_i .$
 The covering $t:U' \to U $ factors through the normalization $W$ of the fibre 
product $V_1 \times _U ...\times _U V_{N} \times _U V_{N+1} .$
Let $\pi _i :W \to V_i$ be the standard projection.

 Since $W$ is smooth around those points $p_i '$ that map 
to $p_i $ via $t_i \circ \pi _i$ and the fibres of  $\cX _i =
{t_i \circ \pi _i}^* \cX$ 
above those points are stable. Hence the same is true for $\cV .$ 

Therefore on  $\cV$  we have an \'etale neighborhood $\cU '$ of 
the zero section $Q '$ in $\cV$ and an \'etale map $T':\cU' \to \cY '$ 
to a toric variety, when we take as $\cY '$ the total space of the 
normal bundle of $\cQ '$ in $\cX :$
$\cN =\cN _{\cQ '/ \cX '} \mid _{U'}.$ 
  
 Now the as in proposition \ref{one fibre}  the action of 
$\Gamma$ on $\cY '$ is linearizable, and therefore it acts as a subgroup of 
the torus. We take as $T: \cU \to \cY$ the quotient of $T' .$ 
The statement that $T^*\cO _{\cY} (D_{\cY})\simeq
 \cO _{\cU}(X + \cX \mid _B )$ follows from 
the analogous statement for t in lemma \ref{toricbasexttwo} , since we 
obviously have a commutative diagram:

$$\begin{array}{ccc} {\cY} & \stackrel{T}{\to} & \cU  \\ 
 {\dar} & &{\dar} \\
  Z &\stackrel{}{\to} &  U
\end{array}  $$

\end{proof}

 Let us keep the notation as in lemma \ref{many fibres} and set $Q_i = 
\cQ \mid _{X_i}.$ We have:

\begin{lemma} \label{toric two} The fan of the toric variety 
$\cU$ is the union of the cones $\sigma _1 \text{ ,}
 \sigma _2 \text{, } ... \sigma _N,$ respectively generated by:
 $\{ \frac{1}{n} e_1, w_1 =e_1+ k_1 e_2 ,e_3  \},$ ... 
$ \{ w_N = e_1 + (k_1 + ...+ k_N) e_2 + (n_1 k_2 +(n_1 +n_2) k_3+ ....+
(n_1  +...+n_{N-1}) k_N ) e_3 , 
    w_{N+1} = e_1 + (k_1 + ...+ k_{N+1}) e_2 + 
(n_1 k_2 + (n_1 +n_2) k_3+....+(n_1+...+n_N) k_{N+1})  e_3, e_3 \}$
in the lattice $L = \frac{1} {n} e_1 \ZZ \oplus e_3 \ZZ \oplus 
(\frac {1}{n} e_2 + 
\frac{a} {h_1} e_3) \ZZ ,$
where the $n_i  =Q_i ^2$ and 
$g.c.m(a,h_1)=1$ are the integers determining the actions on the first fibre 
(hence 
the action on all the other fibres is determined by this datum).

\end{lemma}

\begin{proof}
 
 Let $(x_1, y_1, t) \text{, }(x_2, y_2, t)$....$(x_{N+1} ,y_{N+1},t)$ 
be coordites around $q_1=B_1 \cap C_1$ $p_1= B_2 \cap C_2$ ....
$q_2=C_N \cap B_2,$ so that these a neighborhood of $p_i$ (resp. $q_j$) in 
 $\cC$ has equation $x_i y_i =t^k_i$ (resp.$x_j y_j =t^k_j$). As in 
lemma \ref{many fibres} we have a diagram:

$$\begin{array}{ccc} {\cY} & \stackrel{f}{\to} & \cZ \\ 
\pi ' {\dar} & &{\dar}\pi\\
  V &\stackrel{p}{\to}
 & U
\end{array}  $$

where all the morphisms and varieties are toric and 
${\cY} / \Gamma \simeq   \cZ,$ 
 where $\Gamma$ is the Galois group. 
  
 The fan of $U$ is given by the cones generated by 
$\{f_1 , v_1=f_1 + k_1 f_2 \}$, $\{v_1, v_2=f_1 + (k_1+k_2) f_2 \}$,...,
$\{v_N , v_{N+1}=f_1 + (k_1+...+k_{N+1}) f_2 \},$ in the lattice generated 
by $\{f_1, v_1, ...v_{N+1} \}$ where  $\{ f_1, f_2 \}$
 is the standard basis of $\RR ^2 .$ Hence the fan of $V$ is the union of the 
cones  $\langle f_1 ' , v_1 '=f_1 '+ k_1 f_2 '\rangle $, 
$\langle v_1 ', v_2 '=f_1 '+ (k_1+k_2) 
f_2 '\rangle $,..., $\langle v_N ' , v_{N+1} '=f_1  ' + 
(k_1+...+k_{N+1}) f_2 '\rangle ,$ 
where $f_1 ' =n f_1 $ and $f_2 ' = n f_2.$ The fan of
$\cY$ is the union  $\sigma _1 ' \cup  ... \cup \sigma _{N+1}'$ of the cones:
 $\sigma _1' =\langle e_1, (e_1 + k_1 e_2) ,e_3\rangle  $, 
$\sigma _2 ' =\langle e_1 + k_1 e_2, w_1 ' ,e_3\rangle  ,$
 ...$,\sigma _{N+1} '=\langle w_N ', w_{N+1} ', e_3\rangle  $ 
for some vectors $w_i ',$ in the lattice $L'$ generated by 
$\{e_1, e_1 + k_1 e_2, w_1 ' , ...,w_N ', 
w_{N+1} '\}.$ Our first goal is to find these vectors.

 Let $\pi _{e_3} : \RR ^3 \to \RR ^2$ be the projection along $e_3$ onto 
$e_1 \RR \oplus e_2 \RR $, and let $\pi _{w_i '} : \RR ^3 \to \RR ^2$ 
 be the projection along $w_i '$ onto $e_2 \RR \oplus e_3 \RR .$
 The zero section $\cQ$ is the toric divisor of $\cZ$ corresponding to the 
ray $e_3,$
 hence all the $w_i \text{'s}$ must project to $e_1 +(k_1 +..+k_i) e_2,$ 
since $\pi \mid _{\cQ} : \cQ \to U$ 
is an isomorphism (strictly speaking we should write $\cQ \mid _U$ here, 
but with 
abuse of notation we shall simply write $\cQ$ for the rest of the discussion).

 Analogously the surface $X_1$ corresponds to the ray $e_1 +k_1 e_2,$ 
 the image of the cones $\sigma _1$ and $\sigma _2$ in the lattice 
$L' / \langle e_1 +k_1 e_2\rangle $ should give the fan of a toric varity 
whose only 
complete toric curve 
(the zero section) has self intersection $-n_1 .$ Therefore, 
the vector $ \pi _{ e_1 +k_1 e_2}(w_1 ') =(k_2,z)$  must be proportional to 
$(1,n_1),$ 
since $\sigma _1$ maps to the cone given by the second quadrant,
and so $w_1 ' =(1, k_1 + k_2 , n k_2) .$ 

 Assuming by inductive hypothesis that 
$w_ N ' = (1, k_1 +..+k_N, n_1 k_2 +...+(n_1+...+n_{N-1}) k_N ), $ we want 
to show that 
 $w_ {N +1} ' = (1, k_1 +..+k_{N+1}, n_1 k_2 +...+(n_1+...+n_{N}) k_{N+1} ).$ 
Since $\pi _{e_3} (w_N ')=(1,k_1 +..+k_{N+1}),$ we just have to show that the 
last 
coordinate is the one claimed. Note that $\pi _{w_N '} 
( w_{N+1}')= (k_N,z- n_1k_1+...+(n_1 + ...+n_{N-1})k_N)$ and 
 and  $\pi _{w_N '} 
( w_{N-1})= (-k_{N-1},- (n_1+...+n_{N-2})k_{N-1})$
so $z=n_1 k_2 +...+(n_1+...+n_{N-1}) k_N .$ 
 All is left to do is to identify the action of $\Gamma.$  

 Obviously the vector 
$\frac {1} {n} e_1$ from the lattice of the base 
lifts to an element of the super-lattice $L$ of $L'.$ In fact the kernel of 
the map: 
 $$\ZZ / n \ZZ  \times \ZZ / n \ZZ \to \ZZ / k_1\ZZ $$ 
as in lemma \ref{many fibres}, must be contained in the cone generated by the 
two 
adjacent vectors  $e_1$  and $e_1 + k_1 e_2$, since it acts 
trivially on the fibres above the corresponding point.  Note that this kernel 
is 
generated by $\frac{1}{n} e_1$ and  $\frac{k_1}{n} e_2.$
  We want to know more, namely 
we want to determine the kernel of the map: 

 $$\ZZ / n \ZZ  \times \ZZ / n \ZZ \to \ZZ / h_1\ZZ .$$ 

and in doing so, lift the vector $\frac{1}{n} e_2$ from the lattice of the 
base curve $U.$  Let 
 $\pi _{e_3}$
 as  before, then we are looking for a vector $w$ such that 
$\pi _{e_3} (w) = \frac{1}{n} e_2$ and such that it represents the action of 
$\ZZ / 
h_1 \ZZ $ on the fibre. That is to say, we want:
$w=(0, \frac{1}{n}, z)$  with $z $ such that $h_1 z \in \ZZ.$ Hence
 $w=(0 , \frac {1}{n}, 
\frac{a} {h_1}) $ for some $a$ with $g.c.d(a,h_1)=1;$ the integer $a$ is 
completely 
determined by the action of $\bmu _{h_1} \simeq \ZZ / h_1 \ZZ$ on $\cX' .$
 This proves the lemma.

\end{proof}

\section{Toric log flips and toric log-canonical contractions}

\subsection{Log-flips}

In this section we want to show that we can perform the log-flips
torically.

We will need the following:

\begin{definition}

Given a pair $(X,D)$ consisting of a variety $X$ and a divisor
$D=D_1 + ...+D_n$ (where the $D_i$ are irreducible and reduced), we
say that a a pair $(T, D_T)$ consisting of a toric variety $T$ and its
toric divisor, is a {\bf toric \'etale neighborhood} of a subvariety
$Y \subset X$ if there is an \'etale neighborhood $u:U \to X$ of $Y$
in $X$ such that:

\begin{enumerate}

\item there is an \'etale map $t: U \to T ;$
\item $t^* \cO _T (D_T) \simeq \cO _U (D);$ 
\end{enumerate}

where $\cO _U (D ) := u^* \cO _X (D).$
\end{definition}

 Let $(\cX \to \cC \to \Delta, \cQ , )$ be a 
family of {\it strictly prestable}  
elliptic surfaces over a 
DVR scheme $\Delta .$ Assume that the special fibre 
$\cX _0 \to \cC _0$ contains a surface $X_1 \to C_1$ 
over a rational curve $C_1$ and a surface $X_2 \to C_2$ along a fibre 
$G_1$ which is either stable or a twisted curve. Let $Q_i = X_i \cap \cQ.$
 In what follows, we will say thar a triple $(\cY \to S \to \Delta ,
 \cD , \cY \to \bM _{1,1})$ is {\it locally stable} if the naturally
 induced triple $(\cY \to S \to \Delta ,\cD , \cY \to S \times \bM
 _{1,1})$ is stable.

We have:

\begin{theorem} \label{log-flip}
 The log-flip of $Q_1$ gives 
rise to a locally stable family $(\cX ^+ ,\cQ ^+) \to \cC ^+ \to \Delta ,$
or, which is the same the triple $ (\cX ^+ ,\cQ ^+, \cX ^+ \to \cC ^+ \times 
\bM _{1,1})$ is stable. 
Moreover the central fibre decomposes as a union of surfaces $X_1 ^+ \cup 
X_2 ^+$ where $X_1 ^+ $ is obtained from $X_1$ by contracting $Q_1$ and 
log-flipping produces a new rational curve $Q_1^+$ which meets 
$Q_2 ^+$ (the zero section of $X_2^+$) within the smooth 
locus of $X_2 ^+ .$ Furthermore, there exists {\it an \'etale toric
  neighborhood} $(T,D_T)$ whose fan $F^+ \subset N$ is given by (keeping the notation
as in lemma \ref{toric one fibre}) the union of the two cones: $\sigma
^+ _1 + \langle e_1 , w , e_3 \rangle$ and $\sigma ^+ _2 = \langle e_1
, e_1 + e_2 , w \rangle .$ This singularity is canonical, and in
particular $\cX _0 ^+$ is {\it semilog-canonical}.
\end{theorem}

\begin{proof}
First, we want to show that it is enough to perform the log-flip in an 
\'etale neighborhood of $Q $ in $\cX .$ 

In fact, if $\mathcal U$ is an \'etale 
neighborhood of $Q$ in $\cX ,$ let ${\mathcal U }^+$ be
the log-flip of $Q.$ 
It is clear that $\cX \setminus {\mathcal U}$ and  patch 
together to form an Artin algebraic (or analytic) space $\cX ^+ ,$   
so we just need to show that we can find an ample line bundle on it.

Let $S \subset \cC$ any effective horizontal divisor meeting $C_1$ 
transversally. After possibly an \'etale 
  base change $\Delta ' \to \Delta$ we can make sure that $S$ is a section 
around $C_1 ,$ so we may assume it is to begin with. 
So the bundle $\omega _{\cX / \Delta} \otimes \cO _{\cX}(\cQ) \otimes 
\cO _{\cX}(\cX \mid _S)$ is $\QQ \text{-Cartier}$ and contracts $Q,$ and is 
$\Delta \text{-ample}$ 
on $\cX ^+ ,$ since by hypothesis $Q$ is the only curve on which 
$\omega _{\cX / \Delta} \otimes \cO _{\cX} (\cQ)$ fails to be positive. 
 It therefore suffices to construct $\cX ^+$ as an algebraic space, because 
the fact that it possesses an ample divisor makes it automatically a scheme.

 From lemmas \ref{one fibre} and \ref{toric one fibre} we know 
that there is an \'etale neighborhood of $Q_1$ in $\cX$ that is toric and 
whose fan $F$ is the union of the two cones $\sigma _1=\langle e_1  
  e_1+e_2 , e_3 \rangle $ and $\sigma _2=
\langle e_3 ,e_1 +e_2 ,  w \rangle  $  
in the lattice 
$N := \bigoplus _{i=1,3}  e_i \ZZ \oplus w \ZZ ,$ where 
 $\{e_1 ,e_2 ,e_3 \}$ is the standard basis of $\RR ^3$
and $w =\frac{1}{k} (  e_2 + n e_3).$ 
 The log-flip is thus constructed by taking the fan 
$F^+=\sigma _1 ^+ \cup \sigma _2 ^+$ 
where  $ \sigma _1 ^+ =\langle e_1  ,w , e_3\rangle  $ and 
$\sigma _2 ^+ =\langle e_1, e_1 +e_2 , w \rangle $  
(see \cite{R1} Theorem 2.4), and now 
the curve $Q^+$ corresponds to the face generated by $\{e_1 , w  ,w- e_1 \}$ 
and $X_2 ^+$ to the vector $e_1 +e_2.$ 

The matrix $A$ associated to the 
vectors $\{ e_1 , e_3 , w \}$ of the lattice $N$ is :
$$ \bordermatrix {&&& \cr
& 1 &0 &0 \cr
 & 0 &0 &1 \cr
& 0 &\frac{1}{k} &\frac{n}{k}} .$$

Since $det(A)=- \frac{1}{k},$ (i.e., $A\in SL_3 ({\frac {1}{k}}
\ZZ )$) these vectors  form a lattice basis for $N,$ 
 hence the toric variety corresponding to the cone $\sigma _1 ^+$ is smooth 
(see \cite{F} page 29). The rest of the proposition is obvious.
\end{proof}

 In particular, what this says is that if we have a chain of 
 rational curves $C_1 \cup...\cup C_n \subset \cC _0$ above which the
 zero-section is an extremal ray, we can contract them one by one and 
perform the log-flips above them inductively.

\begin{remark} We can "straighten up" the lattice $N$ in which the fans $F$ 
and $F^+$ 
live. Indeed, by applying the transformation: 

$$ \bordermatrix {&&& \cr
& 1 &0 &0 \cr
 & 0 &k &0 \cr
& 0 &-n &1 } ,$$
we can send $N$ to the lattice $f_1 \ZZ \oplus f_2 \ZZ \oplus f_3 \ZZ, $ and 
the fan 
$F$ becomes $\langle f_1,f_3 ,  (f_1+ kf_2 -nf_3)\rangle \cup \langle
(f_1 + kf_2 -nf_3),
 f_2 , f_3 \rangle $
and the fan $F^+$ becomes $\langle f_1,f_2, f_3\rangle  \cup 
\langle (f_1 + kf_2 -nf_3), f_1 ,  f_2 \rangle .$
 In particular the singular point of the threefold $\cX ^+$ 
through which $Q_1 ^+$
 goes is an $\frac{1}{n} (1,k, 1)$ threefold singularity.
\end{remark}

We need to understand what the different surfaces after the log-fip look like, 
in particular what their singularities are. We have:

\begin{lemma} With the notation as in theorem \ref{log-flip}.
Let $n'$ and $a$ be the unique integer such that $nn' \equiv 1 
\text{(mod k)}$ and 
$0 \leq ak- n' <n .$ Then  the singularity of $X_2 ^+$ at 
$Q_1 ^+ \cap G^+$ is an 
$A_{k,m}$ singularity, where $m=ak-n' .$ 
\end{lemma}
\begin{proof} 
Straightforward, using the description of the fan $F^+$ given in the remark 
above, 
since the surface $X_2 ^+$ corresponds to the ray $f_1 \RR _{\geq 0} ,$
 that is to say $F^+ = \langle e_1 , e_2 \rangle \cup \langle e_2 , 
-ke_1 + n e_2 \rangle .$ 

\end{proof}
and:
\begin{lemma}
Let $k'$ be the unique integer such that $kk' \equiv 1 \text{ (mod n)}$
and $0 \leq an- k' <n .$
Then the 
singularity of $X_1 ^+$ at the point to which $Q_1$ gets contracted is an 
$A_{n,k'} .$ 
singularity
\end{lemma}

\begin{proof} 
For this is more convinient to look at the description of $F^+$ in theorem 
\ref{log-flip}.
 The surface $X_1 ^+$ corresponds to the ray $\langle e_1 + e_2\rangle .$
Let $\pi : \RR ^3 \to \RR ^3$ be the projection along that ray onto 
$e_2 \RR \oplus e_3 \RR .$ Then $\pi (w) = (\frac{1}{k}, \frac{n}{k}) ,$
and $\pi (e_1)=
(-1,0) .$   The cone $\langle \pi(w) , \pi(e_1)\rangle $ in the 
lattice 
$e_3 \ZZ \oplus \pi(w) \ZZ$ is equivalent, via the transformation 
$\pi(w) \to f_1$ 
and $e_3 \to f_2 ,$ to the cone $\langle f_1, (nf_2 -kf_1)\rangle ,$ 
in the fan $f_1 \ZZ \oplus f_2 \ZZ .$ 
The former gives rise to an affine toric variety 
isomorphic to the one in the statement (sending $f_1$ to $f_2$ and viceversa). 
\end{proof}

We will also need:

\begin{lemma} \label{selfintinv}

Keeping the notations as above, one has that: 
$$ {Q_1 ^+ }^2 = \frac{1}{Q_1 ^2} .$$

\end{lemma}

\begin{proof}
The fan of a toric neighborhood of $Q_1 $ in $X_1 $ is given by:
$$
F = \langle e_1 , e_2 \rangle \cup \langle e_2 , -k e_1 + n e_2 \rangle 
$$ 

in the lattice $N = e_1 \ZZ \oplus e_2 \ZZ ;$ the one of a toric neighborhood 
of $Q_1 ^+$ in $X_2 ^+$ is given by: 

$$
F ^+= 
\langle e_1 , e_2 \rangle \cup \langle e_2 , k e_1 - n e_2 \rangle .
$$ 

Hence we have that $Q_1 ^2 = -\frac {n}{k}$ and ${Q_1 ^+} ^2 = -\frac {k}{n}.$

\end{proof}

When performing a log-flip according to theorem \ref{log-flip}, even if the 
surface $X_2$ were {\it standard } to begin with, the surface $X_2 ^+$ in 
general is not. In fact it is some toric blow-up of one. 
We are therefore naturally led to the following:

\begin{definition}\label{logstand}

A {\bf log-standard elliptic surface} $(Y \to C,Q,G+F)$ is a triple consisting 
of an elliptic surface 
$Y\to C, $ its zero section $Q$ and a marking of $s$ curves 
$G=\bigcup G_i$ and and $N$ of fibres $F=\cup F_j$ with 
a regular birational map $g:Y \to Y',$ called {\bf structure morphism},
 to a {\em standard } 
elliptic 
surface $Y'\to C$ with zero section $Q'$ such that:
\begin{enumerate}
\item the exceptional divisor $E=\bigcup E_i$ is the disjoint union of  
smooth and irreducible rational curves $E_i$  meeting the $G_i$ 
transversally at one point,
and $g(E_i)=p_i =g(G_i) \cap Q \in Y' ;$
\item for each $i$ and $j$ the curves 
$G_i ':=g( G_i)$ and $F_j '= g( F_j)$ are either stable or a 
twisted fibres of $Y' $ and $G_i $ is the proper transform of $G_i '.$
\item for each $i$ there exist \'etale neighborhoods $U \to Y$ of 
$E_i \subset Y$ and $U' \to Y'$ of $p_i \in Y'$ and 
  morphisms to toric varieties 
$t:U\to T$ and $t' : U' \to T';$
\item the fan of $T$ is the union of the two cones $\langle e_1 , e_2 \rangle 
\cup \langle e_1 , ke_1 -n e_2 \rangle $ in the lattice $N= e_1 \ZZ \oplus 
e_2 \ZZ , $ and $t^* \cO (D_T ) \simeq \cO_U ( E_i + G_i + Q );$ 
\item the fan of $T'$ is the cone $\langle e_2 , ke_1 -n e_2 \rangle $ 
in $N,$ and ${t'}^* \cO (D_T') \simeq \cO_{U'} ( G_i' + Q ');$
\item the morphims $U \to U'$ is induced by the toric blow-up $T \to T'$ 
induced by the subdivision $\langle e_1 , e_2 \rangle 
\cup \langle e_1 , ke_1 -n e_2 \rangle $ of 
$\langle e_2 , ke_1 -n e_2 \rangle. $

\end{enumerate}   

Here $D_Z$ is the toric divisor of a toric variety $Z.$  

 We shall call a divisor like  $G_i$ a {\bf splice} and one like $F_i$
 a {\bf scion}.
\end{definition}

\begin{picture}(250,250)(0,0)
\put (180,30){\makebox(100,200) {\includegraphics{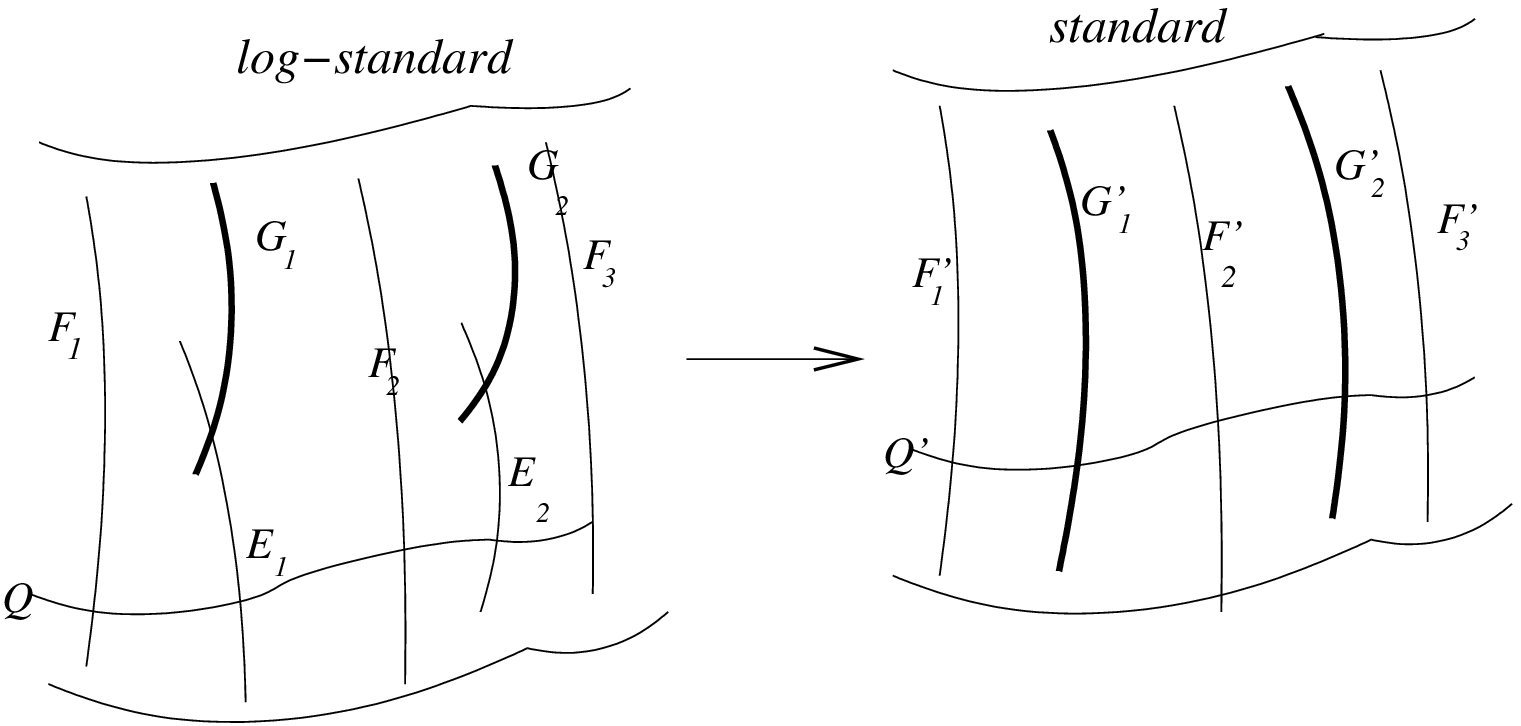}}}
\put(150,20){{fig.6}}
\end{picture}

Something analogous happens to the surface $X_1 ^+ ,$ except that we now loose 
the fibration structure. 
Let us keep the same notation as above.
We have:

\begin{definition}\label{I}

A {\bf ~N-pseudoelliptic surface} is a pair $(Y,G+F)$ 
consisting of a surface $Y,$ $N$ marked curves $F=\bigcup F_i$ 
and $s$ marked curves $G= \bigcup G_i$
with a regular birational 
map $g:Y' \to Y,$ called {\bf structure morphism}, 
from a {\em log-standard }elliptic surface 
$(\pi ' :Y' \to \PP ^1, Q' , G'+F')$ 
such that:
\begin{enumerate}

\item The proper transform of the $F_i$ are fibres of $Y' \to \PP ^1 ;$

\item the exceptional divisor of $g$ is the zero-section $Q$ of $Y'$
and $g$ maps $Q$ to a point $p;$

\item there exist a \'etale neighborhoods $V\to Y$ of $p \in Y$ and 
$V'\to Y'$ of $Q'\subset Y'$ and morphisms 
to toric varieties $\tau : V \to Z$ and $\tau ': V' \to Z'$
and a toric morphism $b: Z' \to Z$ such that $ b \circ \tau ' = \tau \circ g ;$

\item   the fan of $Z'$ is the union of the two cones 
$\langle e_1 , e_2 \rangle 
\cup \langle e_1 , ke_1 -n e_2 \rangle $ in the lattice $N$
and $\tau ^* \cO (D_Z ) \simeq \cO_U ( F_1' +S '+ Q ')$
for some fibre $S'$ of $Y' \to \PP ^1 ;$ 

\item the fan of $Z$ is the cone 
$\langle e_2 ,   ke_1 -n e_2 \rangle $ in the lattice $N$
and $\tau ^* \cO (D_{Z'} ) \simeq \cO_U ( F_1 +S )$
where $S=g_* S' ;$

\item the toric morphism $b: Z' \to Z$ is induced by the subdivision 
$\langle e_1 , e_2 \rangle 
\cup \langle e_1 , ke_1 -n e_2 \rangle $ of 
$\langle e_2 ,   ke_1 -n e_2 \rangle .$ 

\end{enumerate}

We call such a $Y$ {\bf isotrivial }if $Y' \to \PP ^1$ is isotrivial.
Furthermore, we call {\bf pseudolleptic surface of type I}
a {\em ~1-pseudoellitpic surface } $(Y,G+F_1);$ 
the component $G_i$ is still called a {\bf splice} and $F_i$ is still
called a {\bf scion}.

\end{definition}

\begin{remark}

An $n$-pseudoelliptic surface is {\it log-canonical} if and only if $n \leq 2 .$
\end{remark}


\begin{picture}(250,250)(0,0)
\put (180,30){\makebox(100,200) {\includegraphics{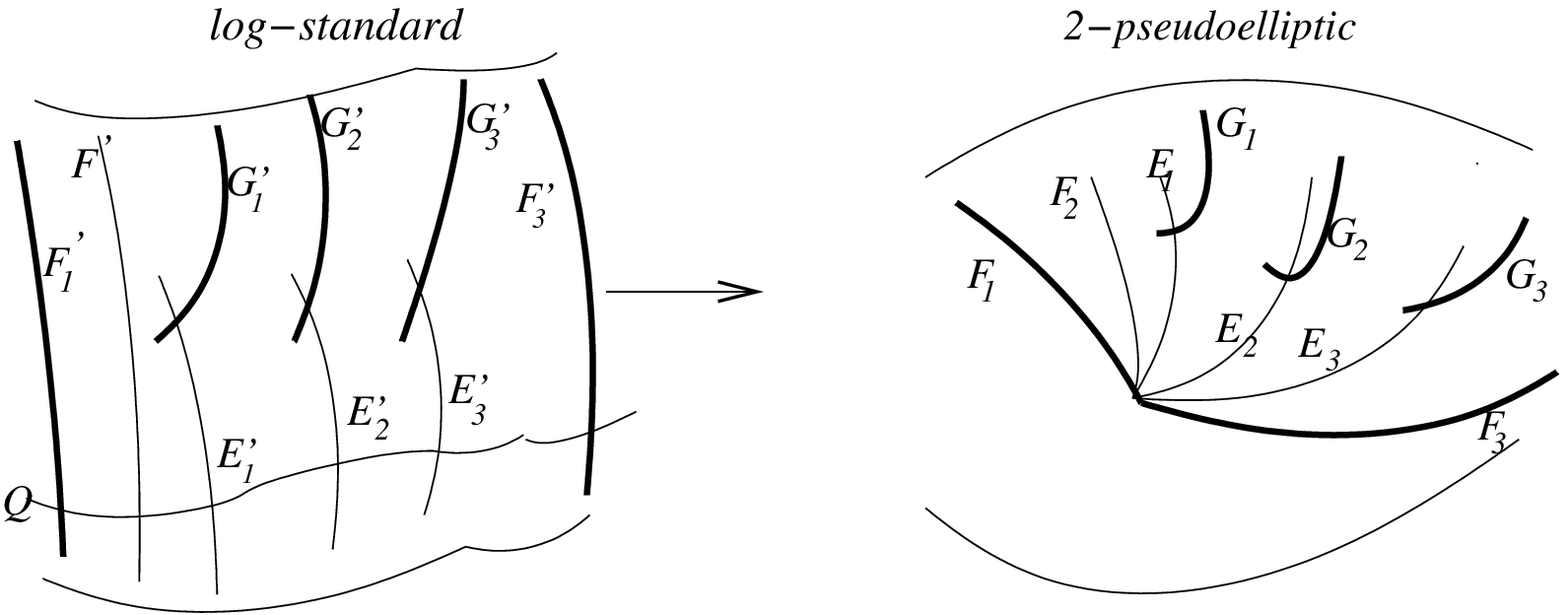}}}
\put(150,20){{fig.7}}
\end{picture}

 After performing any number of flips, we want to make sure that we know what happens to 
$X_1 ^+$ and $X_2 ^+$ in particular that we know when 
the restriction of the log-canonical divisor is nef and big on them.
This is taken care by the following  proposition:

\begin{proposition} \label{ample}

\begin{enumerate} 

\item Let $Y$ be a log-standard surface, $g: Y \to Y'$ its structure
  morphism with $Y'$ standard. Let 

$$ L_{Y'} = K_{Y'} + \sum _{i=1} ^s G_i ' + \sum _{j=1}^r F_j ' +Q'$$
and

$$L_Y = K_Y + \sum _{i=1} ^s G_i  + \sum _{j=1}^r F_j+Q ; $$
then:
$$g ^* L_{Y '} =L_{Y} + \sum _{i=1}^s E_i .$$

\item
Let $Y$ be an $n$-pseudoelliptic surface with structure morphism
$g : Y' \to Y,$ with $Y'$ is log-standard. 
Let:

$$ L_{Y'} = K_{Y'} + \sum _{i=1}^s G_i ' + \sum _{j=1}^r F_j ' +Q$$
and

$$L_Y = K_Y + \sum _{i=1}^s G_i  + \sum _{j=1}^r F_j; $$
then:
 $$g ^* L_{Y } =L_{Y'} + \frac{2-n}{Q ^2} Q .$$

\end{enumerate}
\end{proposition}

\begin{proof}
We start with proving  part (1).

We have that:

$$ g ^* L_{Y' } =L_{Y} +  \sum _{i=1} ^s a_i E_i$$

\noindent
given that $\sum _{i=1} ^s E_i$ is the exceptional divisor of $g.$
 By an easy inductive argument (i.e., by performing the toric blow-ups
 one at the time), one can easily convince oneself that infact all the
 $a_i's$ must be equal to each other; let us indicate this rational
 number by $a.$

Let $\pi :Y\to C$ be the projection to the base curve of $Y,$ and let
$g$ be the genus of $C.$ Furthermore, let $k$ be the least common
multiple of all the orders of monodromy around the fibres $F_i .$

 We can take a base change of oder $k$ as in proposition \ref{genlogcan}
to untwist the possible twisted fibres and get a diagram:

$$\begin{array}{ccc} S &\stackrel{f} {\rightarrow} & Y \\
 \dar & &\dar \\
B &\stackrel{\phi} {\rightarrow} &C
\end{array}$$

\noindent
 Now the zero section $\overline Q$ of the surface $S$ is entirely
 contained in the smooth locus of $S,$ and therefore, if we indicate
 respectively by ${\overline F}_i$ and ${\overline G} _j $ the proper
 transforms of $F_i$ and $G_j ,$ by a
 computation analogous to the one in proposition \ref{genlogcan}, we
 get:

$$L_Y \cdot Q= \frac {1}{k} (K_S +\sum _i {\overline F} _i + Q + k
\sum _j {\overline G} _j)\cdot {\overline Q}= 2g-2+s .$$

\noindent
Since according to proposition \ref{genlogcan}, we have: 
$$L_{Y'} \cdot Q'= 2g-2+s+r ,$$ 
we can infer that $a=1 .$ This completes the proof of part (1).

As for part (2), we can again write:

$$ g^* L_{Y} = g ^* (K_{Y} +\sum _{i=1}^s G_i + \sum _ {j=1} ^r F_j
 )  =(K_{Y'} + \sum _{i=1}^s   G_i '+  \sum _ {j=1} ^r F_j  ' +Q ')+ a Q=
=L_{Y '} + a Q .$$

\noindent
But $ g ^* L_{Y} \cdot Q  =0$ by the projection formula and 
$L_{Y '} \cdot Q  =n-2 $ by lemma \ref{genlogcan}, therefore: 

$$ n-2 + a Q ^2 = 0 .$$
This proves part (2).

\end{proof}

We need to say something about the positivity of the 
log-canonical bundle of a log-canonical surfaces and of
n-pseudoelliptic ones. We have:

\begin{proposition} \label{flipextremalray}

\begin{enumerate}

\item Let $Y\to C$ be log-standard with base curve $C$ of genus $g,$
  with $r$ $F_j$ 's and $n$ $G_i$'s. Moreover, assume that $ -1<E_i
  ^2<0 .$ Then the following hold:

\begin{enumerate}

\item if $2g-2+r >0,$ then $L_Y$ is ample;

\item if $2g-2+r=0 ,$ then $L_Y$ is semiample, and for any irreducible
  curve $D:$
$$L_Y \cdot D=0 \text{ iff }  D=Q ;$$

\item if $2g-2+r <0$ then $[Q] \in {\overline {NE}} (Y)$ is an
  extremal ray for $L_Y .$
 
\end{enumerate} 

\item let $Y$ be $n-$pseudoelliptic, then:

\begin{enumerate}

\item if $n\geq 2$ then $L_Y$ is ample;

\item if $n=1$ and $Q ^2 < -1 ,$ then $L_Y$ is ample

\item if $n=0,$ then $L_Y$ is ample if and only if $Q^2 <-4 .$

\end{enumerate}

\end{enumerate}

\end{proposition}

\begin{proof}

We will start proving part (1).

{\bf Proof of part (1)}
Recall from Proposition \ref{ample}, part (1), with the same notation
as therein, that:

$$g^*L_{Y'}=L_Y + \sum _{i=1}^s E_i .$$

\noindent
Therefore, in order to show that $L_{Y}^2 >0,$ it suffices to show
that $L_{Y'} ^2 >0$ and $2L_Y \cdot \sum _i E_i + \sum _i E_i ^2>0 .$
The first  inequality is a
consequence of corollary \ref{twistedampleness}, and:
$$2L_{Y} \cdot \sum _{i=1}^s E_i + \sum _{i=1}^s E_i ^2 = - \sum
_{i=1}^s E_i ^2 >0,$$

\noindent
where we have used that $L_Y \cdot \sum _{i=1}^s E_i  =  g^*L_{Y'}
\cdot  \sum _{i=1}^s E_i +
(\sum _{i=1}^s E_i)\cdot (\sum _{i=1}^s E_i) .$

Now,if $D$ is any irreducible curve on $Y$ other than one of the
$G_i$'s or one of the $E_i$'s, then $L_Y \cdot D= L_{Y'}
\cdot g(D)$ which is positive according to corollary
\ref{twistedampleness}. In case the irreducible curve is $E_i,$ then
we have that:
$$ L_Y \cdot E_i = -E_i ^2 >0 ;$$   

\noindent 
and in case it is $G_i,$ then, since $-G_i ^2 = G_i \cdot E_i,$ we
have:

$$  L_Y \cdot G_i = L_{Y'} \cdot G_i ' - E_i \cdot G_i= Q' \cdot G_i '
- E_i \cdot G_i ,$$

and since $G_i ' \cdot Q' =\frac{1}{-E_i ^2} E_i \cdot G_i$ (this is
obtained by writing $g^* Q ' = G_i + \sum a_i E_i$ and intersecting with
$E_i$ to obtain $a_i$ and then intersect $g^*Q '$ with $G_i$), we have that:

$$  L_Y \cdot G_i= (\frac{1}{-E_i ^2}-1) E_i \cdot G_i >0$$
since by lemma \ref{selfintinv}:
 $$-1 <E_i ^2 <0 .$$

\noindent
We can then conclude part (1) with the aid of corollary \ref{twistedampleness}.

We now prove part (2).

{\bf Proof of part (2)}
Recall from Proposition \ref{ample}, part (2)
$$g^* L_Y  =L_{Y'} + \frac {2-n}{Q^2} Q .$$

This implies that $L_Y \cdot G_i>0.$ 
If $n\geq 2$ $L_{Y'}$ is ample for the previous part, and $ \frac
{2-n}{Q^2} Q$ is effective (or $0$ if $n=2$), so part (a) is proved.

If $n=1,$ 

$$g^* L_{Y} = L _{Y'} +\frac {1}{Q^2} Q .$$

\noindent
Let $g'':Y' \to Y''$ the structure morphism of the log-standard
surface $Y' .$ Then, if we set:
$$ M:= K_{Y''} + \sum G_i '' +F'' +(1+ \frac{1}{Q^2}) Q'', $$
where $Q'' =g'' (Q')$, $F'' =g''(F')$ and ${G_i }'' = g'' ({G_i} '),$ then:

$$ g^* L_Y = {g''}^* M +\sum E_i$$ 

\noindent
and again we can conclude by means of part (1) and corollary \ref{twistedampleness}, given that  
$0< a:= (1+\frac {1}{Q^2})<1  $ and that $\lambda >1 .$ 

Let us finally analyze the case in which  $n=0.$ As before, we have: 

$$ g^* L_{Y} ={g''}^* M +\sum E_i . $$

where now:

$$ M:= K_{Y''} +  \sum G_i '' +(1+ \frac{2}{Q^2}) Q'', $$
\noindent
Therefore one can once again conclude, with the aid of corollary
\ref{twistedampleness} (given that $Q^2 < -4 $) that $g^* L_{Y}$ is
positive on every curve that is different from $Q$ and from one of the
$E_i$'s. Indeed, every other curve, except for the $G_i$'s (and for
these the computation is exactly like in part (1) of this
proposition),  will not meet the $E_i$'s. Also, since:
$$L_{Y} ^2 = M^2 +\sum _{i=1} ^s E_i ^2$$
from the very same corollary, we get that $L_Y ^2$ is positive as soon
as $Q^2 <-4 ,$ in fact, we have:
$$ L_Y ^2 = 2 a (s-2) + a (2-a) \lambda + \sum _{i=1} ^s E_i ^2 >0$$
\noindent
with $a:= (1+ \frac{2} {Q^2})$ and $\lambda = - Q^2,$ and therefore, since $E_i ^2 >-1,$ as
soon as:

$$2 a (s-2) + a (2-a) \lambda -s >0$$
and this occurs as soon as $Q^2 < -4$ as claimed.

So we only need to check that $L_Y \cdot E_i >0 ,$ but this is clearly
true, due to the fact that $E_i ^2 <0$ and $M\cdot E_i =0 .$

This concludes the proof, again thanks to Nakai-Mosheizon.

\end{proof}

We also want to know how $X_2 ^+$ has changed. For one thing we know that 
now it  must be a {\it log-standard} elliptic surface, even if it were 
only {\it standard } to begin with.

We want to know whether $X_2 ^+$ is {\it strictly pre-stable} when $X_2 $ is.
We have:

\begin{proposition} \label{self-int2} 
Mantaining the notation as in proposition \ref{flipextremalray} part (1),  
 the self intersection $Q^2$ of $Q$ in $Y$ is:
$$Q ^2 = {Q'} ^2 + \frac{1}{E_i ^2}  $$
 where  ${Q'} ^2$ 
is taken in $Y'.$ In particular, if $Y'$ is a strictly prestable
standard surface (i.e., $ {Q'} ^2 <-1$) and $-1<E_i ^2 <0$ for each $i,$
then $Q^2 <-1$  as well.
\end{proposition}

\begin{proof}
Write $g^* Q= Q' + \sum a_i E_i ,$ and intersect with $E_i ,$ to see
that it must be that $a_i = \frac{1}{-E_i ^2}$ (since $Q\cdot E_i
=1$).
 This yields:
$$ {Q'} ^2 = Q \cdot b^* Q' = Q^2 + \frac{1}{-E_i ^2}  $$ 

\noindent
which concludes the proof.

\end{proof}

We can therefore make the following:

\begin{definition}
We call a {\em triple } $(X \to C, Q, G+F)$ consisting of a 
log-standard elliptic surface with a given marking {\bf strictly 
prestable} if:
\begin{enumerate} 

\item $(Q \mid _{X \mid _B})^2 <-1 $ for each rational component $B \subset C;$
\item if $g: Y:= X \mid _B \to Y'$ is the structure morphism, then for
  each irredecubile component $E_i$ of the exceptional divisor $E$ of
  $g$ we have that $-1<E_i ^2 < 0 .$ 

\end{enumerate}
Analogously we call a {\em quadruple}  $(X \to C, Q+F, G, f:X \to C 
\stackrel {j}{\to}  \bM _{1,1})$ consisting of a log-standard elliptic 
surface with marking and map to moduli {\bf strictly prestable} if
the same conditions are asked only of those rational components $B \subset C$ 
for which $j\mid _B : B \to \bM _{1,1}$ is constant.
\end{definition}

\subsection{Small log-canonical contractions}

Let $(\pi: \cX \to \cC \to \Delta , \cQ)$ be a family of {\it strictly 
prestable} log-standard elliptic surfaces  
with section over a DVR scheme (or a polydisk, if one favour the
analytic flavor)) $\Delta,$ such that the special fibre 
$\cX _0 $ contains a chain of surfaces $X:=X_1 \cup ...\cup X_N$ attached
transvesally along stable or twisted fibres, with base curve a chain
or rational curves $C:= C_1 \cup ...\cup C_N .$ Furthermore assume
that $X$ is attached transversally to $\overline { \cX _0 \setminus X}$ along one
(twisted or stable) fibre each end. Let $Z_1 \to B_1$ and $Z_2 \to
B_2$ be the adjacent surfaces (that is to say the irreducible
components of $\overline { \cX _0 \setminus X}$ that meet $X$ at each end).

Let $\cC'  $ be obtained from $\cC $ by contracting the curve $C,$
$\rho : \cC \to \cC '$ the contraction map and $p= \rho (C).$

Let also $w_1 =e_1+ k_1 e_2 \text{, ..., } w_i = e_1 + (k_1 + ...+
k_i) e_2  + (n_1+...+n_{i-1}) k_i  e_3) \text{, ..., } w_{r+1} = e_1 +
(k_1 + ...+ k_{r+1}) e_2 + (n_1+...+n_{r}) k_{r+1}  e_3 ,$ and $L$ be
the lattice $L = \frac{1} {n} e_1 \ZZ \oplus e_3 \ZZ \oplus (\frac {1}{n} 
e_2 + 
\frac{a} {h_1} e_3)\ZZ .$

According to lemma \ref{extr ray} $\cQ \mid _C$ must be contracted by the 
log-canonical map, 
and in fact the contraction is described by the following:

\begin{theorem} \label{logcancont}

 The stable model of $\cX$ is a family $\cX ^c \stackrel {\pi '} {\to}
 \cC ' \to \Delta$ of surfaces 
such that the generic fibre and $\overline{\cX^c  _0 \setminus X}$
have not changed and $X^c = {\pi '} ^{-1} (p)$ 
 is attached to the rest of the central fibre along marked curves. The
 singularity of this point in $\cX ^c$ is toric, and a toric
 neighborhood of $p$ in $(\cX ^c , Z_1 ^c+ Z_2 ^c)$ (where $Z_i ^c:= \cX ^c
 \mid _ {\rho (B_i)}$) is given by the cone:
$\sigma =\langle \frac{1}{n} e_1 , w_1 \text{, ... } 
 w_r 
 w_{N+1} , e_3 \rangle$
in the lattice $L  .$ This singularity is canonical, and in particular
$\cX ^c _0$ is {\it semi-logcanonical}.

\end{theorem}

\begin{proof}As in the proof of theorem \ref{log-flip}
 we can reduce ourselves to finding the log-canonical 
contraction on an \'etale neighborhood, since by hypothesis 
$\omega _{\cX / \Delta} (\cQ)$ is $\Delta \text{-ample}$ except for 
contracting $Q .$
The rest is a direct consequence of \cite{R1}, proposition \ref{many
  fibres} and lemma \ref{toric two}. 
\end{proof}

 The normalization ${X^c}^{\nu}$ of $X^c = X_i ^c$ consists of a union of 
 normal surfaces $(Y_i ^c, {G_1^c} ^i , {G_2^c} ^i)$ with marked double 
curves ${G_1^c} ^i$ and ${G_2^c} ^i,$
and if $(Y_i ^c, {G_1^c} ^i , {G_2^c} ^i)$ denote components with 
double curves of the normalization of $X,$ we have a map $Y_i \to Y_i ^c$ 
which contracts the zero section to a point $p' \in {X^c}^{\nu}.$

\begin{lemma}\label{fantypeII}

There is an \'etale neighborhood U of $p' \in Y_i ^c $ with a map to a 
toric variety $g:U \to Z$ such that:
\begin{enumerate}
\item the fan of $Z$ is the cone $\Delta = \langle e_1 , w =-e_1+ ne_2 \rangle$
in the lattice $N = ({\frac {1}{n}} e_1 + \frac {a}{h_1}e_2) \ZZ \oplus 
 (\frac{1}{n}e_1 + \frac{n-a}{h_1}e_2)\ZZ,$ where $a$ and $h_1$ are completely 
determined by the monodromy of the action around one of the two ``marked''
fibres $G_1$ and $G_2 .$
\item  $g^* \cO _Z (D_Z) \simeq \cO _U (Q +G_1 +G_2),$ where 
$D_Z$ is the toric divisor.
\end{enumerate} 
 
\end{lemma}

\begin{proof}

The lemma is a consequence of the following two observations. 
On the one hand that the projection along 
$e_3$ onto $e_1$ must determine the action on $Q\simeq \PP ^1,$ 
which is given by choosing the lattice $\frac{1}{n}e_1.$ 

On the other hand, let $h_1$ be the integer such that $\bmu _{h_1}$ is the 
subgroup of $\bmu _n$ acting nontrivially on $G_1.$ The kernel of the map
$$\bmu _n \times \bmu _n \to \bmu _{h_1}$$ 
is determined by a vector $v=(\frac{1}{n}, z)$ such that $h_1 z=a \in \ZZ.$
 Thus the claim.

\end{proof}

What one has to prove now is that the log-canonical divisor on such 
components, which originally was only nef and big, has become ample after 
contracting the zero-section. In fact:

\begin{proposition} \label{ample2}
Let $X ^c $ be a $2-$pseudoelliptic surface,
$\alpha : X \to X ^c ,$ its structure morphism, and assume that $X \to
C$ is a {\it strictly stable} log-standard surface ($C$ rationl), with
zero section $Q.$
Then: $\alpha ^* L_{X ^c } =L_{X} .$ 
 In particular, $L_{X ^c }$ is ample.

\end{proposition} 

\begin{proof}

Write:
$$\alpha ^* L_{X ^c } =L_{X} + a Q ,$$ 

and intersect with $Q.$  
 By the projection formula, $\alpha ^* L_{X ^c } \cdot Q = 0$ and 
 $L_{X} \cdot Q =0$ according to lemma \ref{extr ray}, therefore 
the conclusion, on accounts of proposition \ref{flipextremalray} part (2).
\end{proof}

\begin{definition} \label{II}
 We call {\bf pseudoelliptic surface of type II}
 a {\em pair}  $(Y ^c, G^c +F_1 ^c +F_1 ^c )$  
consisting of a {\em ~2-pseudoelliptic surface} as in definition 
\ref{I} with the extra conditions that $S$ in point $5)$ is $F_2 ^c ,$ 
and that we replace the fan of point $5)$ 
with the fan of lemma \ref {fantypeII}
 
Such a triple is called {\bf isotrivial} if the surface $Y\to Y^c$ as in lemma 
\ref{fantypeII} is isotrivial.
\end{definition}

\section {The Stable reduction theorems}

\subsection{Stable reduction of triples}

We are now ready to prove the stable reduction theorem in the relative case 
of elliptic surfaces with sections and endowed with a regular map to 
$\bM _{1,1}.$

\begin{theorem} \label{stabred}
 Let ${\cX}_\eta\to {\cC}_\eta\to \eta$ be a stable 
elliptic surface over a smooth curve ${\cC}_\eta .$
 Then there is a finite extension of discrete valuation rings
$R\subset R'$ and a triple 
$( \cX ' \to \cC ' , \cQ ' , f':\cX' \to \bM _{1,1})$ over $S'$ 
such that:
\begin{enumerate}
\item
$\cX ' \to \cC '$ gives rise to an extension:
$$
\begin{array}{ccc}
{\mathcal X}_\eta\times_{\Delta} \Delta' &\subset & {\mathcal X}' \\
\down            &       & \down \\
{\mathcal C}_\eta \times_{\Delta} {\Delta}' &\subset & {\mathcal C}' \\
\down            &       & \down \\
\{\eta'\} & \subset & {\Delta}', \end{array}
$$

compatible with the extension $\cQ _{\eta} \times _S S' \subset \cQ ' ;$

\item $\cC '\to \bM _{1,1}$ is Konstevich stable.
\item the components of $\cX '_0$ that dominate the components of $\cC' _0$ are 
{\em log-standard stricly prestable} quadruples 
$(X,Q,G, f: X \to \bM _{1,1})$
where $G$ consists of either one or two fibres, which are either stable or twisted; if
 $X\to C$  in such a quadruple turns out to be isotrivial, then $C$ is not rational.

\item the components of $\cX '_0$ that are mapped to a point of $\cC' _0$ 
are {\em isotrivial log-pseudoelliptic surfaces} either of type I or II. 
In the former case 
they are attached 
to the rest of the central fibre according to lemma \ref{toric one fibre} 
and in the latter 
according to lemma \ref{toric two}.  
\end{enumerate}
 The extension is unique up to a unique 
isomorphism, and its
formation commutes with further finite extensions of discrete valuation
rings.
\end{theorem}

\begin{proof}

By the {\it strictly prestable reduction theorem} (theorem  \ref{divcont})
 we can find a finite base change 
$\Delta ' \to \Delta$ 
, $\Delta ' \text{-schemes} $ and $\Delta ' \text{-morphisms}$ 
unique up to unique isomorphisms 
$(\cX '\to \cC ' , \cQ ', f: \cX '\stackrel{\pi '}{\to} \cC ' 
\stackrel{j '}{\to} \bM _{1,1})$ that extend 
$({\cX}_\eta\to {\cC}_\eta ,\cQ _{\eta} , f_\eta ),$ which is {\it strictly 
prestable}. The extension commutes with further base changes and 
the log-canonical divisor 
$\omega _{\cX ' / \Delta '}(\cQ ') 
\otimes {f'}^* \cO _{\bM _{1,1}} (3)$ is ample away 
from those isotrivial components of the central fibre $\cX ' _0$ the meet 
the rest of the central fibre in one or two fibres (twisted or stable).

 Let $X\to C$ such a component. Then 
${j'}^* \cO _{\bM _{1,1}} \simeq \cO _C,$ 
and therefore 
$$\omega _{\cX ' / \Delta '} (\cQ ')
\otimes f^*\cO _{\bM _{1,1}} (3) \mid _X 
\simeq \omega _{\cX  / \Delta '}(\cQ ') \mid _X.$$

According to lemma 
\ref{extr ray}
we know that if an irreducible rational curve meets the rest of the central 
fibre in one 
point, we need to log flip the zero section above it. 
According to 
theorem \ref{log-flip} we can then perform the log flip to get 
$(\cX ^+ \to \cC ^+ , \cQ ^+ ).$ In doing so, one produces an 
isotrivial pseudoelliptic 
surface of type I $(X_1 ^+, G_1)$ attached to a log-standard elliptic surface 
$(X_2 ^+ ,Q _2 ^+, G_2^+)$ according to the fan in theorem \ref{log-flip}. According to 
the same theorem, the zero section $Q_2 ^+$ of $X_2 ^+$  
misses the singular point where it meets $X_1 ^+ ,$ so we can iterate 
the process, and 
prune the tree of all the $j \text{-trivial}$ rational curves meeting 
the rest of the central 
fibre in only one point. 

The log-canonical bundle is now nef and big
 after Proposition \ref{flipextremalray}, since the family $\cX \to \cC$ was 
{\it strictly prestable} to begin with (so it did not have any component 
whose zero section had self-intersection $\geq -1$), and according to 
proposition \ref{self-int2} it stays such. Furthermore, according to
theorem \ref{log-flip}, this way we only produce semi-logcanonical singularities.

 In order to make the log-canonical bundle {\it ample}, we need to contract 
all the chains of rational curves that meet the rest of the central 
fibre in two ends, 
according to lemma \ref{extr ray}. But this is taken care by theorem 
\ref{logcancont}, and will produce 
{\it isotrivial pseudoelliptic surfaces of type II}. The log-canonical
bundle is now ample on accounts of proposition \ref{ample2}. On
accounts of theorem \ref{logcancont}, the singularities we thus obtain
are at most semi-loganonical.
 This ends the proof.

\end{proof}

\subsection{Stable reduction for pairs} \label{stablepairs}

Here we want to deal with the absolute case. As it has been mensioned earlier, 
the steps of the MMP in a one parameter family are going to be similar to the 
ones performed in the case of triples, except that we now need to perform the 
flips and the small contractions also in cases in which the $j$-map is not 
constant.

\begin{theorem} \label{stabredpairs}
 Let ${\cX}_\eta\to {\cC}_\eta\to \eta$ be a stable 
elliptic surface over a smooth base curve 
$\cC _\eta$ of genus $g\geq 2 .$ Then there is a finite extension of
discrete valuation rings
$R\subset R'$ and a pair 
$( \cX ' \to \cC ' , \cQ ' )$ over $S'$ 
such that:
\begin{enumerate}
\item
$\cX ' \to \cC '$ gives rise to an extension:
$$
\begin{array}{ccc}
{\mathcal X}_\eta\times_SS' &\subset & {\mathcal X}' \\
\down            &       & \down \\
{\mathcal C}_\eta \times_S S' &\subset & {\mathcal C}' \\
\down            &       & \down \\
\{\eta'\} & \subset & S', \end{array}
$$

compatible with the extension $\cQ _{\eta} \times _S S' \subset \cQ ' ;$

\item the components of $\cX '_0$ that dominate the components of $\cC' _0$ are 
{\em log-standard strictly prestable } triples $(X,Q,G, )$
where $G$ consists of either one or two fibres, which are either stable or twisted.

\item the components of $\cX '_0$ that are mapped to a point of $\cC' _0$ 
are {\em log-pseudoelliptic surfaces} either of type I or II. In the former case 
they are attached 
to the rest of the central fibre according to lemma \ref{toric one fibre} 
and in the latter 
according to lemma \ref{toric two}.  
\end{enumerate}
 The extension is unique up to a unique 
isomorphism, and its
formation commutes with further finite extensions of discrete valuation
rings.

\end{theorem}

\begin{proof}

By the {\it strictly prestable reduction theorem} (theorem  \ref{divcont})
 we can find a finite base change 
$\Delta ' \to \Delta$ 
, $\Delta ' \text{-schemes} $ and $\Delta ' \text{-morphisms}$ 
unique up to unique isomorphisms 
$(\cX \to \cC , \cQ , f: \cX \stackrel{\pi}{\to} \cC 
\stackrel{j}{\to} \bM _{1,1})$ that extend 
$({\cX}_\eta\to {\cC}_\eta ,\cQ _{\eta} , f_\eta ),$ which {\it strictly 
prestable}. The extension commutes with further base changes.  

The log-canonical divisor 
$\omega _{\cX ' / \Delta '}(\cQ ')$ is ample away from those  
components fibred over a rational curve that meet the central fibre $\cX ' _0$
 along one or two fibres (stable or twisted). 

 In fact, let $r: Z \to B$ be a component of $\cX ' _0 .$ 
On the one hand, if $B$ is not rational, 
it is obvious that  $\omega _{\cX ' / \Delta '}(\cQ ')\otimes \cO _Z
\simeq \omega _X (\cQ ' \mid _Z ) \otimes \cO _Z (D),$ where $D$ is the 
dual curve, is ample, since $\omega _{Z/B} (\cQ ' \mid _Z)$ is relatively 
ample (Kollar semipositivity theorem) and $\omega _B$ is ample.
 On the other hand, if $B$ is rational but meets the rest 
of $\cC ' _0$ in at least three points, then $r^*\omega _B (D)$ is ample.
 The rest of the proof can be translated word by word 
from theorem \ref{stabred}.

\end{proof}

\begin{remark} It is worth noting that if the base curve $\cC _\eta$ is 
rational or elliptic, then the log-canonical bundle is not ample: one needs to contract 
all the base curves. In this sense it is probably more natural, 
in the rational and elliptic 
base curve case, to consider the moduli of triples (with map to $\bM _{1,1}$,)
 at least if one wants to preserve the fibration structure.

\end{remark}

\subsection{The rational  base case}

Here we deal with one of the two cases left out from section \ref{stablepairs}:
namely the case in which the base curve of the elliptic surface 
$\cX _{\eta} \to \cC _{\eta} ,$ is rational.
 In this case, the log-canonical bundle is not ample even on the surface 
$\cX _{\eta}$ itself: indeed we need to contract the zero 
section $\cQ _{\eta}, $ to make the log-canonical bundle ample.

We perform the stable reduction theorem for the triple 
$(\cX _{\eta} \to \cC _{\eta},
\cQ _{\eta}, j_{\eta} : \cC _{\eta}\to \bM _{1,1}).$ We may then assume to 
despose of a stable $\Delta$-triple $(\cX \to \cC \to \Delta, \cQ \to \Delta , 
j: \cC \to \bM _{1,1})$ over some DVR scheme $\Delta$ to begin with.
We want to be able to contract the base curve in the general member.

We can now state and prove:

\begin{theorem} \label{rationalbase}

Let $({\cX}_\eta\to {\cC}_\eta\to \eta , \cQ _{\eta}\to \eta ,
 {\cX}_\eta \to \bM _{1,1} )$ be a stable 
triple over a {\em rational} smooth base curve 
$\cC _\eta \simeq \PP ^1.$ Furthermore, let ${\cY}_\eta \to {\cX}_{\eta}$ be 
the log-canonical contraction of $\cQ _{\eta} .$

 Then there is a finite extension of 
discrete valuation rings
$R\subset R'$ and a pair 
$( \cX '  , \cQ ' )$ over $\Delta ' =Spec(R')$ 
such that:
\begin{enumerate}

\item
$\cX '$ gives rise to an extension:
$$
\begin{array}{ccc}
{\mathcal X}_\eta\times_{\Delta } \Delta ' &\subset & {\mathcal X}' \\
\down            &       & \down \\

\{\eta'\} & \subset & \Delta ', 
\end{array}
$$

compatible with the extension $\cQ _{\eta} \times _{\Delta } \Delta '
 \subset \cQ ' ;$

\item the components of $\cX '_0$ 
are {\em n-pseudoelliptic surfaces} with either $n=0$ or $n=1$.  
In the former case 
they are attached 
to the rest of the central fibre according to lemma \ref{toric one fibre} 
and in the latter 
according to lemma \ref{toric two}.  
\end{enumerate}
 The extension is unique up to a unique 
isomorphism, and its
formation commutes with further finite extensions of discrete valuation
rings.

\end{theorem}

\begin{proof}

We can apply Theorem \ref{stabred} to obtain a stable {\it triple} 
$({\cY}\to {\cC} \to \Delta , \cQ \to \Delta ,
 {\cX} \to \bM _{1,1} ), $ after a possible base change (unique up to 
a unique isomorphism) $\Delta ' \to \Delta,$ that satisfies the analogous 
of property ~1.

Either $\cC _0$ consists of a simple tree (i.e., a chain of rational
curves meeting transversally each only one consecutive curve at one
point), or every {\it chain-like} component $C
\subset \cC _0$ (i.e., a component that consists of a chain) will meet
another {\it tree-like component} in a leaf
which is not an extremity (i.e., this leaf will meet two more leaves
of tree-like component it belongs to). 

In the first case there are two
possibilities: the tree (i.e., $\cC _0$ has either an even or an odd
number of leaves. If it has an odd number of components, there is a
well-defined {\it central leaf}, that is to say the leaf which
disconnets the tree in two components of equal lenght. In case there
is an even number of components, then there is, analogously, a well
defined concept of {\it central pair of leaves}.

So in the case of a simple tree of odd lenght, let $B$ be the central
curve. We can prune, starting
from the two ends (by
log-flipping according to theorem \ref{log-flip}) all the leaves that belong
to the two chain-like components that $B$ disconnets from $\cC _0 .$
Let as call $\cY' \to \cC ' \to \Delta$ the family thus obtained from
$\cY \to \cC  \to \Delta ,$ and $\cQ '$ the new zero-section. 
At this point the zero section $\cQ ' \mid _B$ of the log-standard
surface $Y:= \cY ' \mid _B \to B$ has the property that $L _Y \cdot Q
=-2, $ which is the same as in the general fibre. We can now
divisorially contract the zero section $\cQ '$ in $\cY$  (e.g., by
means of the line bundle $\omega _{\cY '} ((1+a) \cQ ')$ where
$a=\frac {2}{{\cQ '_{\eta} }^2}$). The log-canonical bundle is now
ample, according to proposition \ref{flipextremalray}. Also, on
accounts of theorems \ref{log-flip} and \ref{logcancont} the
singularities thus produced are at most semi-logcanonical. 

The result in the special fibre is
a configuration of two chains of $0$-pseudoelliptic surfaces attached
to another $0$-pseudoelliptic surface (the surface obtained from $Y$
by contracting the zero-section).

Analogously, if we a simple tree of even lenght, do the same
operations word by word as above, barring that now $B$ is replaced by
the {\it central pair} $B_1 \cup B_2$. 

The result, now consists of two chains of $0$-pseudoelliptic surfaces
as above, only now attached to a union of two $1$-pseudoelliptic
surfaces (the result of contracting the zero sections of $\cY ' \mid
_{B_1} \cup \cY ' \mid _{B_2}$).

In the event that there are two chain-like components , we can now
individuate a {\it spine}, namely the one component that is attached to
the other at an extremity. In this case we first prune (by means of
theorem \ref{log-flip}) the other component and reduce ourselves to
considering only the spine, in other words reducing the problem to the
previous case.

We can now conclude by a simple induction argument.
\end{proof}
 
\subsection{The elliptic  base case}

Here, at last, we deal with the final case.

We first need the following:

\begin{definition}

A surface $Y$ with a structure morphism $g: Y'\to Y$ is said to be a
{\bf pseudoelliptic surface of type $E_0$} (resp. {\bf $E_{I_N}$}) if
the surface $(Y' \to E,Q,F) $ is a log-standard elliptic surface, with
one marked fibre $F$ and a zero section $Q,$ mapping to an
irreducible elliptic curve $E$ (resp. to a closed chain of rational
curves $E$) as base curve, and if $g$ has the zero section $Q$ as
exceptional curve. The singularity of $Y$ at $g(Q)$ is an {\it elliptic}
(resp. {\it degenerate cusp}) singularity  

\end{definition}
We have: 

\begin{theorem} \label{ellipticbase}

Let $({\cX}_\eta\to {\cC}_\eta\to \eta , \cQ _{\eta}\to \eta ,
{\cX}_\eta \to \bM _{1,1} )$ be a stable 
triple over an {\em elliptic} smooth base curve 
$\cC _\eta \simeq E.$ Furthermore, let ${\cY}_\eta \to {\cX}_{\eta}$ be 
the {\it pseudoelliptic surface of type $E_0$ } obtained by contracting $\cQ _{\eta} .$

 Then there is a finite extension of 
discrete valuation rings
$R\subset R'$ and a pair 
$( \cX '  , \cQ ' )$ over $\Delta ' =Spec(R')$ 
such that:
\begin{enumerate}

\item
$\cX '$ gives rise to an extension:
$$
\begin{array}{ccc}
{\mathcal X}_\eta\times_{\Delta } \Delta ' &\subset & {\mathcal X}' \\
\down            &       & \down \\

\{\eta'\} & \subset & \Delta ', 
\end{array}
$$

compatible with the extension $\cQ _{\eta} \times _{\Delta } \Delta '
 \subset \cQ ' ;$

\item the components of the central fibre $\cX '_0$ consist of a {\it pseudoelliptic
    surface} either of {\it type $E_0$} or of type {\it $E_{I_N}$},
  attached to a configuration of $0$-pseudolliptic and
  $1$-pseudoelliptic surfaces as in theorem \ref{rationalbase}.

\end{enumerate}
 The extension is unique up to a unique 
isomorphism, and its
formation commutes with further finite extensions of discrete valuation
rings.

\end{theorem}

\begin{proof}

We may assume, after theroem \ref{rationalbase}, that we have a family
$\cY \to \cC ' \to \Delta$ whose generic fibre is isomorphic to $\cX
_{\eta} \to {\cC '} _{\eta},$ and whose central fibre is a a surface $\cY
_0$ whose components are a log-standard elliptic surface $(Y \to E ,
Q+F)$ with one marked fibre $F$ and having as a base curve $E$ either an elliptic curve or a
closed chain of rational curves, attached along $F$ to a configuration
of $0$-pseudolliptic and  $1$-pseudoelliptic surfaces as in theorem \ref{rationalbase}.

We can now divisorially contract the zero section $\cQ '$ of $\cY \to
\cC'$ to obtain the result stated.

\end{proof}


\begin{thebibliography}{MMMM}
\bibitem[$\aleph$-V1]{A-V2} D. Abramovich and A. Vistoli, 
{\em Compactifying the space of stable maps}, preprint 1999.

\bibitem[$\aleph$-V2]{A-V} D. Abramovich and A. Vistoli, {\em
Complete moduli for fibered surfaces}, preprint 1997.





\bibitem[A1]{A1}
V.Alexeev,{\em Moduli spaces $M \sb {g,n}(W)$ for surfaces,} in 
{\it Higher-dimensional complex varieties (Trento, 1994)}, 1--2, 
de Gruyter, Berlin.

\bibitem[A2]{A2}
V.Alexeev, {\em Boundedness and $K^2$ for log surfaces,} Int.J.Math. 
{\em 5} (1994), no. 6.

\bibitem[{\em SGA} 1]{SGA1} A.Grothendieck,
 {\em `` Revetements etales et Groupe Fondamental'', Lectures Notes in Math. 224
,}, Springer 
Verlag, New York, {\bf 1971}

\bibitem[B-P-VdV]{B} W.Barth, C.Peters, A.Van de Ven 
{\em Compact Complex Surfaces} Springer-Verlag 1984

\bibitem [D-M]{D-M} P. Deligne, D. Mumford, {\em The irreducibility of
    the moduli space of curves of given genus}, Publication I.H.E.S
  {\bf 36} (1969).
\bibitem[dJ-O]{dj-o} A. J. de Jong and F. Oort, {\em On extending families of
curves}, J. Alg. Geom. {\bf 6} no. 3 (1997), 545-562.



\bibitem[F]{F} W.Fulton, {\em Introduction to Toric Varieties,}
 Annals of Mathematics Studies N.$131$, Princeton University Press 

\bibitem [H-M] {H-M} J. Harris, I. Morrison, {\em Moduli of curves},
 GTM {\bf 187}, Springer-Verlag
\bibitem[H] {H} 
R. Hartshorne, {\em Algebraic Geometry}, GTM {\bf 52}, Springer-Verlag (1977)

\bibitem[Ha] {Ha}
B. Hasset, {\em Stable log surfaces and limits of quartic plane curves,} 
   Manuscripta Mathematica 100 (1999), 469-497

\bibitem[KMM]{KMM}

Y. Kawamata; K. Matsuda; K. Matsuki, {\em Introduction to the minimal
  model program}, Adv. Stud. Pure Mat. 10, 283-360 (1987)

\bibitem[K-S]{K-S}
J.Koll\'ar and N. Shepherd-Barron, {\em Threefolds and deformations of 
surface singularities,} Invent. Math. {\em 91} (1988), 299-338.



\bibitem[K]{K}
J. Koll\'ar et al., {\em Flips and abundance for algebraic threefolds},
 Asterisque, vol. 211, 1992.


\bibitem[M1]{M}
R. Miranda, {\em The moduli of Weiertstrass Fibrations Over $\PP ^1$}, 
Math. Ann. {\em 255} (1981), 379-394

\bibitem[M2]{M2}
R. Miranda, {\em The basic theory of elliptic surfaces } 
Universit\`a di Pisa, Dipartimento di Matematica, Pisa 1982

\bibitem[Mo]{Mo} 
S. Mori, {\em Threefolds whose canonical bundles are not
  numerically effective} Ann. Math. {\em 116} (1982), 133-176

\bibitem[M-S]{Mu} 
D. Mumford and K. Suominen, {\em Introduction to the theory of moduli }, 
5-th Nordic Summer School in Mathematics, Oslo, 1970

\bibitem[R1]{R1}
M. Reid, {\em Decomposition of Toric Morphisms}, volume in honor of Shafarevich.

\bibitem[R2]{R2}
M. Reid, {\em Minimal Models of Canonical Threefolds}, Advanced Studies in Pure 
Mathematics 1, edited by S.Iitaka.

\bibitem[R3]{R3} 
M. Reid, {\em Young person's guide to canonical singularities }, 
Algebraic Geometry Bowdoin 1985, Proceedings of Symposia in Pure Mathematics


\bibitem[S]{S} 
J. Silvermann, {\em Advanced Topics in the Arithmetic of Elliptic Curves},
Springer-Verlag




\end{thebibliography}
\end{document}